\documentclass[10pt]{article}

\DeclareMathAlphabet{\mathpzc}{OT1}{pzc}{m}{it}

\usepackage{amsmath}
\usepackage{amssymb}
\usepackage{amsfonts}
\usepackage{amsthm}
\usepackage{mathrsfs}
\usepackage{enumerate}
\usepackage[pdftex]{color, graphicx}
\usepackage{hyperref}
\newtheorem{thm}{Theorem}[section]
\newtheorem{lem}[thm]{Lemma}

\newtheorem{prop}[thm]{Proposition}
\newtheorem{cor}[thm]{Corollary}
\newtheorem{conj}{Conjecture}
\numberwithin{equation}{section}

\begin{document}

\title{{\bf On Ramanujan's cubic continued fraction, II}}

\author{Sushmanth J. Akkarapakam and Patrick Morton}
\maketitle

\begin{abstract}
We show that all the periodic points with period $n \ge 3$ of the algebraic function defined by the relation $g(x,y) = x^2+2xy^2-y = 0$ are expressible in terms of values of Ramanujan's
cubic continued fraction $c(\tau)$.  It is shown that the values $c(w/2), c(w/6), c(w/6+1/2)$ are periodic points, for specific elements $w$ in the order of discriminant $-d$ in 
$K = \mathbb{Q}(\sqrt{-d})$, with $-d \equiv 1$ (mod $8$); and all the periodic points are conjugates of these values or specific linear fractional expressions in their conjugates, times cube roots of unity.  These periodic points are determined by relating them to the periodic points of the algebraic function defined by $k(x,y)= -xy^2+x^2+y^2+3x = 0$; and these periodic points (for periods $\ge 2$) generate ring class fields over $K$.
\end{abstract}

\section{Introduction.}

This paper is part of a larger project to investigate the arithmetic properties of the periodic points of algebraic functions.  As in previous papers we focus on
several algebraic functions which arise from the relations between specific modular functions.  \medskip

In \cite{akm} we investigated properties of Ramanujan's cubic continued fraction
\begin{align}
c(\tau)=\cfrac{q^{1/3}}{1+\cfrac{q+q^2}{1+\cfrac{q^2+q^4}{1+\cfrac{q^3+q^6}{1+\cdots}}}}, \ \ q = e^{2\pi i \tau},
\label{eqn:1.1}
\end{align}
which has the product representation
\begin{equation}
c(\tau)=q^{1/3}\prod_{n=1}^\infty\frac{(1-q^{6n-1})(1-q^{6n-5})}{(1-q^{6n-3})^2} = \frac{\eta(\tau) \eta^3(6\tau)}{\eta(2\tau) \eta^3(3\tau)}
\label{eqn:1.2}
\end{equation}
in terms of the Dedekind $\eta$-function.  We focused on the algebraic function $y = \mathfrak{f}(z)$ defined by the polynomial
$$f(x,y) = x^3(y^2+y+1)-y(y^2-2y+4);$$
this gives the modular identity $f(2c(\tau),2c(3\tau)) = 0$ between $c(\tau)$ and $c(3\tau)$.  We determined all the periodic points of $\mathfrak{f}(z)$
in terms of values of $c(\tau)$, one of which is the value $2c(w/3)$, where
$$w = \frac{k+\sqrt{-d}}{2}, \ \ -d \equiv 1 \ (\textrm{mod} \ 3), \ \ k^2 + d \equiv 0 \ (\textrm{mod} \ 9)\footnote{$k$ must also satisfy $k \equiv 2$ (mod $12$) or $k \equiv 1$ (mod $6$), according as $-d$ is even or odd.}.$$
We also showed that $\mathbb{Q}(c(w/3))$ equals the ring class field $\Omega_{2f}$ of ring conductor $2f$
over $K = \mathbb{Q}(\sqrt{-d})$, where $-d = d_K f^2$ and $d_K$ is the discriminant of $K$.  See the more detailed statement of results in Theorem \ref{thm:7.1} in the Appendix (Section 7). \medskip

In this paper we consider the identity $g(c(\tau),c(2\tau)) = 0$ and the algebraic function $y = \mathfrak{g}(x)$ defined by $g(x,y) = 0$, where
$$g(x,y) = x^2+2xy^2-y, \ \ \mathfrak{g}(x) = \frac{1 \pm \sqrt{1-x^3}}{x}.$$
We will determine all the periodic points of $\mathfrak{g}(x)$ in terms of values of $c(\tau)$, one of which is $c(w/2)$; and will show that $\mathbb{Q}(c(w/2))$ equals the ring class field $\Omega_{3f}$ of conductor $3f$ over $K$, where again $-d = d_K f^2$ and
\begin{equation}
w = \frac{a+\sqrt{-d}}{2}, \ \ -d \equiv 1 \ (\textrm{mod} \ 8), \ \ a^2 + d \equiv 0 \ (\textrm{mod} \ 8).
\label{eqn:1.3}
\end{equation}
This investigation turns out to be much more difficult than the investigation into the periodic points for $f(x,y) = 0$, in spite of the fact that $g(x,y)$ is only 
quadratic in $y$.  It requires us to consider the auxiliary functions
$$s(\tau) = \frac{\eta^4(\tau) \eta^2(6\tau)}{\eta^2(2\tau) \eta^4(3\tau)}, \ \ r(\tau) = \frac{\eta^2(\tau) \eta^4(6\tau)}{\eta^4(2\tau) \eta^2(3\tau)},$$
which satisfy the relations
$$s^2(\tau) = 1-8c^3(\tau), \ \ r^2(\tau) = \frac{1-s^2(\tau)}{9-s^2(\tau)}.$$
In order to determine the periodic points for $\mathfrak{g}(z)$, we first determine all the periodic points for the algebraic function defined by
$$k(x,y) = -xy^2+x^2+y^2+3x = 0,$$
where $k(-s(\tau),s(2\tau)) = 0$ is the modular relationship between $s(\tau)$ and $s(2\tau)$.  This turns out to be interesting in its own right, and gives us
the opportunity to use some classical tools from the theory of complex multiplication (see \cite{deu2}). \medskip

\begin{thm} \begin{enumerate}[(a)]
\item If $w$ is given by \eqref{eqn:1.3}, then $s(w/2)$ and $\pm 3r(w/2)$ (for a specific sign $\pm$) are periodic points for the algebraic function defined by
$k(x,y) = 0$.  Their minimal periods with respect to this algebraic function are always even.

\item If $w$ is given by \eqref{eqn:1.3} and, additionally, $3 \mid N(w)$ and $\left(\frac{d_K}{3}\right) = 0, 1$, then $s(w/6), s(w/6 + 1/2)$ are periodic points for $k(x,y) = 0$.

\item The only periodic points in $\overline{\mathbb{Q}}$ for the algebraic function defined by $k(x,y) = 0$ are $0, -1, 3$, and the values and conjugates over 
$\mathbb{Q}$ of the values in (a) and (b), as $-d$ runs over all negative quadratic discriminants $-d \equiv 1$ mod $8$. 

\item For the discriminant  $-d = d_K f^2$, $\mathbb{Q}(s(w/2)) = \mathbb{Q}(r(w/2)) = \Omega_{3f}$ and $\mathbb{Q}(s(w/6)) = \mathbb{Q}(s(w/6+1/2)) = \Omega_f$.  The sign $\varepsilon = \pm 1$ of $3\varepsilon r(w/2)$ in (a) is given by
$$\varepsilon = -r(w/2) s\left(\frac{w}{2}+\frac{1}{2}\right)^{\bar \tau^{-1}}, \ \ \bar \tau = \left(\frac{\Omega_{3f}/K}{\mathfrak{p}_2}\right),$$
where $(2) = \mathfrak{p}_2 \mathfrak{p}_2'$ and $\mathfrak{p}_2 \mid w$.
\end{enumerate}
\label{thm:1.1}
 \end{thm}

We will show that the value $s(w/2)$ is a periodic point for $k(x,y) = 0$ by showing that $k(s(w/2), s(w/2)^{\bar \tau}) = 0$ (Corollary \ref{cor:3.13}).
The function $s(\tau)$ is an example of a modular function satisfying $k(-s(\tau),s(2\tau)) = 0$, but for which the values $s(w/2)$ are periodic points not for
$k(-x,y) = 0$, but for $k(x,y) = 0$.  
\medskip

A major role in our discussion is played by the identities
\begin{align}
\label{eqn:1.4} u^8(\tau) &= \frac{(s(\tau)-1)(s(\tau)+3)^3}{(s(\tau)+1)(s(\tau)-3)^3},\\
\label{eqn:1.5} u^8(3\tau) &= \frac{(s(\tau)+3)(s(\tau)-1)^3}{(s(\tau)-3)(s(\tau)+1)^3},
\end{align}
where $u(\tau)$ is the modular function
$$u(\tau) = \sqrt{2}q^{1/8} \prod_{n = 1}^\infty{\left(\frac{1+q^{2n}}{1+q^{2n-1}}\right)}, \ \ q = e^{2\pi i \tau},$$
equal to the continued fraction
$$u(\tau) = \sqrt{2} \ \frac{q^{1/8}}{1+} \ \frac{q}{1+q+} \ \frac{q^2}{1+q^2+} \ \frac{q^3}{1+q^3+}\dots, \ \ q = e^{2\pi i \tau},$$
for $\tau \in \mathcal{H}$ (the upper half-plane).  See \cite[p. 153]{duk}.  The above relations represent an indirect connection between Ramanujan's cubic 
continued fraction $c(\tau)$ and $u(\tau)$.  Since
$$s(\tau) = \frac{2c^2(\tau)}{c(2\tau)}-1,\ \ \tau \in \mathcal{H},$$
(Corollary \ref{cor:2.6}) these identities imply the relationships
\begin{align*}
u^8(\tau) &= \frac{(c^2(\tau)-c(2\tau))(c^2(\tau)+c(2\tau))^3}{c^2(\tau)(c^2(\tau)-2c(2\tau))^3},\\
u^2(\tau) u^2(3\tau) &= \frac{c^4(\tau)-c^2(2\tau)}{c^2(\tau)(c^2(\tau)-2c(2\tau))},
\end{align*}
between these two continued fractions (see Corollary \ref{cor:4.5}).  Equations \eqref{eqn:2.2} and \eqref{eqn:3.4} below then yield that
\begin{align*}
r(\tau) &= \frac{c^2(\tau)}{c(\tau/2)+c^2(\tau)} = \frac{u(\tau)u(3\tau)c^2(\tau)}{c^2(\tau)+c(2\tau)}, \\
u(\tau) u(3\tau) &= \frac{c^2(\tau)+c(2\tau)}{c^2(\tau)+c(\tau/2)}, \ \ \tau \in \mathcal{H}.
\end{align*}

The identity \eqref{eqn:1.4} leads us to consider the polynomial
$$\mathcal{Q}(x,y) = (1-x^2)^2(y - 1)(y + 3)^3 - (1+ x^2)^2(y + 1)(y - 3)^3$$
and the related resultant
\begin{equation}
\textrm{Res}_x( b_d(x), \mathcal{Q}(x,y)) = \mathcal{Q}_d(y),
\label{eqn:1.6}
\end{equation}
where $b_d(x)$, for $-d \equiv 1$ (mod $8$), is the minimal polynomial over $\mathbb{Q}$ of solutions $\pi, \xi$ of the quartic 
Fermat equation in the ring class field $\Omega_f$ over $K = \mathbb{Q}(\sqrt{-d})$ (see \cite{lm}, \cite{mdeu}), where $-d = d_K f^2$:
 $$\pi^4 + \xi^4 = 1, \ \ \xi = \frac{\pi^{\varphi^2}+1}{\pi^{\varphi^2}-1}, \ \ \varphi = \left(\frac{\Omega_f/K}{\wp_2}\right).$$
 (It can be shown that these conditions determine the minimal polynomials $b_d(x)$ uniquely.)
  In addition, $(2) = \mathfrak{p}_2 \mathfrak{p}_2'$ in $K$, with $(\pi) = \mathfrak{p}_2, (\xi) = \mathfrak{p}_2'$ in $\Omega_f$.  
 See \cite[Thm. 1.1, Section 10]{lm} and \cite{mdeu}.
These solutions are themselves periodic points of the algebraic function defined by
$$G(x,y) = (y^2+1) x^2 +2y = 0, \ \ y = \frac{-1\pm \sqrt{1-x^4}}{x^2}.$$

 \begin{thm}
Let $\mathcal{S}$ denote the set of periodic points for $k(x,y) = 0$ in $\overline{\mathbb{Q}}$.  Then  $\mathcal{S} \setminus \{0,-1,3\}$ 
coincides with the collection of all roots of $\mathcal{Q}_d(y)$, as $-d$ ranges over all negative quadratic discriminants $-d \equiv 1$ (mod $8$).
 The three excluded values are the roots of $\textrm{Res}_x(\mathcal{Q}(x,y),x(x+1))$.  The polynomial $\mathcal{Q}_d(y)$ factors into $3, 2$ or $1$ 
 distinct irreducible polynomials, according as $\left(\frac{-d}{3}\right) = 1, 0$ or $-1$.  The periodic point $s(w/2)$ is always a root of the irreducible factor 
 of $\mathcal{Q}_d(y)$ of largest degree. 
 \label{thm:1.2}
 \end{thm}
 
The polynomial $\mathcal{Q}_d(y)$ plays the same role here that the polynomial $C_d(x)$ plays in \cite[Section 4]{akm}.  
 The identities for $u(\tau)^8$ and $u(3\tau)^8$ also imply the following theorem.
 
 \begin{thm}
If $\left(\frac{d_K}{3}\right) = 0$ or $1$ and $w = \frac{a+\sqrt{-d}}{2}$, with $a^2 + d \equiv 0$ mod $24$, where $-d = d_Kf^2, 3 \nmid f$, 
and $\frak{p}_2 \mid w$, then 
$$\upsilon = \frac{(s(w/6)-1)(s(w/6)-3)}{(s(w/6)+1)(s(w/6)+3)}$$
is the fourth power of a unit in the ring class field $\Omega_f$ over $K = \mathbb{Q}(\sqrt{-d})$.
\label{thm:1.3}
\end{thm}
 
 For the algebraic function defined by $g(x,y) = 0$, we set
 $$g_1(x,y) = 4g(x/2,y/2) = x^2+xy^2-2y$$
 and determine the periodic points for $g_1(x,y) = 0$.  Let $\omega = (-1+\sqrt{-3})/2$.
 
 \begin{thm}
The only periodic points for the algebraic function defined by $g_1(x,y) = 0$ are $0,\omega^i,-2\omega^i$ for $i \in \{0,1,2\}$, and the values and conjugates
over $\mathbb{Q}$ of
\begin{equation}
2\omega^i c(w/2), \ 2\omega^i c(w/6), \ 2\omega^i c(w/6+1/2), \ 2\omega^ic(-1/3w),
\label{eqn:1.7}
\end{equation}
where
$$w = \frac{a+\sqrt{-d}}{2}, \ \ a^2+d \equiv 0  \ \textrm{mod} \ 2^3,$$
and $-d = d_Kf^2$ runs over all negative quadratic discriminants satisfying $-d \equiv 1$ (mod $8$).
The last periodic point can also be written as
$$2\omega^ic(-1/3w) = \omega^i \frac{1-2c(3w/2)}{1+c(3w/2)} = \omega^i \mathcal{A}(2c(3w/2)) = \omega^i B\left(\frac{-1}{c(3w/2)}\right),$$
where
$$\mathcal{A}(x) = \frac{2-2x}{2+x}, \ \ B(x) = \frac{x+2}{x-1};$$
and $\frac{-1}{c(3w/2)}$ and $2c(3w/2)$ are periodic points for $g_1(x,y)=0$.  All of the periodic points for $g_1(x,y) = 0$ are algebraic integers; and the 
conjugates of the first three listed values in \eqref{eqn:1.7} lie in the ring class field $\Omega_{3f} = \mathbb{Q}(c(w/2))$ over $K = \mathbb{Q}(\sqrt{-d})$, 
while the last value lies in $\mathbb{Q}(c(3w/2)) = \Omega_{9f}$.
\label{thm:1.4}
\end{thm}

Like the results of \cite{akm}, this theorem supplements the result of \cite[Theorem 13]{ckp} since, for example,  $c(w/2)$ is evaluated at a basis quotient for the ideal
$\mathfrak{p}_2^{(-d)} = (2, w) \cap \textsf{R}_{-d}$, where $\textsf{R}_{-d}$ is the order of discriminant $-d = d_k f^2$, and the result of \cite{ckp} is only given for the maximal order and ideals that are relatively prime to $6$.  In particular, the values $c(w/2)$ considered here do always generate the ray class field $\Sigma_6$
of conductor $6$ over $K$, when $f = 1$, since $\Sigma_6 = \Omega_3$ when $d_K \equiv 1$ (mod $8$).  On the other hand, the values 
$\omega^i c(w/6) = c(w/6+i)$ in Theorem \ref{thm:1.4} sometimes only generate the Hilbert class field $\Sigma$ over $K$, when $f = 1$ (in particular, when $-d \equiv 1$ mod $3$; see Proposition \ref{prop:5.7}).  Thus, the conclusion of Theorem 13 of Cho-Koo-Park \cite{ckp} does not always hold when $c(\tau)$ is evaluated at a basis quotient for an ideal which is not relatively prime to $6$. \medskip
 
The paper is laid out as follows.  In Section 2 we prove identities for $s(\tau), r(\tau)$ and $c(\tau)$ that are fundamental for the whole paper, including the fact that 
 $k(-s(\tau), s(2\tau)) = 0$.  In Section 3 we prove that $s(w/2)$ is a periodic point for $k(x,y) = 0$ using \eqref{eqn:1.4} and facts from \cite{akm0} and \cite{mdeu}.
 In Section 4 we determine the other periodic points for $k(x,y)$ mentioned in Theorem \ref{thm:1.1}({\it a,b}).  Then we use this to prove Theorem \ref{thm:1.1} in Section 6.  In Section 5 we prove that the values in \eqref{eqn:1.7} are periodic points for $g_1(x,y) = 0$.  Then we derive the rest of Theorem \ref{thm:1.4} from 
 Theorem \ref{thm:1.1}(c) in Section 6.  Minimal polynomials of periodic points of small periods for $k(x,y) = 0$ and $g_1(x,y) = 0$ are given in Tables \ref{tab:2}-\ref{tab:5}  at the end of Section 6.  In the Appendix (Section 7), we state the main theorem that we proved in \cite{akm} and fill in some small details that we neglected 
 to mention in \cite{akm}.  We also fill in a gap in the proof of \cite[Prop. 18(a), p. 1751]{akm} in Lemma \ref{lem:7.2}.
 
\section{Identities for modular functions.}

\begin{thm}
The modular functions
\begin{equation}
\alpha(\tau) = \frac{1}{c(\tau)}, \ \ \beta(\tau) = \frac{c(\tau/2)+c^2(\tau)}{c^2(\tau)} = \frac{-1}{c_1(\tau/2)c(\tau)}+1
\label{eqn:2.1}
\end{equation}
parametrize the elliptic curve
$$\mathcal{E}: \ \beta^2 = \alpha^3+1.$$
\label{thm:2.1}
\end{thm}

\begin{proof}
Putting $\tau/2$ for $\tau$ in the identity connecting $c(\tau)$ and $c(2\tau)$ from \cite[Prop. 1]{akm} gives
$$c^2(\tau/2)+2c^2(\tau)c(\tau/2)-c(\tau) = 0,$$
and completing the square yields
$$(c(\tau/2)+c^2(\tau))^2= c^4(\tau)+ c(\tau).$$
Hence,
\begin{align*}
\beta^2(\tau) &= \left(\frac{c(\tau/2)+c^2(\tau)}{c^2(\tau)}\right)^2 = \frac{c^4(\tau)+c(\tau)}{c^4(\tau)}\\
&=\frac{1}{c^3(\tau)} + 1 = \alpha^3(\tau) + 1.
\end{align*}
This proves the assertion.  For this note that $c(\tau) = -c(\tau/2)c_1(\tau/2)$ from \cite[Eq. (4.2)]{akm}.
\end{proof}

From Cooper \cite[pp. 368-371]{coo}, we take the identities
$$\frac{1}{r_b(q)^2}-\frac{1}{r_c(q)^3}=1,$$
where
\begin{align}
\label{eqn:2.1a} r_c(q) &= \frac{\eta(\tau) \eta^3(6\tau)}{\eta(2\tau) \eta^3(3\tau)} = c(\tau),\\
\label{eqn:2.1b} r_b(q) &= \frac{\eta^2(\tau) \eta^4(6\tau)}{\eta^4(2\tau) \eta^2(3\tau)} = r(\tau),
\end{align}
where the first identity is from \cite[Thm. 6.6, p. 367]{coo} or \cite[p. 1714]{akm}.  The above identity gives that
\begin{equation}
\frac{1}{r(\tau)^2} = \frac{1}{c^3(\tau)}+1.
\label{eqn:2.1c}
\end{equation}
Hence, Theorem \ref{thm:2.1} gives that
\begin{align}
\notag \frac{1}{r(\tau)} &= \frac{c(\tau/2)+c^2(\tau)}{c^2(\tau)},\\
\label{eqn:2.2} r(\tau) &= \frac{c^2(\tau)}{c(\tau/2)+c^2(\tau)},
\end{align}
since the leading term of both sides is $q^{1/2}$. \medskip

Now define
\begin{equation}
s(\tau) = \frac{\eta^4(\tau) \eta^2(6\tau)}{\eta^2(2\tau) \eta^4(3\tau)},
\label{eqn:2.3}
\end{equation}
so that \cite[Eq. (6.15), p. 370]{coo}, gives that
\begin{equation}
1-8c^3(\tau) = s^2(\tau).
\label{eqn:2.4}
\end{equation}
Writing
$$s(\tau) = \left(\frac{\eta(\tau)}{\eta(3\tau)}\right)^4 \left(\frac{\eta(6\tau)}{\eta(\tau)}\right)^2 \left(\frac{\eta(\tau)}{\eta(2\tau)}\right)^2$$
and using Newman \cite[Thm. 1, p. 375]{mnew} shows that $s(\tau)$ is a modular function for $\Gamma_0(12)$.

\begin{lem}
We have the identities
\begin{equation}
s\left(\frac{-1}{\tau}\right) = 3r\left(\frac{\tau}{6}\right), \ \ s\left(\frac{\tau}{6}\right) = 3r\left(\frac{-1}{\tau}\right).
\label{eqn:2.5}
\end{equation}
\label{lem:2.2}
\end{lem}

\begin{proof}
This follows easily from the identity $\eta(-1/\tau) = \sqrt{-i \tau} \eta(\tau)$:
\begin{align*}
s\left(\frac{-1}{\tau}\right) &= \frac{\eta^4(-1/\tau) \eta^2(-6/\tau)}{\eta^2(-2/\tau) \eta^4(-3/\tau)}\\
&= \frac{(-i \tau)^2 \eta^4(\tau) (-i \tau/6) \eta^2(\tau/6)}{(-i \tau/2) \eta^2(\tau/2) (-i \tau/3)^2 \eta^4(\tau/3)}\\
&= \frac{3\eta^4(\tau) \eta^2(\tau/6)}{\eta^2(\tau/2) \eta^4(\tau/3)}\\
&= 3r(\tau/6).
\end{align*}
The second identity follows by putting $-6/\tau$ in for $\tau$ in the first.
\end{proof}

\begin{thm}
We have the identity
\begin{equation}
s(\tau) s^2(2\tau) +s^2(\tau)+s^2(2\tau)-3s(\tau) = 0.
\label{eqn:2.6}
\end{equation}
\label{thm:2.3}
\end{thm}

\begin{proof}
From \cite[Prop.1]{akm} we know that $u=c(\tau)$ and $v=c(2\tau)$ satisfy the equation $2uv^2=v-u^2$. On cubing both sides we get
\[8u^3v^6=v^3-u^6-3vu^2(v-u^2).\]
Then substituting $v-u^2=2uv^2$ on the right hand side and simplifying we get
\[8u^3v^6=v^3-u^6-6u^3v^3.\]
Let $x = s(\tau), y = s(2\tau)$.  Since $s^2(\tau)=1-8c^3(\tau)$, we have $u^3=\frac{1}{8}(1-x^2)$ and $v^3=\frac{1}{8}(1-y^2)$. Substituting these into the above equation and clearing the denominators gives us
$$(1-x^2)(1-y^2)^2=8(1-y^2)-(1-x^2)^2-6(1-x^2)(1-y^2),$$
which simplifies to
$$x^4+y^4+8x^2y^2=x^2y^4+9x^2,$$
or
$$(x^2+y^2)^2=x^2y^4+9x^2-6x^2y^2.$$
Hence
$$(x^2+y^2)^2=(xy^2-3x)^2,$$
giving us that either $x^2+y^2-xy^2+3x=0$ or $x^2+y^2+xy^2-3x=0$. But from the infinite products for $x$ and $y$, it is clear that $x, y$ are both $1+O(q)$ as $q$ tends to $0$. Hence, it is the second expression above (and not the first) that vanishes for $q$ sufficiently small. By the identity theorem of complex analysis, it vanishes for $\tau$ in the upper half-plane.
\end{proof}

\begin{cor}
The functions $x = r(\tau), y = r(2\tau)$ satisfy the identity
$$3x^2y+x^2+y^2-y = 0,$$
and $x = 3r(\tau), y = 3r(2\tau)$ satisfy the equation
$$x^2y+x^2+y^2-3y = 0.$$
\label{cor:2.4}
\end{cor}

\begin{proof}
These formulas follow from putting $-1/\tau$ for $\tau$ in Theorem \ref{thm:2.3}, using Lemma \ref{lem:2.2} to express everything in terms of $r(\tau)$, 
and then replacing $\tau$ by $12\tau$.
\end{proof}

We set
$$t(x) =  \frac{-x+3}{x+1}$$
and note the identity
\begin{align}
\label{eqn:2.9a} &k(x,y) = -xy^2+x^2+y^2+3x\\
\label{eqn:2.9b} &(x+1)^2(y+1)^2 k(t(x),t(y)) = 16k(y,x),
\end{align}
which is easy to verify.

\begin{prop} \begin{enumerate}[(a)]

\item We have that
\begin{equation}
s\left(\tau + \frac{1}{2}\right) = \frac{-s(\tau)+3}{s(\tau)+1}.
\label{eqn:2.7}
\end{equation}

\item If $w = \frac{k+\sqrt{-d}}{2}$, where $k \in \mathbb{Z}$ is odd, then
$$t(s(w/2)) = \overline{s(w/2)}$$
is the complex conjugate of $s(w/2)$, and the mapping $x \rightarrow t(x) = \frac{-x+3}{x+1}$ acts on the roots of the minimal polynomial of $s(w/2)$.

\item $s^2(w) = s(w/2) \overline{s(w/2)} = |s(w/2)|^2$.

\item If $w = \frac{3a+\sqrt{-d}}{2}$, where $a$ is odd, $-d \equiv 1$ (mod $8$), and $24 \mid 9a^2+d$, then
$$t(s(w/6)) = s\left(\frac{w}{6}+\frac{1}{2}\right) = \overline{s(w/6)}$$
\end{enumerate}
\label{prop:2.5}
\end{prop}

\begin{proof}
Putting $\tau$ and $\tau + \frac{1}{2}$ in the identity \eqref{eqn:2.6} gives the equations
\begin{align*}
s(\tau) & s^2(2\tau) +s^2(\tau)+s^2(2\tau)-3s(\tau) = 0,\\
s\left(\tau + \frac{1}{2}\right) & s^2(2\tau) +s^2\left(\tau + \frac{1}{2}\right)+s^2(2\tau)-3s\left(\tau + \frac{1}{2}\right) = 0;
\end{align*}
since $s(2\tau+1) = s(2\tau)$.
These equations imply successively that
\begin{align*}
s^2(2\tau) &= \frac{-s^2(\tau)+3s(\tau)}{s(\tau)+1} = s(\tau) t(s(\tau)),\\
\left(s(\tau) - s\left(\tau + \frac{1}{2}\right)\right) &s^2(2\tau) = -s^2(\tau) + s^2\left(\tau + \frac{1}{2}\right) + 3\left(s(\tau) - s\left(\tau + \frac{1}{2}\right)\right),\\
s^2(2\tau) &= -s(\tau) - s\left(\tau + \frac{1}{2}\right) + 3,\\
s\left(\tau + \frac{1}{2}\right) &=  \frac{s^2(\tau)-3s(\tau)}{s(\tau)+1} -s(\tau)+3\\
&= \frac{-s(\tau)+3}{s(\tau)+1}.
\end{align*}
This proves (a).  To prove (b), note that
$$\frac{w+1}{2} = \frac{k+2+\sqrt{-d}}{4} = -\frac{\overline{w}}{2}+\frac{k+1}{2},$$
so that
$$t(s(w/2)) = s(w/2+1/2) = s\left(-\frac{\overline{w}}{2}+\frac{k+1}{2}\right) = s\left(-\frac{\overline{w}}{2}\right) = \overline{s(w/2)}.$$
The second assertion of (b) follows immediately, since $(x+1)^\delta f(t(x))$ ($\delta = \textrm{deg} f$) and $f(x) = \textrm{minpol}(s(w/2))$ have a root in common.
From above we also have
$$s^2(w) = s(w/2) t(s(w/2)) = s(w/2) \overline{s(w/2)} = |s(w/2)|^2,$$
which gives (c).\medskip

The proof of (d) is similar to the proof of (b):
$$t(s(w/6)) = s\left(\frac{w}{6}+\frac{1}{2}\right) = s\left(-\frac{\overline{w}}{6}+\frac{a+1}{2}\right) = s\left(-\frac{\overline{w}}{6}\right) = \overline{s(w/6)}.$$
\end{proof}

\begin{cor}
We have the idenitities
\begin{align*}
c\left(\tau+\frac{1}{2}\right) &= \frac{2e^{2\pi i/6}c(\tau)}{s(\tau)+1},\\
c_1(\tau) = c\left(\tau+\frac{3}{2}\right) &= \frac{-2c(\tau)}{s(\tau)+1},\\
s(\tau)+1 &= \frac{2c^2(\tau)}{c(2\tau)}.
\end{align*}
\label{cor:2.6}
\end{cor}

\begin{proof}
Part (a) of the proposition implies that
\begin{align*}
8c^3\left(\tau+\frac{1}{2}\right) &= 1-s^2\left(\tau+\frac{1}{2}\right)\\
&= 1- \left(\frac{-s(\tau)+3}{s(\tau)+1}\right)^2\\
&= \frac{8(s(\tau)-1)}{(s(\tau)+1)^2}\\
&= \frac{8(s(\tau)^2-1)}{(s(\tau)+1)^3} = \frac{-64c^3(\tau)}{(s(\tau)+1)^3}.
\end{align*}
From this we obtain the first identity, since the $q$-expansion of $c(\tau+1/2)$ starts with $e^{2\pi i/6} q^{1/3}$; and the $q$-expansion of $2c(\tau)/(s(\tau)+1)$
starts with $q^{1/3}$.  The second identity is immediate from the first by replacing $\tau$ by $\tau+1$ and noting $s(\tau+1) = s(\tau)$ and 
$c(\tau+1) = e^{2\pi i/3} c(\tau)$.  The third identity follows from the second and $c_1(\tau) = -c(2\tau)/c(\tau)$, which is \cite[Eq. (4.2)]{akm}.
\end{proof}

\noindent {\bf Remark.} The last relation in this corollary is equivalent to the relation
\begin{equation*}
\Big[\frac{1^3, 12}{3,4}\Big] = -\Phi+2\tilde \Phi, \ \Phi = \Big[\frac{2^2, 3^3, 12}{1, 4, 6^2}\Big], \ \tilde \Phi = \Big[\frac{1, 6^7}{2, 3^3, 12^2}\Big],
\end{equation*}
from K\"ohler's treatise \cite[Eqs. (26.65), (26.66), pp. 502-503]{ko}.

\begin{conj}
The quantity $s(w)$ is real and $\mathbb{Q}(s(w)) = \Omega^+_{12f}$, the real subfield of $\Omega_{12f}$.
\label{conj:1}
\end{conj}

\section{Values $s(w/2)$ as periodic points of an algebraic function.}

In terms of the function $s=s(\tau)$, the $j$-function has the representation
\begin{equation}
j(\tau) = \frac{(s^2 - 3)^3(s^6 - 9s^4 + 243s^2 - 243)^3}{s^{12}(s^2 - 9)^3(s^2 - 1)} = J_s(s), \ \ s = s(\tau).
\label{eqn:3.1}
\end{equation}
This follows from the fact that $s(\tau) = \sqrt{1-8c^3(\tau)}$ and
$$J_s(\sqrt{1-8c^3}) = \frac{(4c^3 + 1)^3(64c^9 + 48c^6 + 228c^3 + 1)^3}{(8c^3 - 1)^6(c^3 + 1)^3c^3}, \ \ c = c(\tau),$$
where the expression on the right is $j(\tau)$, by \cite[Eq. (3.5), p. 1718]{akm}.  The relation \eqref{eqn:3.1} shows that the field of 
modular functions $\textsf{K}_{\Gamma_0(12)}$ for $\Gamma_0(12)$ is generated over $\mathbb{C}$ by $s(\tau)$, since 
$[\textsf{K}_{\Gamma_0(12)}:\mathbb{C}(j(\tau))] = 24$.

\begin{lem} \begin{enumerate}[(a)]

\item If $-d = d_K f^2 \equiv 1$ (mod $8$) and $w = \frac{k+\sqrt{-d}}{2}$, with $k^2+d \equiv 0$ (mod $8$),
 then $r(w/2) s(w/2) \in \Omega_{3f}$ over $K = \mathbb{Q}(\sqrt{-d})$.  

\item With the same assumptions as in (a), $s^3(w/2)$ and $r^3(w/2)$ lie in $\Omega_{3f}$.

\item If $\left(\frac{d_K}{3}\right) = 0$ or $1$, where $-d = d_Kf^2 \equiv 1$ (mod $8$) and $w = \frac{k+\sqrt{-d}}{2}$, with $k^2+d \equiv 0$
 (mod $24$), then $r(w/6)s(w/6), s^3(w/6)$, and $r^3(w/6)$ lie in the ring class field $\Omega_f$.
\end{enumerate}
\label{lem:3.1}
\end{lem}

\begin{proof}
(a) From the definitions
\begin{align*}
r(\tau) &= \frac{\eta^2(\tau) \eta^4(6\tau)}{\eta^4(2\tau) \eta^2(3\tau)},\\
s(\tau) &= \frac{\eta^4(\tau) \eta^2(6\tau)}{\eta^2(2\tau) \eta^4(3\tau)},
\end{align*}
we have that
\begin{equation*}
r(\tau) s(\tau) = \frac{\eta^6(\tau) \eta^6(6\tau)}{\eta^6(2\tau) \eta^6(3\tau)}.
\end{equation*}
Thus, inserting values of the function $\gamma_3(\tau) = (j(\tau)-1728)^{1/2}$ gives that
\begin{align}
\notag r(w/2) s(w/2) &= \frac{\eta^6(w/2)}{\eta^6(3w/2)} \frac{\eta^6(3w)}{\eta^6(w)}\\
\label{eqn:3.2} &= \left(\frac{\eta(w/2)}{\eta(3w/2)}\right)^6 \gamma_3(3w/2)\\
\label{eqn:3.3} & \ \  \times \left(\frac{\eta(3w)}{\eta(w)}\right)^6 \gamma_3(3w)^{-1} \times \frac{\gamma_3(3w)}{\gamma_3(3w/2)}.
\end{align}
From Schertz \cite[p. 159]{sch} the expression in \eqref{eqn:3.2} lies in $\Omega_{3f}$, since $3w/2$ is the basis quotient for an ideal 
in the order $\textsf{R}_{-9d}$ (note that $8 \mid 9k^2+9d$).  In the same way, the first expression in \eqref{eqn:3.3} lies in $\Omega_{3f}$, as well.
For the second expression in \eqref{eqn:3.3}, we use \cite[Thm. 6.3.2, p. 145]{sch} applied to $\alpha = 3w$ and $3w/2$.  The discriminant of $\alpha$ is
$D(\alpha) = -9d = d_K (3f)^2$, which is odd.  Let $k' = 3k$ or $3k + 6$ satisfy $k' \equiv -1$ (mod $4$).  Then $k'^2+d \equiv 0$ (mod $8$).  Further,
$\gamma_3(3w)$ or $\gamma_3(3w+3) = \pm \gamma_3(3w)$ lies in $\Omega_{3f}$, by Schertz's theorem, since $3w$ or
$3w+3$ has the minimal polynomial
\begin{equation*}
m_1(X) = X^2-k' X +\frac{k'^2+9d}{4}
\end{equation*} 
over $\mathbb{Q}$.  Since $\gamma_3(\tau + l) = \pm \gamma_3(\tau)$, for $l \in \mathbb{Z}$, this shows that the quantity $\gamma_3(3w)$ in \eqref{eqn:3.3} 
lies in $\Omega_{3f}$.  To prove that $\gamma_3(3w/2)$ or $\gamma_3((3w+3)/2)$ lies in $\Omega_{3f}$ we need an additional argument, since the minimal
polynomial for these values has the form
$$m_2(X) = 2X^2 -k' X + \frac{k'^2+9d}{8},$$
with $2$ dividing the leading coefficient.  But the discriminant $-9d$ is odd, so in every ideal class in the ring $\textsf{R}_{-9d}$ there are ideals relative prime to 2.  
(See \cite[Cor. 7.17, p. 128]{cox}.)  Suppose that $(2, 3w) \sim \mathfrak{a} = (n, \frac{l+\sqrt{-9d}}{2})$, where $\mathfrak{a} \in \textsf{R}_{-9d}$ and $2 \nmid N(\mathfrak{a})$.  Then $j(3w/2) = j(\mathfrak{a})$, so that $\gamma_3(3w/2) = \pm \gamma_3(\alpha)$, where $\alpha = \frac{l+\sqrt{-9d}}{2n}$ is a basis quotient for the ideal $\mathfrak{a}$, in which $4n \mid l^2+d$.  Then the minimal polynomial of $\alpha$ or $\alpha+3$ has the form
$$m_3(X) = nX^2-l' X + \frac{l'^2+9d}{4n}, \ \ l' = l \ \textrm{or} \ l+6n \equiv -1 \ (\textrm{mod} \ 4).$$
Hence, Schertz \cite[p. 145]{sch} yields that $\gamma_3(\alpha) = \pm \gamma_3(\alpha +3) \in \Omega_{3f}$.  Similarly, if 
$$(2, 3w) \sim \mathfrak{a} = (n, l+\sqrt{-9d}) = (n, l-n + \sqrt{-9d}),$$
we may assume $l$ is odd; and then $(2,3w) \sim \mathfrak{a} = (n, \frac{l+\sqrt{-9d}}{2})$.  Thus $\gamma_3(3w/2)$ lies in $\Omega_{3f}$ in any case. \medskip

Putting all this together shows that \eqref{eqn:3.2} and \eqref{eqn:3.3} are quantities in $\Omega_{3f}$, so the same holds for $r(w/2) s(w/2)$. \medskip

Part (b) follows using the same arguments as in (a), since
\begin{align*}
s^3(w/2) &= \frac{\eta^{12}(w/2) \eta^6(3w)}{\eta^6(w) \eta^{12}(3w/2)} = \left(\frac{\eta(w/2)}{\eta(3w/2)}\right)^{12} \left(\frac{\eta(3w)}{\eta(w)}\right)^6,\\
r^3(w/2) &= \frac{\eta^{6}(w/2) \eta^{12}(3w)}{\eta^{12}(w) \eta^{6}(3w/2)} = \left(\frac{\eta(w/2)}{\eta(3w/2)}\right)^{6} \left(\frac{\eta(3w)}{\eta(w)}\right)^{12}.
\end{align*}
Part (c) is proved in the same way, using that $w/2, w/3$ and $w/6$ are basis quotients for proper ideals in the order $\textsf{R}_{-d}$, of norms $2, 3$, and $6$.  This is because, for example,
$$r(w/6)s(w/6) =  \frac{\eta^6(w/6)}{\eta^6(w/2)} \frac{\eta^6(w)}{\eta^6(w/3)},$$
so that the same arguments apply with $\textsf{R}_{-d}$ in place of $\textsf{R}_{-9d}$.
\end{proof}

\begin{thm}
If $w = \frac{k+\sqrt{-d}}{2}$, where $-d = d_K f^2 \equiv 1$ (mod $8$) and $8 \mid k^2+d$, then $s(w/2), r(w/2) \in \Omega_{3f}$, the ring class field of
conductor $3f$ over $K = \mathbb{Q}(\sqrt{-d})$.
\label{thm:3.2}
\end{thm}

\begin{proof}
Recall the equation
$$\frac{1}{r^2(\tau)} = \frac{1}{c^3(\tau)}+1 = \frac{8}{1-s^2(\tau)}+1 = \frac{9-s^2(\tau)}{1-s^2(\tau)},$$
from which we obtain the identity
\begin{equation}
r^2(\tau) = \frac{1-s^2(\tau)}{9-s^2(\tau)}.
\label{eqn:3.4}
\end{equation}
Multiplying by the denominator of the right side and rearranging gives that
$$9r^2(\tau)+s^2(\tau) = 1+r^2(\tau)s^2(\tau).$$
We conclude from Lemma \ref{lem:3.1}(a) that
$$(3r(w/2)+s(w/2))^2 = 1+6r(w/2)s(w/2)+r^2(w/2) s^2(w/2) \in \Omega_{3f}.$$
Supposing that $r(w/2)$ and $s(w/2)$ do not lie in the field $\Omega_{3f}$, it follows from Lemma \ref{lem:3.1} that $r(w/2)$ and $s(w/2)$ are both 
Kummer elements for the cubic extension $L = \Omega_{3f}(s(w/2))$.  (Recall that $\omega \in \Sigma_3 \subset \Omega_{3f}$.) 
 But then $3r(w/2)+s(w/2) = \lambda$ must lie in $\Omega_{3f}$, since it lies in $L$ and satisfies a quadratic equation over $\Omega_{3f}$.  Therefore,
$$3r(w/2) = \lambda-s(w/2), \ \ \lambda \in \Omega_{3f},$$
which implies in turn that
$$3r(w/2) s(w/2) = \lambda s(w/2) -s^2(w/2).$$
Hence, $s(w/2)$ satisfies a quadratic equation over $\Omega_{3f}$.  This gives the contradiction that $s(w/2)$ generates a quadratic extension of $\Omega_{3f}$
and proves that $s(w/2), r(w/2) \in \Omega_{3f}$.
\end{proof}

\noindent {\bf Remark.} It follows from Proposition \ref{prop:2.5}(c) and Theorem \ref{thm:3.2} that $s(w)$ is real and $s^2(w) \in \Omega_{3f}^+$.  From 
Schertz \cite[Thm. 5.2.1, p. 128]{sch} we have that $s(w) \in \Omega_{12f}$, and Conjecture \ref{conj:1} says that $\mathbb{Q}(s(w)) = \Omega_{12f}^+$.
Note that $[\Omega_{12f}^+:\Omega_{3f}^+] = 2$.  \medskip

Note that the functions $(X,Y) = (s(\tau),r(\tau))$ parametrize the elliptic curve
$$\textsf{E}: \ Y^2(9-X^2) = 1-X^2,$$
whose $j$-invariant is $j(\mathcal{E}) = \frac{2^4 13^3}{3^2}$.  In general, the elliptic curve
$$\textsf{E}_k: \ Y^2(k-X^2) = 1-X^2$$
has $j$-invariant equal to
$$j(\textsf{E}_k) = \frac{16(k^2 + 14k + 1)^3}{k(k - 1)^4};$$
and $\textsf{E}_k$ is $2$-isogenous to the Legendre normal form
$$Y^2 = X(X-1)(X-k).$$
Curiously, the algebraic function $y = \left(\frac{1-x^2}{k-x^2}\right)^{1/2}$ defined by the equation for $\textsf{E}_k$ appears to have at most $4$ periodic points altogether, which are all fixed points; and and there are exactly $4$ fixed points, if $k \neq 1, -3$. \medskip

\begin{thm}
If $w = \frac{k+\sqrt{-d}}{2}$, where $-d = d_K f^2 \equiv 1$ (mod $8$) and $8 \mid k^2+d$, then $c(w/2), c(w/4) \in \Omega_{3f}$, the ring class field of
conductor $3f$ over $K = \mathbb{Q}(\sqrt{-d})$.  Also, $s(w/4) \in \Omega_{3f}$.
\label{thm:3.3}
\end{thm}

\begin{proof}
We use the same idea as in the proof of Lemma \ref{lem:3.1}.  We have, with $\gamma_2(\tau) = j(\tau)^{1/3}$, that
\begin{align}
\notag c(\tau) &= \frac{\eta(\tau) \eta^3(6\tau)}{\eta(2\tau) \eta^3(3\tau)}\\
\notag c^8(w/2) &= \left(\frac{\eta(w/2)}{\eta(w)}\right)^8 \left(\frac{\eta(3w)}{\eta(3w/2)}\right)^{24}\\
\label{eqn:3.5} &= \left(\frac{\eta(w/2)}{\eta(w)}\right)^8 \gamma_2(w)\\
\label{eqn:3.6} & \ \ \times \left(\frac{\eta(3w)}{\eta(3w/2)}\right)^{24} \gamma_2(3w)^{-3} \times \frac{\gamma_2(3w)^3}{\gamma_2(w)}.
\end{align}
From Schertz \cite[p. 159]{sch} the expression in \eqref{eqn:3.5} and the first expression in \eqref{eqn:3.6} lie in $\Omega_{3f}$.  
Since $\gamma_2^3(3w) = j(3w)$ lies in $\Omega_{3f}$, we just have to show that $\gamma_2(w) \in \Omega_{3f}$.  This follows easily from 
Schertz's Theorem 6.3.1 in \cite[p. 144]{sch}.  Let $k' = k+2m \equiv 0$ (mod $3$), with $m \in \{0,1,2\}$.  Then the minimal polynomial of $w+m$
has the form 
$$m_4(X) = X^2-k'X+\frac{k'^2+d}{4}, \ \ k' \equiv 0 \ (\textrm{mod} \ 3).$$
Schertz's theorem then gives that $\mathbb{Q}(\gamma_2(w+m)) \subseteq \mathbb{Q}(j(3w)) \subset \Omega_{3f}$.  But  
$\gamma_2(w+m) = \omega^i \gamma_2(w)$, for some cube root of unity $\omega^i$, and since the latter lies in $\Sigma_3 \subset \Omega_{3f}$, 
this shows that $\gamma_2(w) \in \Omega_{3f}$.  Hence, $c^8(w/2) \in \Omega_{3f}$.  Theorem \ref{thm:3.2} implies that
$$s^2(w/2) = 1-8c^3(w/2) \in \Omega_{3f} \ \Rightarrow \ c^3(w/2) \in \Omega_{3f}.$$
Hence, $c(w/2) = \frac{c^9(w/2)}{c^8(w/2)} \in \Omega_{3f}$.  Now use \eqref{eqn:2.2},
$$r(w/2) = \frac{c^2(w/2)}{c(w/4)+c^2(w/2)},$$
and Theorem \ref{thm:3.2} to conclude that $c(w/4) \in \Omega_{3f}$.  Now we have the relation
$$s^2(w/4) = 1-8c^3(w/4),$$
and equation \eqref{eqn:2.6} gives
$$s(w/4) = -\frac{s^2(w/4)+s^2(w/2)}{s^2(w/2)-3}.$$
Taken together, these facts imply that $s(w/4) \in \Omega_{3f}$, which completes the proof. \medskip
\end{proof}

For the next proposition we need some notation.  By $\alpha \cong p^{m/n}$, for $p$ a rational prime (usually $p = 3$), $(m,n) = 1$, and $\alpha \in L$, a normal extension
of $\mathbb{Q}$, we mean that
$$(\alpha) = \mathfrak{a}^m, \ \ \textrm{where} \ \ \mathfrak{a}^n = (p);$$
and $\mathfrak{a}$ is an integral ideal in $R_L$.  Thus, $n$ divides the ramification index of $p$ in $L/\mathbb{Q}$ and 
$\mathfrak{a} = \prod_{\mathfrak{p}_i \mid p}{\mathfrak{p}_i^{e_1}}$, where the $\mathfrak{p}_i$ are prime ideals and $e_1 n$ is equal to the ramification index
of $p$ in $L/\mathbb{Q}$.  In this situation $\textsf{N}(\alpha)^n \cong \textsf{N}(p^m)$, so $\textsf{N}(\alpha) \cong p^{m[L:\mathbb{Q}]/n}$, if $\textsf{N}$ 
denotes the absolute norm from $L$.

\begin{prop}
Assume that $-d = d_K f^2 \equiv 1$ (mod $8$) and $3^t || f$.  Then in the ring class field $\Omega_{3f}$ over $K = \mathbb{Q}(\sqrt{-d})$, we have
\begin{equation}
s(w/2) \cong \begin{cases}
3, \ &  \ \left(\frac{d_K}{3}\right) = +1,\\
3^{1-\frac{1}{2 \cdot 3^{t+1}}}, \ & \ 3 \mid d_K,\\
3^{1-\frac{1}{4 \cdot 3^t}}, \ & \ \left(\frac{d_K}{3}\right) = -1.
\end{cases}
\label{eqn:3.7}
\end{equation}
Moreover, if $\textsf{N}$ denotes the norm from $\Omega_{3f}$ to $\mathbb{Q}$ and $h = h(-d)$, we have that
\begin{equation}
\textsf{N}(s(w/2)) = \begin{cases}
3^{4h} \ & \ \left(\frac{-d}{3}\right) = +1,\\
3^{h(6-\frac{1}{3^t})}, \ & \ 3 \mid d_K,\\
3^{6h}, \ & \ \left(\frac{d_K}{3}\right) = +1, \ 3 \mid f;\\
3^{h(6-\frac{3}{2 \cdot 3^t})}, \ & \  \left(\frac{d_K}{3}\right) = -1, \ 3 \mid f;\\
3^{6h}, \ & \ \left(\frac{-d}{3}\right) = -1.
\end{cases}
\label{eqn:3.8}
\end{equation}
In particular, $3 \mid s^2(w/2)$, so that $s(w/2)$ is never a unit.
\label{prop:3.4}
\end{prop}

\begin{proof}
Writing $s(\tau)$ in terms of the $\eta$-function, we have
\begin{align*}
s(\tau) &= \frac{\eta^4(\tau) \eta^2(6\tau)}{\eta^2(2\tau) \eta^4(3\tau)},\\
s(w/2) &= \left(\frac{\eta(w/2)}{\eta(3w/2)}\right)^4 \left(\frac{\eta(3w)}{\eta(w)}\right)^2,
\end{align*}
\begin{align*}
s(w/2)^{12} &= \left(\frac{\eta(w/2)}{\eta(3w/2)}\right)^{48} \left(\frac{\eta(3w)}{\eta(w)}\right)^{24}\\
&= \frac{\Delta(w,2)^2}{\Delta(3w,2)^2} \frac{\Delta(3w,1)}{\Delta(w,1)}\\
&= 3^{12} \cdot \left(3^{-12}\frac{\Delta(w,2)}{\Delta(3w,2)}\right)^2 \left(3^{12} \frac{\Delta(3w,1)}{\Delta(w,1)}\right)\\
&= 3^{12} \frac{\varphi_P(w)}{\varphi_P(w/2)^2}.
\end{align*}
Here, $\Delta(w_1,w_2) = (2 \pi)^{12} w_2^{-12} \eta(w_1/w_2)^{24}$ is Hasse's notation \cite[Eq. ($2^\prime$), p. 119]{hbeg} and 
$P = \begin{pmatrix}
   3   &  0  \\
    0  &  1
\end{pmatrix}$
is a primitive matrix of determinant $3$ for which $P\begin{pmatrix} w \\ 1 \end{pmatrix} = \begin{pmatrix} 3w \\ 1 \end{pmatrix}$ and $P\begin{pmatrix} w \\ 2 \end{pmatrix} = \begin{pmatrix} 3w \\ 2 \end{pmatrix}$ are bases for proper ideals in $\textsf{R}_{3f}$.  From Deuring's treatise \cite[p. 43]{deu2}, with $\alpha = w$, resp. $w/2$, we have
\begin{equation*}
\varphi_P(\alpha) \cong \begin{cases} 1, & \left(\frac{d_K}{3}\right) = +1,\\
3^{6/3^{t+1}}, & 3 \mid d_K,\\
3^{12/(4\cdot 3^t)}, & \left(\frac{d_K}{3}\right) = -1.
\end{cases}
\end{equation*}
Since $\varphi_P(w) \cong \varphi_P(w/2)$, this gives \eqref{eqn:3.7} on taking $12$-th roots.  Now \eqref{eqn:3.8} follows from the fact that
\begin{align*}
[\Omega_{3f}:\mathbb{Q}] &= [\Omega_{3f}:\Omega_f] [\Omega_f:K] [K:\mathbb{Q}]\\
&=  [\Omega_{3f}:\Omega_f] \cdot h(-d) \cdot 2\\
&= \begin{cases} 4h, &  \left(\frac{d_K}{3}\right) = +1, 3 \nmid f;\\
6h, & 3 \mid d_K \ \textrm{or} \ 3 \mid f;\\
8h, & \left(\frac{d_K}{3}\right) = -1, 3 \nmid f.
\end{cases}
\end{align*}
It is clear that the exponent of $3$ on the right side of \eqref{eqn:3.7} is $> \frac{1}{2}$, implying that $3 \mid s^2(w/2)$.  Further $s(w/2) \cong 3$ if and only if $\left(\frac{d_K}{3}\right) = +1$.
\end{proof}

\begin{cor}
With assumptions as in the theorem, $r(w/2) s(w/2)$ is always a unit in $\Omega_{3f}$.
\label{cor:3.5}
\end{cor}

\begin{proof}
Since
$$r(\tau) s(\tau) = \frac{\eta^6(\tau) \eta^6(6\tau)}{\eta^6(2\tau) \eta^6(3\tau)},$$
we find as in the above proof that
$$(r(w/2) s(w/2))^4 = \frac{\varphi_P(w)}{\varphi_P(w/2)} \cong 1.$$
\end{proof}

\begin{prop}
Assume that $-d = d_K f^2 \equiv 1$ (mod $8$), with $3^t || f$.  The ramification index of the prime divisors of $3$ in $\Omega_{3f}/\mathbb{Q}$ is $2 \cdot 3^{t+1}$, resp. $4 \cdot 3^t$, according as $\left(\frac{d_K}{3}\right) = 0$ or $-1$.
\label{prop:3.6}
\end{prop}

\begin{proof}
Equation \eqref{eqn:3.7} implies that the ramification index of prime divisors of $3$ in $\Omega_{3f}/\mathbb{Q}$ is at least 
$2 \cdot 3^{t+1}$, resp. $4 \cdot 3^t$, when $\left(\frac{d_K}{3}\right) = 0$, resp. $-1$.  Let $f = 3^t f_0$, with $(3,f_0) = 1$.  
Consider the chain of normal extensions
\begin{equation}
\mathbb{Q} \subset K \subseteq \Omega_{f_0} \subseteq \Omega_{f} \subset \Omega_{3f}.
\label{eqn:3.9}
\end{equation}
If $3 \mid d_K$, then the sequence of ramification indices of prime divisors of $3$ in this chain of extensions is
$$2, 1, e(\Omega_f:\Omega_{f_0}), e(\Omega_{3f}:\Omega_f).$$
Since $e(\Omega_f:\Omega_{f_0}) \mid 3^t$ and $e(\Omega_{3f}:\Omega_f) \mid [\Omega_{3f}:\Omega_f] = 3$,
it follows that $e(\Omega_f:\mathbb{Q}) \mid 2 \cdot 3^t$, whence $e(\Omega_{3f}:\Omega_f) = 3$ and $e(\Omega_f:\Omega_{f_0}) = 3^t$.  
\medskip

If $\left(\frac{d_K}{3}\right) = -1$, then the sequence of ramification indices of prime divisors of $3$ in the chain \eqref{eqn:3.9} is
$$1, 1, e(\Omega_f:\Omega_{f_0}), e(\Omega_{3f}:\Omega_f),$$
where
$$e(\Omega_f:\Omega_{f_0}) \mid [\Omega_f:\Omega_{f_0}] = 4\cdot 3^{t-1}.$$
Thus $e(\Omega_f:\mathbb{Q}) \mid 4 \cdot 3^{t-1}$, giving that $e(\Omega_{3f}:\Omega_f) = 3$ and $e(\Omega_f:\Omega_{f_0}) = 4 \cdot 3^{t-1}$,
 as before.  
\end{proof}

\noindent {\bf Remark.} This shows that $\Omega_f$ and $\Omega_{3f}$ are totally ramified over $\Omega_{f_0}$ at the prime divisors of $3$.  See
Hasse \cite[p. 332]{hzer}, but note that Hasse excludes a finite set of primes, including $2$ and $3$, from consideration. \medskip

Now we make use of the modular functions
\begin{align*}
u(\tau) &= \sqrt{2}q^{1/8} \prod_{n = 1}^\infty{\left(\frac{1+q^{2n}}{1+q^{2n-1}}\right)}, \ \ q = e^{2\pi i \tau},\\
v(\tau) &= q^{1/2} \prod_{n=1}^\infty{(1-q^n)^{\left(\frac{2}{n}\right)}}.
\end{align*}
From \cite[Props. 3.1, 10.1]{akm0} we take the following identities.  Setting $x = u(\tau), y = u(2\tau)$, we first have
\begin{equation}
x^4(y^4+1) = 2y^2, \ \ x = u(\tau), \ y=u(2\tau).
\label{eqn:3.10}
\end{equation}
Also, we have
\begin{equation}
u^4(v^2+1)^2+4v(v^2-1) = 0, \ \ u = u(\tau), \ v = v(\tau).
\label{eqn:3.11}
\end{equation}
See \cite[p. 153]{duk} and \cite[Props. 3.1,10.1]{akm0}.  Finally, if
$$w = \frac{k+\sqrt{-d}}{2}, \ \ k^2+d \equiv 0 \ (\textrm{mod} \ 32), \ \ (N(w),f) = 1,$$
then from \cite[eq. (6.6)]{akm0} the unit $v(w/8)$ is given by
\begin{equation}
v(w/8) = \pm \frac{1 \pm \sqrt{1+\pi^2}}{\pi}, \ \ \pi = \pm \frac{\mathfrak{f}_2^2(w/2)}{\mathfrak{f}^2(w/2)},
\label{eqn:3.12}
\end{equation}
where $\mathfrak{f}(\tau), \mathfrak{f}_1(\tau), \mathfrak{f}_2(\tau)$ are the Schl\"afli-Weber functions.  (See \cite[p. 233]{cox}.)
Furthermore, if $(2,w) = \mathfrak{p}_2$ and $\mathfrak{p}_2 \mathfrak{p}_2' = (2) \ \textrm{in} \ R_K$, then $(\pi) = \mathfrak{p}_2$ in $\Omega_f$ and 
$\mathbb{Q}(\pi) = \mathbb{Q}(\xi) = \mathbb{Q}(\pi^4) = \mathbb{Q}(\xi^4) = \Omega_f$, where
$$\pi^4 + \xi^4 = 1, \ \ \pi \cong \mathfrak{p}_2, \ \xi \cong \mathfrak{p}_2', \ \xi=\pm i \frac{\mathfrak{f}_1^2(w/2)}{\mathfrak{f}^2(w/2)}.$$
These last relations are proved in \cite[Section 10]{lm}.

\begin{prop}
For $\tau$ in the upper half-plane, we have
\begin{equation}
u^8(\tau) = \frac{(s(\tau)-1)(s(\tau)+3)^3}{(s(\tau)+1)(s(\tau)-3)^3}.
\label{eqn:3.13}
\end{equation}
\label{prop:3.7}
\end{prop}

\begin{proof}
Let $F(x)$ denote the rational function
$$F(x) = \frac{16(x^2+14x+1)^3}{x(x-1)^4}.$$
Then equation (9.6) in Duke's paper \cite[p. 153]{duk} implies that $j(\tau) = F(u^8(\tau))$.  We compute using Maple that 
$$F\left(\frac{(s-1)(s+3)^3}{(s+1)(s-3)^3}\right) = \frac{(s^2 - 3)^3(s^6 - 9s^4 + 243s^2 - 243)^3}{s^{12}(s^2 - 9)^3(s^2 - 1)} = j(\tau), \ s = s(\tau),$$
by \eqref{eqn:3.1}.  Hence,
$$F(u^8(\tau)) = F\left(\frac{(s-1)(s+3)^3}{(s+1)(s-3)^3}\right) = F(A(s)), \ \ s=s(\tau).$$
On the other hand, we have the algebraic identity
\begin{align}
\label{eqn:3.13a} F(x)-F(y) &= \frac{16(x-y)(xy-1)}{x(x - 1)^4y(y - 1)^4}G(x,y),\\
\notag G(x,y) &= x^4y^4 - 4x^4y^3 - 4x^3y^4 + 6x^4y^2 - 752x^3y^3 + 6x^2y^4 - 4x^4y \\
\notag & \ \ - 2584x^3y^2 - 2584x^2y^3 - 4xy^4 + x^4 - 752x^3y + 13348x^2y^2 \\
\notag & \ \ - 752xy^3 + y^4 - 4x^3 - 2584x^2y - 2584xy^2 - 4y^3 + 6x^2\\
\notag & \ \  - 752xy + 6y^2 - 4x - 4y + 1.
\end{align}
Now we use that $u^8(\tau) \simeq 16q$ and
$$A(s) = \frac{(s-1)(s+3)^3}{(s+1)(s-3)^3} \simeq 16q,$$
as $q \rightarrow 0$, which follows from
$$s(\tau) = 1 - 4q + 4q^2 + 4q^3 - 12q^4 + 8q^5 + 12q^6 + \cdots = 1 - 4q +O(q^2).$$
Then $G(u^8(\tau), A(s(\tau)))  \rightarrow 1$ and $u^8(\tau) A(s(\tau)) -1 \rightarrow -1$, as $q \rightarrow 0$.  This implies that
the factor $(x-y)$ must be $0$, for $\tau$ in the upper half-plane, i.e., $u^8(\tau) = A(s(\tau))$.
\end{proof}

\begin{prop}
If $w = \frac{k+\sqrt{-d}}{2}$ is chosen so that $k^2+d \equiv 0$ (mod $2^5$), then
\begin{equation}
u^4(w/2) = \frac{1 + \xi^2}{1- \xi^2},
\label{eqn:3.14}
\end{equation}
and $u^8(w/2)$ generates $\Omega_f$ over $\mathbb{Q}$.
\label{prop:3.8}
\end{prop}

\begin{proof}
Starting with \eqref{eqn:3.11} and \eqref{eqn:3.12}, we find first with $v = v(w/8)$ that
\begin{align}
\notag u^4(w/8) &= \frac{-4v(v^2-1)}{(v^2+1)^2}\\
\notag &= \pm \frac{2(1 \pm \sqrt{\pi^2 + 1})^2\pi}{(1+\pi^2 \pm \sqrt{\pi^2 + 1})^2}\\
\label{eqn:3.15} &= \pm \frac{2\pi}{\pi^2+1}.
\end{align}
With this information we use \eqref{eqn:3.10} and solve the equation
$$x^4(y^4+1) = 2y^2, \ \ x = u(w/8), \ y = u(w/4),$$
for $y^4 = u^4(w/4)$ in terms of $x^4 = u^4(w/8)$, and find that
$$u^4(w/4) = \pi^2 \ \textrm{or} \ \frac{1}{\pi^2}.$$
Now we solve the same equation for $y^4 = u^4(w/2)$ in terms of $x^4 = u^4(w/4)$ and find either that
\begin{align*}
u^4(w/2) &= \frac{2}{\pi^4} \pm \frac{2\sqrt{1-\pi^4}}{\pi^4} - 1\\
&= 2\frac{1 \pm \xi^2}{1-\xi^4} - 1\\
&= \frac{2}{1 \mp \xi^2} - 1 = \frac{1 \pm \xi^2}{1 \mp \xi^2}, \ \ \textrm{if} \ u^4(w/4) = \pi^2;
\end{align*}
or
\begin{align*}
u^4(w/2) &= 2\pi^4 -1 \pm 2\pi^2 \sqrt{\pi^4 - 1}\\
&= 2\pi^4-1 \pm 2i \pi^2 \xi^2, \ \ \textrm{if} \ u^4(w/4) = \pi^{-2}.
\end{align*}
However, if the latter formula holds, then \eqref{eqn:3.13} would imply that
$$i = \sqrt{-1} \in \mathbb{Q}(s(w/2),\pi^2,\xi^2) \subseteq \Omega_{3f},$$
which is impossible, since in $\Omega_{3f}$ the prime $2$ is not ramified.  Hence, the first formula must hold, giving
$$u^4(w/2) =  \frac{1 \pm \xi^2}{1 \mp \xi^2}.$$

From above we have the relation
$$u^8(w/4) = \pi^4,  \ \ w = \frac{k+\sqrt{-d}}{2}, \ k^2+d \equiv 0 \ (\textrm{mod} \ 2^5), \ \frak{p}_2 \mid w.$$
From \cite[Eq. (9), p. 333]{mdeu} we take the relation
\begin{equation*}
\pi^\tau = \frac{-1+\xi^2}{\pi^2}, \ \ \tau = \left(\frac{\Omega_f/K}{\frak{p}_2}\right) = \left(\frac{\Omega_{3f}/K}{\frak{p}_2}\right)\Big|_{\Omega_f}.
\end{equation*}
Hence,
\begin{equation*}
(\pi^4)^\tau = \frac{(1-\xi^2)^4}{\pi^8} = \frac{(1-\xi^2)^4}{(1-\xi^4)^2} = \left(\frac{1-\xi^2}{1+\xi^2}\right)^2.
\end{equation*}
In particular, $\pi^{4\tau}$ is an algebraic integer, while its reciprocal is not.  It follows that
\begin{equation}
u^8(w/4)^\tau = \left(\frac{1-\xi^2}{1+\xi^2}\right)^2 = \frac{(s(w/4)^{\bar \tau}-1)(s(w/4)^{\bar \tau}+3)^3}{(s(w/4)^{\bar \tau}+1)(s(w/4)^{\bar \tau}-3)^3},
\label{eqn:3.17a}
\end{equation}
where
$$\bar \tau =  \left(\frac{\Omega_{3f}/K}{\frak{p}_2}\right).$$
On the other hand,
\begin{equation}
u^8(w/2) = \left(\frac{1 \pm \xi^2}{1 \mp \xi^2}\right)^2 = \frac{(s(w/2)-1)(s(w/2)+3)^3}{(s(w/2)+1)(s(w/2)-3)^3}.
\label{eqn:3.17b}
\end{equation}

Now, the function $u^2(\tau)$ has the infinite product expansion
$$u^2(\tau) = 2 q^{1/4} \prod_{n=1}^\infty{\left(\frac{1+q^{2n}}{1+q^{2n-1}}\right)^2}.$$
From \cite{akm0} we see that this is related to the function
$$\mathfrak{p}(\tau) = 2q^{1/16} \prod_{n=1}^\infty{\left(\frac{1+q^{n/2}}{1+q^{n/2-1/4}}\right)^2},$$
for which
$$\mathfrak{p}(4\tau) = 2q^{1/4} \prod_{n=1}^\infty{\left(\frac{1+q^{2n}}{1+q^{2n-1}}\right)^2}.$$
This gives the identity
\begin{equation*}
u^2(\tau) = \mathfrak{p}(4\tau).
\end{equation*}
From this and \cite[Thm. 7.2, p.807]{akm0} we obtain
\begin{equation}
u^2(w/2) = \mathfrak{p}(2w) = \frac{-1}{\pi^{\tau}}, \ \ \tau = \left(\frac{\Omega_f/K}{\mathfrak{p}_2}\right).
\label{eqn:3.18}
\end{equation}
This shows that $u^8(w/2)$ is not an algebraic integer, so we have
\begin{align}
\label{eqn:3.18a} u^8(w/2) &= \left(\frac{1+\xi^2}{1-\xi^2}\right)^2 = u^{-8}(w/4)^\tau;\\
\notag u^4(w/2) &= \frac{1+\xi^2}{1-\xi^2}.
\end{align}
Now,
$$u^8(w/2) = \left(\frac{1 + \xi^2}{1 - \xi^2}\right)^2 = \frac{1+2\xi^2+\xi^4}{1-2\xi^2+\xi^4} = \frac{\xi^2+\frac{1}{\xi^2}+2}{\xi^2+\frac{1}{\xi^2}-2}.$$
This quantity generates the same extension over $\mathbb{Q}$ that $\alpha = \xi^2+\frac{1}{\xi^2}$ does.  However, it is easy to see that $\alpha$ must generate the same extension over $\mathbb{Q}$ as $\xi^2$ does, namely $\Omega_f$.  If, for example, $\sigma$ is an automorphism of $\Omega_{f}$ that fixes $\alpha$, then
$\alpha^\sigma = \alpha$ implies that
$$0 = \xi^{2\sigma} + \frac{1}{\xi^{2\sigma}} - \xi^2 - \frac{1}{\xi^2} = \frac{(\xi^{2\sigma} \xi^2 - 1)(\xi^{2\sigma} - \xi^2)}{\xi^{2+2\sigma}}.$$
This implies $\xi^{2\sigma} = \xi^2$ or $\xi^{-2}$.  The latter is impossible, since $\xi^2 \cong \mathfrak{p}_2'^2$.  Hence, $\xi^{2\sigma} = \xi^2$
implies $\sigma = 1$, since $\mathbb{Q}(\xi^2) = \Omega_f$. 
This completes the proof of the proposition.
\end{proof}

\begin{cor} With $w = \frac{k+\sqrt{-d}}{2}$ as in the proposition,
\begin{equation}
\frac{(s(w/2)-1)(s(w/2)+3)^3}{(s(w/2)+1)(s(w/2)-3)^3} = \frac{(-s(w/4)^{\bar \tau}-1)(-s(w/4)^{\bar \tau}+3)^3}{(-s(w/4)^{\bar \tau}+1)(-s(w/4)^{\bar \tau}-3)^3}.
\label{eqn:3.19}
\end{equation}
\label{cor:3.9}
\end{cor}

\begin{proof}
This follows directly from \eqref{eqn:3.17a}, \eqref{eqn:3.17b}, and \eqref{eqn:3.18a}.
\end{proof}

Setting
$$A(x) = \frac{(x-1)(x+3)^3}{(x+1)(x-3)^3},$$
this corollary gives that
$$A(s(w/2)) = A(-s(w/4)^{\bar \tau}).$$
We would like to conclude that $s(w/2) = -s(w/4)^{\bar \tau}$. For this we need the identity
\begin{equation}
A(x)-A(y) = \frac{16(y - x)(x^3y^3 - 18x^2y^2 + 27x^2 + 27xy + 27y^2)}{(x + 1)(x - 3)^3(y + 1)(y - 3)^3}.
\label{eqn:3.20}
\end{equation}
Let
\begin{equation}
L(x,y) = x^3y^3 - 18x^2y^2 + 27x^2 + 27xy + 27y^2.
\label{eqn:3.22a}
\end{equation}

From Theorem \ref{thm:3.2}, Proposition \ref{prop:3.8} and the identity
\begin{equation}
j(3\tau) = \frac{(s^2 - 3)^3(s^6 - 9s^4 + 3s^2 - 3)^3}{s^4(s^2 -9)(s^2 - 1)^3}, \ \ s = s(\tau),
\label{eqn:3.16}
\end{equation}
we obtain the following.  This identity follows from \cite[Eq. (3.6), p. 1718]{akm} and $s(\tau) = \sqrt{1-8c^3(\tau)}$.

\begin{thm}
If $w = \frac{k+\sqrt{-d}}{2}$, where $-d = d_K f^2 \equiv 1$ (mod $8$) and $8 \mid k^2+d$, then $\mathbb{Q}(s(w/2)) = \mathbb{Q}(r(w/2)) = \Omega_{3f}$.
\label{thm:3.10}
\end{thm}

\begin{proof}
We can assume that $k^2+d \equiv 0$ (mod $32$), for the following reason.  The ideal $\mathfrak{p}_2 = (2,w) = (2, w+2l)$, where
$w + 2l = \frac{k+4l + \sqrt{-d}}{2}$ and $l \in \mathbb{Z}$.  Also, $s(w/2+l) = s(w/2)$.  Given that $k^2+d \equiv 0$ (mod $8$), we want to choose $l$ so that
\begin{align*}
(k+4l)^2+d &\equiv 0 \ (\textrm{mod} \ 32) \ \iff \\
k^2+ 8kl + 16l^2 + d &\equiv 0 \ (\textrm{mod} \ 32) \ \iff\\
2l^2 + kl &\equiv -\frac{k^2+d}{8} \ (\textrm{mod} \ 4).
\end{align*}
To satisfy this congruence, let $k' = \frac{k^2+d}{8}$, and choose
$$l \ \begin{cases}
 \ = 0, \ & \textrm{if} \ 4 \mid k',\\
\ = 2, \ & \textrm{if} \ k' \equiv 2 \ (\textrm{mod} \ 4),\\
\ \equiv \left(-k'+2\right)k \ (\textrm{mod} \ 4), & \textrm{if} \ k' \equiv 1 \ (\textrm{mod} \ 2).
\end{cases} $$
Since $k$ is odd, it is easy to check that $2l^2+kl \equiv -k'$ (mod $4$), as desired.
Now note that \eqref{eqn:3.13} and Proposition \ref{prop:3.8} imply that $\Omega_f \subset \mathbb{Q}(s(w/2))$.  Moreover, $j(3w/2) \in \mathbb{Q}(s(w/2))$
by \eqref{eqn:3.16}.  Also, $3w/2$ is a basis quotient for a proper ideal of the order $\textsf{R}_{-9d}$, so that $K(j(3w/2)) = \Omega_{3f}$.  Hence,
$$\Omega_{3f} = \Omega_{f}(j(3w/2)) \subseteq \mathbb{Q}(s(w/2)) \subseteq \Omega_{3f}.$$
This proves that $\mathbb{Q}(s(w/2)) = \Omega_{3f}$.  We postpone the proof that $\mathbb{Q}(r(w/2)) = \Omega_{3f}$ to the last part of the proof of Theorem \ref{thm:4.3} below.  We will not use this fact in the meantime.
\end{proof}

\begin{lem}
If $w = \frac{a+\sqrt{-d}}{2}$, where $-d = d_K f^2 \equiv 1$ (mod $8$) and $8 \mid a^2+d$, so that $\frak{p}_2 \mid w$, then
$$c(w/2) \cong \mathfrak{p}_2'^{-1}, \ \ 2c(w/2) \cong \mathfrak{p}_2.$$
If $16 \mid a^2+d$, then we also have that $c(w/4) \cong \mathfrak{p}_2'^{-1}$ and $2c(w/4) \cong \mathfrak{p}_2$.  If $\left(\frac{d_K}{3}\right) = 0$
or $1$ and $24 \mid a^2+d$, then $c(w/6) \cong \mathfrak{p}_2'^{-1}$, as well.
\label{lem:3.11}
\end{lem}

\begin{proof}
We have, with the notation $\Delta(w_1,w_2) = (2 \pi)^{12} w_2^{-12} \eta(w_1/w_2)^{24}$, 
\begin{align*}
\notag c(\tau) &= \frac{\eta(\tau) \eta^3(6\tau)}{\eta(2\tau) \eta^3(3\tau)}\\
\notag c^{24}(w/2) &= \left(\frac{\eta(w/2)}{\eta(w)}\right)^{24} \left(\frac{\eta(3w)}{\eta(3w/2)}\right)^{72}\\
&= 2^{12} \frac{\Delta(w,2)}{\Delta(w,1)} \left(2^{-12} \frac{\Delta(3w,1)}{\Delta(3w,2)}\right)^3\\
&= \frac{\varphi_P(w)}{\varphi_P(3w)^3}.
\end{align*}
Here $P = \begin{pmatrix}
   1   &  0  \\
    0  &  2
\end{pmatrix}$
is a primitive matrix for which $P\begin{pmatrix} w \\ 1 \end{pmatrix} = \begin{pmatrix} w \\ 2 \end{pmatrix}$ is a basis for the ideal
$\mathfrak{p}_{2,d} = \mathfrak{p}_2 \cap \textsf{R}_{-d}$ in $\textsf{R}_{-d}$ and 
$P\begin{pmatrix} 3w \\ 1 \end{pmatrix} = \begin{pmatrix} 3w \\ 2 \end{pmatrix}$ is a basis for the ideal $\mathfrak{p}_{2,9d}$ in 
$\textsf{R}_{-9d}$.  From Deuring's treatise \cite[p. 43]{deu2} 
we see that $\varphi_P(w) \cong \mathfrak{p}_2'^{12} \cong \varphi_P(3w)$.  Hence,
$$c^{24}(w/2) \cong \frac{\mathfrak{p}_2'^{12}}{\mathfrak{p}_2'^{36}} \cong \mathfrak{p}_2'^{-24}.$$
Hence, $c(w/2) \cong \mathfrak{p}_2'^{-1}$ in $\Omega_{3f}$.  In the same way, we find that $c(w/4) \cong \mathfrak{p}_2'^{-1}$ if $16 \mid a^2 + d$.  The
last assertion follows in the same way on replacing $w$ by $w/3$ in the above $\eta$-products, so that
$$c^{24}(w/6) = \frac{\varphi_P(w/3)}{\varphi_P(w)^3};$$
and here, $w/3$ and $w$ are both basis quotients for ideals in $\textsf{R}_{-d}$.
\end{proof}

We now prove the following.

\begin{thm}
If $-d \equiv 1$ (mod $8$) and $w = \frac{k+\sqrt{-d}}{2}$, with $k^2+d \equiv 0$ (mod $2^5$) and $\frak{p}_2 \mid w$, then
$$-s(w/4)^{\bar \tau} = s(w/2), \ \ \bar \tau =  \left(\frac{\Omega_{3f}/K}{\frak{p}_2}\right).$$
\label{thm:3.12}
\end{thm}

\begin{proof}
Let (see the equation just before Lemma \ref{lem:3.1}) 
$$J(x) = \frac{(64x^9+48x^6+228x^3+1)^3(4x^3+1)^3}{x^3(x^3+1)^3(8x^3-1)^6}.$$
Equation (3.5) in \cite[p. 1718]{akm} says that
$$J(c(\tau)) = j(\tau).$$
It follows from this that
\begin{align}
\notag j(\tau) &= J_2(8c^3(\tau)),\\
J_2(x) &= \frac{(x^3 + 6x^2 + 228x + 8)^3(x + 2)^3}{x(x + 8)^3(x - 1)^6}.
\label{eqn:3.21}
\end{align}
From \cite[Prop. 3, p. 1715]{akm} we have the identity
$$c^3\left(\frac{-1}{6\tau}\right) = \frac{1/8-c^3(\tau)}{1+c^3(\tau)},$$
whence it follows that
$$8c^3\left(\frac{-1}{6\tau}\right) = 8\frac{1-8c^3(\tau)}{8+8c^3(\tau)}.$$
This inspires us to define
$$B(x) = 8\frac{1-x}{8+x}.$$
Plugging this into $J_2(x)$ yields the expression
$$J_2(B(x)) = -\frac{(x - 4)^3(x^3 + 12x^2 + 48x - 64)^3}{x^6(x - 1)(x + 8)^2} = J_3(x)$$
and the formula
$$j(6\tau) = j\left(\frac{-1}{6\tau}\right) = J_2(B(8c^3(\tau))) = J_3(8c^3(\tau)).$$
Thus,
$$j(3w) = J_3(8c^3(w/2)) \ \textrm{and} \ j(3w/2) = J_3(8c^3(w/4)).$$
It follows that
$$ J_3(8c^3(w/2)) = j(3w) = j(3w/2)^{\bar \tau}  = J_3(8c^3(w/4)^{\bar \tau}) .$$
This allows us to get a second relation satisfied by $s(w/2)$ and $-s(w/4)^{\bar \tau}$, by computing
\begin{align*}
J_3(x) - J_3(y) &= \frac{(y-x)}{y^6(y - 1)(y + 8)^2x^6(x - 1)(x + 8)^2} (x^2y^2 - 48xy + 64x + 64y)\\
& \times (x^3y + x^2y^2 + xy^3 - x^3 + 7x^2y + 7xy^2 - y^3 - 8x^2 - 8xy - 8y^2)\\
& \times (x^6y^6 +16 R(x,y)), \ \ R(x,y) \in \mathbb{Z}[x,y].
\end{align*}
Letting $x = 8c^3(w/2), y = 8c^3(w/4)^{\bar \tau}$, $x$ and $y$ satisfy one of the four relations on the right of the above equation.  They cannot
satisfy the second or fourth relation, because both of these relations would imply that $2 \mid \mathfrak{p}_2^{12}$ or $2 \mid \mathfrak{p}_2^{36}$,
by Lemma \ref{lem:3.11}; and both are impossible.  Hence, we either have
\begin{equation}
x = 8c^3(w/2) = 8c^3(w/4)^{\bar \tau} = y,
\label{eqn:3.22}
\end{equation}
{\it or}
\begin{equation}
\mathcal{K}(x,y) = x^3y + x^2y^2 + xy^3 - x^3 + 7x^2y + 7xy^2 - y^3 - 8x^2 - 8xy - 8y^2 = 0.
\label{eqn:3.23}
\end{equation}
Using the relation $8c^3(\tau) = 1-s^2(\tau)$, we replace $x,y$ in \eqref{eqn:3.23} by $x=1-s^2, y=1-t^2$, giving
\begin{align}
\label{eqn:3.24} \mathcal{K}(1-s^2,1-t^2) &= (s^3t - s^2t^2 + st^3 - 6st - 3)(s^3t + s^2t^2 + st^3 - 6st + 3),\\
\notag s &= s(w/2), \ \ t = -s(w/4)^{\bar \tau}.
\end{align}
Taking the resultant of $L(s,t)$ (from \eqref{eqn:3.22a}) and the first factor in \eqref{eqn:3.24} gives
\begin{align}
\notag \textrm{Res}_s(L(s,t),&s^3t - s^2t^2 + st^3 - 6st - 3) = -27t(t^8 - 16t^6 + 86t^4 - 144t^2 + 81)\\
& \ \ \times (t^8 - 20t^6 + 182t^4 - 180t^2 + 81);
\label{eqn:3.25}
\end{align}
while the resultant of $L(s,t)$ with the second factor is
\begin{align}
\notag \textrm{Res}_s(L(s,t),&s^3t + s^2t^2 + st^3 - 6st + 3) = -27t(t^2 - 1)(t^2 - 9)(t^4 - 2t^2 + 9)\\
 & \ \ \times (t^4 - 4t^3 + 14t^2 + 12t + 9)(t^4 + 4t^3 + 14t^2 - 12t + 9).
\label{eqn:3.26}
\end{align}
For each of the polynomials $q(t)$ in \eqref{eqn:3.25} and \eqref{eqn:3.26}, excluding trivial factors, compute the resultant
\begin{equation}
R_d(X) = \textrm{Res}_t(q(t),t^{12}(t^2 - 9)^3(t^2 - 1)X - (t^2 - 3)^3(t^6 - 9t^4 + 243t^2 - 243)^3).
\label{eqn:3.27}
\end{equation}
By \eqref{eqn:3.1}, this yields the corresponding class equation $H_{-d}(X)$ and discriminant $-d$.  These are listed in the adjoining Table \ref{tab:4}, and show that none of these polynomials occur for discriminants $-d \equiv 1$ mod $8$.  (See \cite{cox} and \cite{mcong} for the class equations occurring in this table.) Hence, if $L(s,t) = 0$, then $\mathcal{K}(s,t) \neq 0$, which implies 
that \eqref{eqn:3.22} must hold.  Then $x = y$ in that equation implies that $s^2(w/2) = s^2(w/4)^{\bar \tau}$, hence $s = s(w/2) = \pm s(w/4)^{\bar \tau} = \mp t$.  
If the minus sign holds, then
$$L(s,t) = L(-t,t) =  -t^2(t^4 + 18t^2 - 27)$$
and the resultant \eqref{eqn:3.27} for the quartic yields $R_d(X) = 2^{24} 3^{54}(X-1728)^4$, hence $-d = -4$.  Thus, the plus sign must hold and $s = t$, i.e.,
$s(w/2) = -s(w/4)^{\bar \tau}$.  If $L(s,t) \neq 0$, then we get the same conclusion from \eqref{eqn:3.20}.
\end{proof}

From Theorem \ref{thm:2.3} we have the equation
$$s(w/4) s^2(w/2) + s^2(w/4) + s^2(w/2) - 3s(w/4) = 0.$$
Hence, setting
\begin{equation}
k(x,y) = -x y^2 +x^2 + y^2 + 3x,
\label{eqn:3.28}
\end{equation}
we find that
\begin{align*}
0 = k(-s(w/4),&s(w/2)) = k(s(w/2)^{\bar \tau^{-1}}, s(w/2)) \ \Rightarrow \\
& k(s(w/2),s(w/2)^{\bar \tau}) = 0.
\end{align*}

\begin{cor}
With notation as in Theorem \ref{thm:3.12}, $k(s(w/2),s(w/2)^{\bar \tau}) = 0$ and the quantity $s(w/2)$ is a periodic point of the algebraic function $y = \mathfrak{k}(x)$ defined by
$$k(x,y) = -x y^2 +x^2 + y^2 + 3x = 0.$$
\label{cor:3.13}
\end{cor}

\noindent {\bf Remark.} By the remark at the beginning of the proof of Theorem \ref{thm:3.10}, it suffices to assume that $w = \frac{k+\sqrt{-d}}{2}$, with $k^2+d \equiv 0$ (mod $2^3$), for this corollary.

\begin{proof}
The relation $k(s(w/2),s(w/2)^{\bar \tau}) = 0$ implies that
$$k(s(w/2)^{\bar \tau^i},s(w/2)^{\bar \tau^{i+1}}) = 0.$$
Thus, if $a_i = s(w/2)^{\bar \tau^{i-1}}$ and 
$\bar \tau^n = 1$, then $a_{n+1} = a_1$ and 
$$k(a_1, a_2) = k(a_2,a_3) = \cdots = k(a_i, a_{i+1}) = \cdots = k(a_n,a_1) = 0.$$
This proves that $a_1 = s(w/2)$ is a periodic point for $\mathfrak{k}(x)$, and that its minimal period divides $n$.
\end{proof}

\begin{table}
  \centering 
  \caption{Class equations for polynomials $q(t)$} \label{tab:4}

\begin{tabular}{| c | c | c |}
\hline
  &   &    \\
$q(t)$   &   $R_d(X)$   &   $-d$  \\
  &  &  \\
\hline
 &   &    \\
$t^8 - 16t^6 + 86t^4 - 144t^2 + 81$ & $2^{24} 3^{72}(X-8000)^8$ & $-8$ \\
 &  &   \\
$t^8 - 20t^6 + 182t^4 - 180t^2 + 81$  &  $2^{48}3^{72}(X^2 - 1264000X - 681472000)^4$ & $-20$  \\
 &   &   \\
 $t^4 - 2t^2 + 9$ & $2^{24}3^{36}(X-8000)^4$ & $-8$  \\
 &  &    \\
 $t^4 - 4t^3 + 14t^2 + 12t + 9$ & $2^{36} 3^{36}(X^2 - 52250000X + 12167000000)^2$ & $-32$  \\
 &  &   \\
 $t^4 + 4t^3 + 14t^2 - 12t + 9$ & $2^{36} 3^{36}(X^2 - 52250000X + 12167000000)^2$ & $-32$ \\
 &  &  \\
 \hline
\end{tabular}
\end{table}

\section{Other periodic points for $k(x,y)$.}

We can use the polynomial $L(x,y)$ from \eqref{eqn:3.22a} to exhibit other periodic points for $k(x,y) = 0$, as follows.

\begin{lem}
Let $-d \equiv 1$ (mod $8$) and $w = \frac{k+\sqrt{-d}}{2}$, with $k^2 + d \equiv 0$ (mod $2^3$).  If $x$ is a conjugate of $s(w/2)$, and
$$L(x,u) = 0, \ \ \textrm{for} \ u \in \Omega_{3f},$$
then $u$ is also a periodic point for $k(x,y)$.
\label{lem:4.1}
\end{lem}

\begin{proof}
First, the relation $k(s(w/2), s(w/2)^{\bar \tau}) = 0$ implies that if $\sigma \in \textrm{Gal}(\Omega_{3f}/\mathbb{Q})$ is an automorphism, then
$$0 = k(s(w/2)^\sigma, s(w/2)^{\bar \tau \sigma}) = k(s(w/2)^\sigma, s(w/2)^{\sigma \bar \tau}) \ \textrm{or} \  k(s(w/2)^\sigma, s(w/2)^{\sigma \bar \tau^{-1}}),$$
according as $\sigma \in \textrm{Gal}(\Omega_{3f}/K)$ or not.  Thus we assume $x = s(w/2)^\sigma, y = s(w/2)^{\bar \tau \sigma}$, so that $k(x,y) = 0$, and let $v$ be such that $L(y,v) = 0$.  We compute the resultant
\begin{equation*}
\textrm{Res}_y(\textrm{Res}_x(k(x,y),L(x,u)), L(y,v)) = 3^{18} (-uv^2 + u^2 + v^2 + 3u)^3 \Lambda(u,v)^2,
\end{equation*}
\begin{align*}
& \Lambda(u,v) = u^6v^6 + 9u^5v^6 + 27u^4v^6 - 270u^5v^4 + 27u^3v^6 - 378u^4v^4 + 972u^5v^2\\
& \ - 567u^3v^4 + 4131u^4v^2 - 243u^2v^4 - 729u^5 + 2187u^3v^2 + 729uv^4 - 3645u^4\\
& \  - 5103u^2v^2 + 729v^4 - 2187u^3 - 2187uv^2 + 6561u^2.
\end{align*}
We use that $x^2 = 1 + (-2c(w/2))^{3\sigma} = 1+\pi^3$, where $\pi \cong \mathfrak{p}_2^\sigma$ from Lemma \ref{lem:3.11}.   \medskip

Using Newton's expansion for $(1+z)^{1/2}$ in an unramified extension 
$\tilde \Omega = \mathbb{Q}_2(\pi)$ of the $2$-adic field $\mathbb{Q}_2$ yields that
\begin{align*}
x &= (1+\pi^3)^{1/2} = \sum_{i=0}^\infty{\binom{1/2}{i} \pi^{3i}}\\
&= 1+ \frac{1}{2} \pi^3 - \frac{1}{8} \pi^6 + \sum_{i=3}^\infty{\binom{1/2}{i}\pi^{2i-1} \pi^{i+1}}.
\end{align*}
Now $\pi \cong 2$ in $\tilde \Omega$.  Furthermore, $2^{2i-1} \binom{1/2}{i} \in \mathbb{Z}$, so $\pi^{2i-1} \binom{1/2}{i} \in \widetilde{\mathfrak{o}}$, the ring of
integers in $\tilde \Omega$, which we will call $\pi$-adic integers.  Thus, we see that
$$x = 1 + \pi ^2 \left(\frac{\pi}{2} - \frac{\pi^3}{8} \pi + O(\pi^2)\right) = 1 + \pi^2 \kappa, \ \ \kappa \in \widetilde{\mathfrak{o}}^\times.$$
In the same way, we find that $y = 1 + \pi^2 \lambda$, with $\lambda \in \widetilde{\mathfrak{o}}^\times$.
Now consider the congruence
\begin{align*}
L(1+\pi^2 \kappa, u) &= L(1+\pi^2 \kappa, 1+u_1) \equiv \kappa^3(u_1^3 + u_1^2 + u_1 + 1)\pi^6\\
& + \kappa^2(u_1^3 + u_1^2 + u_1)\pi^4 + \kappa(u_1^3 + u_1^2)\pi^2 + u_1^3 \ (\textrm{mod} \ 2).
\end{align*}
Since $L(x,u) = 0$ and $\pi^2$ divides all of the terms except for $u_1^3$, we must have $\pi \mid u_1$.  Now write $u = 1+\pi u_2$ and compute
\begin{align*}
0 = L(x,u) &= L(1+\pi^2 \kappa,1+\pi u_2) = \kappa^3\pi^9u_2^3 + 3\kappa^3\pi^8 u_2^2 + (3\kappa^2 u_2^3 + 3\kappa^3 u_2)\pi^7\\
& + (-9\kappa^2 u_2^2 + \kappa^3)\pi^6 + (3 \kappa u_2^3 - 27 \kappa^2 u_2)\pi^5 + (-27 \kappa u_2^2 + 12 \kappa^2) \pi^4\\
& + (u_2^3 - 36 \kappa u_2) \pi^3 + (12 u_2^2 + 48 \kappa) \pi^2 + 48 \pi u_2 + 64.
\end{align*}
In this equation, the term whose coefficient has the smallest $\pi$-adic valuation is $u_2^3 \pi^3$; all other terms are divisible by $\pi^4$.  Hence we must have $\pi \mid u_2$.  
This gives that $u = 1+l \pi^2$, and similarly, $v=1+k \pi^2$, where $k,l \in \widetilde{\mathfrak{o}}$. \medskip

Now we show that $\Lambda(u,v) \neq 0$.  But a Maple calculation shows that
\begin{align*}
\Lambda(u,v) &= \Lambda(1+l \pi^2, 1+ k \pi^2) = 64 + (384k - 192l)\pi^2 + (960k^2 + 15264k l + 240l^2)\pi^4\\
&  + (1280k^3 - 6768k^2 l + 46368 k l^2 - 160 l^3)\pi^6 + \pi^7 \Xi,
\end{align*}
where $\Xi \in \mathbb{Z}[k,l,\pi]$, so that $\Xi \in \widetilde{\mathfrak{o}}$.  It is easily checked that all the terms except $64$ are divisible by at least $\pi^7$.  
But $\pi^6 || 64$, so that $\Lambda(u,v) \neq 0$.  Hence, we have shown that
\begin{equation}
k(x,y) = 0 \ \wedge \ L(x,u) = L(y,v) = 0 \ \Rightarrow \ k(u,v) = 0.
\label{eqn:4.1}
\end{equation}
Let $\psi = \bar \tau$ or $\bar \tau^{-1}$ satisfy $k(x,y) = k(x,x^\psi) = 0$.  Then $L(x,u) = 0$ implies $L(x^{\psi^i},u^{\psi^i}) = 0$.  Thus, the sequence
$$k(x,x^\psi) = \dots = k(x^{\psi^i}, x^{\psi^{i+1}}) = \dots = k(x^{\psi^{n-1}}, x) = 0, \ \psi^n = 1,$$
implies by\eqref{eqn:4.1} that
$$k(u,u^\psi) = \dots = k(u^{\psi^i}, u^{\psi^{i+1}}) = \dots = k(u^{\psi^{n-1}}, u) = 0.$$
This proves the lemma.
\end{proof}

Solving for $y$ in $L(x,y) = 0$, we find that one root is given by
$$y_1 = \frac{(3x^2 - 9)(1-x^2)^{1/3} - 9(1-x^2)^{2/3} + 6x^2 - 9}{x^3}.$$
Putting $x=s=s(\tau), c=c(\tau)$ and $s^2 = 1-8c^3$ yields
$$y_1 = \frac{(3s^2 - 9)2c - 9\cdot 4c^2 + 6s^2 - 9}{s^3} = \frac{-3(4c^2+2c+1)^2}{s^3}.$$
Two other roots are given by
\begin{align*}
y_2 &= \frac{\omega (3x^2 - 9)(1-x^2)^{1/3} - 9 \omega^2(1-x^2)^{2/3} + 6x^2 - 9}{x^3}\\
&= \frac{\omega (3s^2 - 9)(2c) - 9 \omega^2 \cdot 4c^2 + 6s^2 - 9}{s^3}\\
&= \frac{-3(4\omega^2 c^2+2\omega c+1)^2}{s^3}
\end{align*}
and
$$y_3 = \frac{-3(4\omega c^2+2\omega^2 c+1)^2}{s^3}.$$
All three of these roots lie in $\Omega_{3f}$, for $\tau = w/2$, and are periodic points for $k(x,y)$, by Lemma \ref{lem:4.1}. \medskip

Another expression for the root $y_1$ is
\begin{align*}
y_1 &= \frac{-3(4 c^2+2 c+1)^2}{s^3} = \frac{-3(4 c^2+2 c+1)^2}{s(1-8c^3)}\\
&= \frac{-3(4 c^2+2 c+1)}{s(1-2c)} = \frac{-3s}{(1-2c)^2}.
\end{align*}
We can verify this by noting that
\begin{equation*}
L\left(s, \frac{-3s}{(1 - 2c)^2}\right) = \frac{27s^2(8c^3 + s^2 - 1)(8c^3 - 24c^2 - s^2 + 24c - 7)}{(1-2c)^6},
\end{equation*}
which is equal to $0$, identically in $\tau$. \medskip

Since $(X,Y) = (-2c(\tau),s(\tau))$ is a parametrization of the elliptic curve $E: Y^2 = X^3+1$, on which $(-1,0)$ is a point of order 2, and
\begin{equation}
(X,Y) + (-1,0) = \left(\frac{-X+2}{X+1},\frac{-3Y}{(X+1)^2}\right) \ \textrm{on} \ E,
\label{eqn:4.2}
\end{equation}
we see that $y_1$ is exactly the $Y$-coordinate of the point
$$(-2c(\tau),s(\tau)) + (-1,0) = \left(\frac{2c(\tau)+2}{-2c(\tau)+1},\frac{-3s(\tau)}{(-2c(\tau)+1)^2}\right) \ \textrm{on} \ E.$$
The roots $y_2$ and $y_3$ are just $y_1$ evaluated at $\tau+1$, resp. $\tau+2$.  The mapping $\varphi: (X,Y) \rightarrow (X,Y) + (-1,0)$ 
in \eqref{eqn:4.2} is an automorphism of the function field of $E$ with order $2$. \medskip

\begin{thm}
Let $-d \equiv 1$ (mod $8$) and $w = \frac{k+\sqrt{-d}}{2}$, with $k^2 + d \equiv 0$ (mod $2^3$).   The conjugates over $\mathbb{Q}$ of the quantities
$$s(w/2), \ \frac{-3s(w/2)}{(1-2\omega^i c(w/2))^2}, \ \omega = \frac{-1+\sqrt{-3}}{2}, \ i = 0,1,2,$$
are periodic points for the algebraic function defined by $k(x,y) = 0$.
\label{thm:4.2}
\end{thm}

This theorem and Lemma \ref{lem:4.1} show that there is an interesting relationship between the curves $k(X,Y) = 0, Y^2 = X^3+1$, and $L(X,Y) = 0$.  Moreover, 
$L(X,Y) = 0$ is an elliptic curve which is isomorphic to $Y^2 = X^3+1$ over $\mathbb{Q}$.  The curve $k(X,Y) = 0$ is also an elliptic curve whose $j$-invariant is
$j = \frac{2^4 \cdot13^3}{3^2}$.  In addition, the latter curve is isomorphic to the curve $\textsf{E}_9$ whose equation is given before 
Theorem \ref{thm:3.3}. \medskip

Now we note the following identity.  If
$$U(x,y) = (9-x^2)y^2 - (1-x^2),$$
so that $U(s(\tau),r(\tau)) = 0$, then
\begin{align}
\notag \textrm{Res}_y(\textrm{Res}_x(&k(x,y),U(x,z)),U(y,w)) \\
\label{eqn:4.3} &= 2^{12}(-3wz^2 + w^2 + z^2 + w)^2(3wz^2 + w^2 + z^2 - w)^2\\
\label{eqn:4.4} &=2^{12}3^{-8}k(3w,3z)^2k(-3w,-3z)^2.
\end{align}
Furthermore, if
$$U_1(x,y) = (1-x^2) y^2-(9-x^2) = -y^2 U\left(x,\frac{1}{y}\right),$$
then
\begin{align}
\notag \textrm{Res}_y(\textrm{Res}_x(&k(x,y),U_1(x,z)),U_1(y,w)) \\
\label{eqn:4.5} &= 2^{12}(-wz^2 + w^2 + z^2 + 3w)^2 (wz^2 + w^2 + z^2 - 3w)^2\\
\label{eqn:4.6} &=2^{12}k(w,z)^2k(-w,-z)^2.
\end{align}

\begin{thm}
With $w$ as in Theorem \ref{thm:4.2}, $\pm 3r(w/2)$ (for some sign) and $\mp \frac{1}{r(w/2)}$ are periodic points for the algebraic function defined by $k(x,y) = 0$.  The values $\pm 3r(w/2)$ are not conjugate to $s(w/2)$ over $\mathbb{Q}$, but $\mp \frac{1}{r(w/2)}$ is conjugate to $s(w/2)$.
\label{thm:4.3}
\end{thm}

\begin{proof}
We know that $k(s(w/2), s(w/2)^{\bar \tau}) = 0$ from Corollary \ref{cor:3.13}.  Letting
$$x=s(w/2), z = r(w/2), y =  s(w/2)^{\bar \tau}, w =  r(w/2)^{\bar \tau},$$
equations \ref{eqn:4.3}, \ref{eqn:4.4} and $U(s(w/2), r(w/2)) = 0$ imply that
$$k(3r(w/2)^{\bar \tau}, 3r(w/2)) = 0 \ \ \textrm{or} \ \ k(-3r(w/2)^{\bar \tau}, -3r(w/2)) = 0.$$
These equations imply, respectively, that
$$k(3r(w/2), 3r(w/2)^{\bar \tau^{-1}}) = 0 \ \ \textrm{or} \ \ k(-3r(w/2), -3r(w/2)^{\bar \tau^{-1}}) = 0.$$
This implies that $3r(w/2)$ or $-3r(w/2)$ is a periodic point for $k(x,y)$.
The second assertion follows from the identity
\begin{equation}
x^2y^2 k\left(\frac{-3}{x}, \frac{-3}{y}\right) = 9k(x,y),
\label{eqn:4.7}
\end{equation}
which implies that $-3/\xi$ is a periodic point for $k(x,y)$ whenever $\xi$ is.  To prove that $\pm 3r(w/2)$ is not conjugate to $s(w/2)$, we use 
Proposition \ref{prop:3.4} and Corollary \ref{cor:3.5}.  By that proposition, $s(w/2) \cong 3^{1-a}$, where $0 \le a < 1/2$ and $a$ is a rational number.  Hence,
the corollary gives that $\pm 3r(w/2) \cong 3^a$.  Divisors of the form $3^{1-a}, 3^a$ can only be conjugate over $\mathbb{Q}$ if $a = 1-a$, i.e. $a = 1/2$, which is not the case. This proves the claim. \medskip

Appealing to Corollary \ref{cor:3.13} once again, $k(s(w/2),s(w/2)^{\bar \tau}) = 0$ (see \eqref{eqn:3.28}) implies that
$$s^2(w/2)^{\bar \tau} = \frac{s(s+3)}{s-1}, \ \ s = s(w/2).$$
Hence,
\begin{align*}
\frac{1}{r^2(w/2)^{\bar \tau}} &= \frac{s^2(w/2)^{\bar \tau}-9}{s^2(w/2)^{\bar \tau} -1}\\
&= \frac{\frac{s(s+3)}{s-1}-9}{\frac{s(s+3)}{s-1}-1} = \left(\frac{s-3}{s+1}\right)^2.
\end{align*}
This implies that
$$\frac{\pm 1}{r(w/2)^{\bar \tau}} = \frac{-s+3}{s+1} = s\left(\frac{w}{2}+\frac{1}{2}\right),$$
and therefore
$$\frac{\pm 1}{r(w/2)} = s\left(\frac{w}{2}+\frac{1}{2}\right)^{\bar \tau^{-1}}.$$
In particular, by Proposition \ref{prop:2.5}(b), $\pm 1/r(w/2)$ is a conjugate of $s(w/2)$, for some choice of the sign.  Since $s(w/2)$ generates 
$\Omega_{3f}$ over $K = \mathbb{Q}(\sqrt{-d})$, so does $r(w/2)$.  This completes the proof of the final statement of Theorem \ref{thm:3.10}.
\end{proof}

\noindent {\bf Remark.} This proof shows that the sign $\varepsilon = \pm 1$ for which $3\varepsilon r(w/2)$ is a periodic point for $k(x,y)$ is given by
$$\varepsilon = -r(w/2) s\left(\frac{w}{2}+\frac{1}{2}\right)^{\bar \tau^{-1}}.$$
We will see later in the proof of Theorem \ref{thm:5.4} (and the remark that follows)
that the sign $\varepsilon$, for which $3\varepsilon r(w/2)$ is a periodic point, is unique.

\subsection{The polynomial $\mathcal{Q}_d(y)$.}

Let $\mathcal{Q}(x,y)$ be the polynomial given by
\begin{align}
\label{eqn:4.8} \mathcal{Q}(x,y) &= (1-x^2)^2(y - 1)(y + 3)^3 - (1+ x^2)^2(y + 1)(y - 3)^3\\
\notag &= 4\left(4x^4y^3 + (-y^4 - 18y^2 + 27)x^2 + 4y^3\right)\\
\notag &= (1-x^2)^2 (y+1)(y-3)^3 \left(A(y)-\left(\frac{1+x^2}{1-x^2}\right)^2\right).
\end{align}
This polynomial arises from combining Propositions \ref{prop:3.7} and \ref{prop:3.8}.  Let $b_d(x)$ be the 
minimal polynomial of the quantities $\pi$ and $\xi$ introduced before Proposition \ref{prop:3.7}.  (See \cite{mdeu}.) The resultant 
\begin{equation}
\mathcal{Q}_d(y) = \textrm{Res}_x(\mathcal{Q}(x,y),b_d(x))
\label{eqn:4.9}
\end{equation}
will play the role of the polynomial $C_d(x)$ in \cite{akm}, and gives us a hint about where the other factors are coming from in Tables \ref{tab:2} and \ref{tab:3}, whose roots are not conjugate to $s(w/2)$ or $\pm 3 r(w/2)$.  We will show that the roots of $\mathcal{Q}_d(y)$ are all periodic points for $k(x,y) = 0$; and conversely, every periodic point for $k(x,y) = 0$ is a root of $\mathcal{Q}_d(y)$, for some discriminant $-d$.  \medskip

Note that $\mathcal{Q}(x,y)$ satisfies the identity
\begin{equation}
(x - 1)^4(y + 1)^4\mathcal{Q}\left(\frac{x + 1}{x - 1}, \frac{-y + 3}{y + 1}\right) = 2^6 \mathcal{Q}(x,y),
\label{eqn:4.10}
\end{equation}
where
\begin{equation}
(x-1)^{2h(-d)} b_d\left(\frac{x+1}{x-1}\right) = 2^{h(-d)}b_d(x)
\label{eqn:4.11}
\end{equation}
from \cite{mdeu}.

\begin{lem}
\begin{enumerate}[(a)]

\item $u(\tau)$ and $s(\tau)$ satisfy the identity
\begin{equation}
u^8(3\tau) = \frac{(s(\tau)+3)(s(\tau)-1)^3}{(s(\tau)-3)(s(\tau)+1)^3}.
\label{eqn:4.13}
\end{equation}

\item The functions $x=s(\tau), y = s(3\tau)$ satisfy the identity
\begin{equation}
\mathcal{C}(x,y) = y^3  -y^2 x^3 +  6y^2 x + 3yx^2 - 9x = 0.
\label{eqn:4.12}
\end{equation}
We also have that $\mathcal{C}(3r(3\tau),3r(\tau)) = 0$.

\end{enumerate}
\label{lem:4.4}
\end{lem}

\begin{proof}
First, part (a) follows from
\begin{align*}
F\left(\frac{(s(\tau)+3)(s(\tau)-1)^3}{(s(\tau)-3)(s(\tau)+1)^3}\right) &= \frac{(s^2 - 3)^3(s^6 - 9s^4 + 3s^2 - 3)^3}{s^4(s^2 - 1)^3 (s^2 - 9)}\\
&= j(3\tau) = F(u^8(3\tau)),
\end{align*}
using \eqref{eqn:3.16} and similar arguments as in the proof of Proposition \ref{prop:3.7}.  Now put $3\tau$ in for $\tau$ in \eqref{eqn:3.13}.  Setting
$x=s(\tau), y = s(3\tau)$ gives that
\begin{equation*}
u^8(3\tau) = \frac{(y-1)(y+3)^3}{(y+1)(y-3)^3} = \frac{(x-1)^3(x+3)}{(x+1)^3(x-3)}.
\end{equation*}
Subtracting the two expressions in $x$ and $y$ yields
$$\frac{(x-1)^3(x+3)}{(x+1)^3(x-3)} -\frac{(y-1)(y+3)^3}{(y+1)(y-3)^3} = \frac{16(xy + 3)(y^3 - x^3y^2 + 6xy^2 + 3x^2y - 9x)}{(y + 1)(y - 3)^3(x - 3)(x + 1)^3}.$$
Now $xy = s(\tau) s(3\tau) = 1 + O(q)$, as $q \rightarrow 0$, which implies that \eqref{eqn:4.12} holds.  The second assertion in (b) follows easily from 
Lemma \ref{lem:2.2}.
\end{proof}

\begin{cor} 
\begin{enumerate}[(a)]

\item We have the identities
\begin{align}
\notag u^2(\tau) u^2(3\tau) &= \frac{(s(\tau)+3)(s(\tau)-1)}{(s(\tau)-3)(s(\tau)+1)},\\
\label{eqn:4.14} \frac{u^4(3\tau)}{u^4(\tau)} &= \frac{(s(\tau)-1)(s(\tau)-3)}{(s(\tau)+1)(s(\tau)+3)}.
\end{align}

\item The functions $x=u(\tau), y = u(3\tau)$ satisfy the identity
\begin{equation*}
x^4+2x^3y^3 = y^4+2xy.
\end{equation*}
\end{enumerate}
\label{cor:4.5}
\end{cor}

\begin{proof}
For part (a), both identities are immediate from part (a) of Lemma \ref{lem:4.4} and Proposition \ref{prop:3.7}.  Note that both sides of the 
first identity are asymptotic to $4q$, as $q \rightarrow 0$, while both sides of the second are asymptotic to $q$.
To prove (b), it will suffice to prove that
$$(x^4-y^4)^2 = (-2x^3 y^3+2xy)^2,$$
or
$$x^8+6x^4 y^4 + y^8 = 4x^6 y^6 + 4x^2y^2.$$
Proposition \ref{prop:3.7}, part (a) and \eqref{eqn:4.13} imply that the last equation is equivalent to
\begin{align*}
\frac{(s-1)(s+3)^3}{(s+1)(s-3)^3} &+ 6 \frac{(s-1)^2(s+3)^2}{(s+1)^2(s-3)^2} + \frac{(s-1)^3(s+3)}{(s+1)^3(s-3)}\\
&= 4\frac{(s-1)^3(s+3)^3}{(s+1)^3(s-3)^3} +4 \frac{(s-1)(s+3)}{(s+1)(s-3)}, \ \ s = s(\tau).
\end{align*}
Both sides of this equation are equal to
$$\frac{8(s-1)(s+3)(s^4-2s^2+9)}{(s+1)^3 (s-3)^3}.$$
Thus,
$$x^4-y^4 = \pm(-2x^3 y^3+2xy).$$
Since
$$x^4-y^4 =4q^{1/2} + O(q^{3/2}) = -2x^3y^3 + 2xy, \ \ \textrm{as} \ q \rightarrow 0,$$
the plus sign must hold.  This proves part (b).
\end{proof}

It is easy to check that
\begin{equation}
x^3 y^3 \mathcal{C}\left(\frac{-3}{x},\frac{-3}{y}\right) = -3^3 \mathcal{C}(y,x).
\label{eqn:4.15}
\end{equation}

\begin{lem} 
Assume $-d \equiv 1$ (mod $8$), and let $w = \frac{a+\sqrt{-d}}{2}$, with $a^2+d \equiv 0$ (mod $24$).
\begin{enumerate}[(i)]

\item If $\left(\frac{-d}{3}\right) = +1$ and $\frak{p}_3 \mid w$, then $s(w/6) \cong s\left(\frac{w}{6}+\frac{1}{2}\right) \cong \frak{p}_3'$.

\item If $\left(\frac{d_K}{3}\right) = 0$ and $\frak{p}_3 \mid w$, then $s(w/6) \cong s\left(\frac{w}{6}+\frac{1}{2}\right) \cong 3^{1/2} \cong \frak{p}_3$.

\item If $\left(\frac{d_K}{3}\right) = 0$ or $1$ and $3 \nmid f$, then $s(w/6), r(w/6) \in \Omega_f$, where $-d = d_K f^2$.
\end{enumerate}
\label{lem:4.6}
\end{lem}

\begin{proof}
As in the proof of Proposition \ref{prop:3.4}, we write
$$s(w/6) = \left(\frac{\eta(w/6)}{\eta(w/2)}\right)^4 \left(\frac{\eta(w)}{\eta(w/3)}\right)^2,$$
so that
\begin{align*}
s(w/6)^{12} &= \left(\frac{\eta(w/6)}{\eta(w/2)}\right)^{48} \left(\frac{\eta(w)}{\eta(w/3)}\right)^{24}\\
&= 3^{24} \frac{\Delta(w,6)^2}{\Delta(w,2)^2} \cdot 3^{-12} \frac{\Delta(w,1)}{\Delta(w,3)}\\
&= \left(3^{12}\frac{\Delta(w,6)}{\Delta(w,2)}\right)^2 \left(3^{-12} \frac{\Delta(w,1)}{\Delta(w,3)}\right)\\
&= \frac{\varphi_P(w/2)^2}{\varphi_P(w)},
\end{align*}
where here $P= \begin{pmatrix}
   1   &  0  \\
    0  &  3
\end{pmatrix}$.  In (i), we have $\varphi_P(w/2) \cong \varphi_P(w) \cong \frak{p}_3'^{12}$, so that $s(w/6) \cong \frak{p}_3'$.
In (ii), $\varphi_P(w/2) \cong \varphi_P(w) \cong 3^6$, so $s(w/6) \cong 3^{1/2} \cong \frak{p}_3$.  The same arguments work with
$w$ replaced by $w+3$, since $w+3 = \frac{a+6+\sqrt{-d}}{2}$ and $(a+6)^2+d \equiv 0$ (mod $24$).  Hence $s\left((w+3)/6\right) = 
s\left(\frac{w}{6}+\frac{1}{2}\right) \cong \mathfrak{p}_3' \cong s(w/6)$. \medskip

Now we prove (iii).  If $3 \mid d_K$, then the cube roots of unity are contained in $\Omega_f$, since $\sqrt{-3} \in \Sigma$.  Then the argument
in the proof of Theorem \ref{thm:3.2} works, with $w/2$ replaced by $w/6$.  Now assume $\left(\frac{d_K}{3}\right) = 1$.  By Lemma \ref{lem:3.1}(c) we know that
$$r(w/6) s(w/6), \ s^3(w/6), \ r^3(w/6) \ \in \Omega_f.$$
As in the proof of Theorem \ref{thm:3.2}, assume first that $L = \Omega_f(s(w/6)) = \Omega_f(r(w/6))$ is a cubic extension of $\Omega_f$.  Since $\Omega_f$ 
does not contain the cube roots of unity, it is not a Kummer extension.  The previous argument yields that
\begin{align*}
&\lambda = 3 r(w/6) + s(w/6) \in \Omega_f,\\
&3r(w/6)s(w/6) = \lambda s(w/6) - s^2(w/6).
\end{align*}
Since $s^3(w/6) = \mu \in \Omega_f$, it follows that $x^3-\mu$ is reducible, and hence that $\rho s(w/6) \in \Omega_f$, for some cube root of unity $\rho$.
If $\rho \neq 1$, then $\rho^2 r(w/6) \in \Omega_f$. \medskip

On the other hand, if $L = \Omega_f(s(w/6))$ is only quadratic over $\Omega_f$, then $x^3 - s^3(w/6)$ is reducible and $\rho s(w/6) \in \Omega_f$
 for some cube root of unity $\rho$.  If $\rho \neq 1$, again we find that $\rho^2 r(w/6) \in \Omega_f$.  Thus, in either case we have
 \begin{align*}
 (\rho^2 r(w/6))^2 &= \rho r^2(w/6) =  \frac{\rho-\rho s^2(w/6)}{9-s^2(w/6)};\\
 9 \rho r^2(w/6) -\rho r^2(w/6) s^2(w/6) &= \rho - \rho s^2(w/6)\\
  &= \rho - \rho^2 (\rho s(w/6))^2 = \rho + (1+\rho)(\rho^2 s^2(w/6))\\
  &= \rho^2 s^2(w/6)+ (1+\rho^2 s^2(w/6)) \rho.
  \end{align*}
  Equating coefficients of $1$ gives that $9\rho r^2(w/6) = \rho^2 s^2(w/6)$, implying that $9r^2(w/6) \cong s^2(w/6)$; hence, $\frac{(9)}{\frak{p}_3'^2} \cong \frak{p}_3'^2$, or $\frak{p}_3^2 \cong \frak{p}_3'^2$, which is a contradiction.  Hence, $\rho = 1$ and $s(w/6) \in \Omega_f$ in any case.
\end{proof}

\begin{cor}
In the cases represented in the lemma, $\textsf{N}(s(w/6)) = 3^h$, where $h = h(-d) = [\Omega_f:K]$.
\label{cor:4.7}
\end{cor}

\begin{conj}
A root of an irreducible factor $q(x)$ of $\textsf{P}_n(x)$ (the polynomial whose roots are the periodic points of minimal period $n$; see Corollary \ref{cor:6.7} below) generates the field $\Omega_f$ over $K = \mathbb{Q}(\sqrt{-d})$, with $-d = d_K f^2$ and $3 \nmid f$, if and only if its constant term is $q(0) = 3^h$, where $\textrm{deg}(q(x)) = 2h$.  In that case, $h = h(-d)$.  This only happens when $\left(\frac{d_K}{3}\right) = 0$ or $1$, i.e., when $3$ has a first degree prime ideal factor in $K$.
\label{conj:2}
\end{conj}

Note that $s(w/2) \cong 3^a$, where $a > \frac{1}{2}$, by Proposition \ref{prop:3.4}, so its norm $3^{2ah}$ could not be $3^h$ in Conjecture \ref{conj:2}, if $s(w/2)$ happened to lie in the field generated by a root of $q(x)$.  Therefore, $s(w/2)$ cannot be a root of any such factor $q(x)$.  The same holds for 
$\pm 3r(w/2) \cong 3^{1-a}$, since $1-a < \frac{1}{2}$.

\begin{thm}
Assume $\left(\frac{d_K}{3}\right) = 0$ or $1$ and $3 \nmid f$.  If $w = \frac{a+\sqrt{-d}}{2}$, with $24 \mid a^2+d$, then $s\left(\frac{w}{6}\right)$ and $s\left(\frac{w}{6}+\frac{1}{2}\right)$ are roots of $\mathcal{Q}_d(x)$ which generate $\Omega_f$ over $K$; and both are periodic points for $k(x,y) = 0$.  If $\left(\frac{d_K}{3}\right) = 1$, then $s(w/6)$ and $s\left(\frac{w}{6}+\frac{1}{2}\right)$ are not conjugate over $\mathbb{Q}$.
\label{thm:4.8}
\end{thm}

\noindent {\bf Remark.}  If $w = \frac{a+\sqrt{-d}}{2}$, with $24 \mid a^2+d$ and $3 \mid f$, then $d = 9d'$, so $3 \mid a$ and $w = 3w'$, where $8 \mid N(w')$.  In this case, we know that $s(w/6) = s(w'/2)$ is a periodic point for $k(x,y)$ which generates $\Omega_{3f'} = \Omega_f$ over $\mathbb{Q}$, by Theorem \ref{thm:3.10} and Corollary \ref{cor:3.13}, since $f' = f/3$.

\begin{proof}
We start with $u = s(w/2), v = s(w/2)^{\bar \tau}, x = s(w/6), y = s(w/6)^{\bar \tau}$ and the identity
\begin{align*}
\textrm{Res}_v(&\textrm{Res}_u(k(u,v),\mathcal{C}(x,u)),\mathcal{C}(y,v)) \\
&= 3^8(-xy^2 + x^2 + y^2 + 3x)^3(x^4y^6 + x^5y^4 + x^3y^6 + x^6y^2 - 7x^4y^4\\
& \ \ - 5x^2y^6 - 3x^5y^2 - 9x^3y^4 + 3xy^6 - 3x^4y^2 + 51x^2y^4 + 21x^3y^2\\
& \ \ - 36xy^4 - 42x^2y^2 - 9x^3 + 90xy^2 + 27x^2 - 27x + 9)^2\\
&= 3^8 k(x,y)^3 \tilde k (x,y)^2.
\end{align*}
Since $k(u,v) = k(s(w/2),s(w/2)^{\bar \tau}) = 0$, Lemma \ref{lem:4.4} shows that
$$k(x,y) = k(s(w/6),s(w/6)^{\bar \tau}) = 0 \ \textrm{or} \ \tilde k(x,y) = \tilde k(s(w/6),s(w/6)^{\bar \tau}) = 0.$$
Suppose that $\tilde k(s(w/6),s(w/6)^{\bar \tau}) = 0$.  By Lemma \ref{lem:4.6}, we know that $s(w/6) \cong s(w/6)^{\bar \tau} \cong \frak{p}_3'$ or $\frak{p}_3$. 
It is easily checked that $\tilde k(x,y)$ is equal to
\begin{align*}
\tilde k(x,y) &= A + B + C +9,\\
A &=x^4y^6 + x^5y^4 + x^3y^6 + x^6y^2 - 7x^4y^4- 5x^2y^6 \equiv 0 \ (\textrm{mod} \ \frak{p}_3'^8 \ \textrm{or} \ \frak{p}_3^8);\\
B &= - 3x^5y^2 - 9x^3y^4 + 3xy^6 - 3x^4y^2 + 51x^2y^4 + 21x^3y^2 - 36xy^4 - 42x^2y^2\\
&\equiv 0 \ (\textrm{mod} \ 3\frak{p}_3'^4 \ \textrm{or} \ 3\frak{p}_3^4);\\
C &= - 9x^3 + 90xy^2+ 27x^2 -27x \equiv 0 \ (\textrm{mod} \ 9\frak{p}_3'^2 \ \textrm{or} \ 9\frak{p}_3^2).
\end{align*}
However, $9 \cong \frak{p}_3^2 \frak{p}_3'^2$ or $\frak{p}_3^4$.  Since $A+B+C \equiv 0$ (mod $\frak{p}_3'^4)$ in case $\left(\frac{d_K}{3}\right) = 1$, 
and is equivalent to $0$ (mod $\frak{p}_3^6$) in case $\left(\frac{d_K}{3}\right) = 0$, it follows that $\tilde k(x,y)$ cannot be 0.  Hence, $k(x,y) = 0$, which means
that $k(s(w/6),s(w/6)^{\bar \tau}) = 0$, and therefore $s(w/6)$ is a periodic point for $k(x,y)$.  From \eqref{eqn:2.9b} and Proposition \ref{prop:2.5}(a) it follows 
that $s(w/6+1/2)$ is also a periodic point for $k(x,y)$.  \medskip

As in the proof of Theorem \ref{thm:3.10} we can assume $32 \mid a^2+d$, where $l$ is chosen so that $3 \mid l$ and $s((w+2l)/6) = s(w/6)$.
Thus, putting $\tau = w/6$ in \eqref{eqn:4.13} gives that
\begin{align}
\notag \left(\frac{1+\xi^2}{1-\xi^2}\right)^2 &= u^8(w/2) = \frac{(s(w/6)+3)(s(w/6)-1)^3}{(s(w/6)-3)(s(w/6)+1)^3}\\
&=  \frac{(\frac{-3}{s(w/6)}-1)(\frac{-3}{s(w/6)}+3)^3}{(\frac{-3}{s(w/6)}+1)(\frac{-3}{s(w/6)}-3)^3} = A\left(\frac{-3}{s(w/6)}\right).
\label{eqn:4.16}
\end{align}
This formula implies that $\frac{-3}{s(w/6)}$ is a root of $\mathcal{Q}_d(x)$, and that $\Omega_f \subseteq \mathbb{Q}(s(w/6))$ (Proposition \ref{prop:3.8}).
Thus, Lemma \ref{lem:4.6}(iii) gives that $\mathbb{Q}(s(w/6)) = \Omega_f$ over $K = \mathbb{Q}(\sqrt{-d})$.  Now Proposition \ref{prop:2.5}(a) implies that
$\mathbb{Q}\left(s\left(\frac{w}{6} + \frac{1}{2}\right)\right) = \Omega_f$, as well. \medskip

It remains to show that $s(w/6)$ is a root of $\mathcal{Q}_d(x)$.  To do this we use the rational functions
$$F(x) = \frac{16(x^2+14x+1)^3}{x(x-1)^4}, \ \ A(x) = \frac{(x-1)(x+3)^3}{(x+1)(x-3)^3},$$
from the proofs of Proposition \ref{prop:3.7} and Lemma \ref{lem:4.4}, for which we know that
\begin{align*}
F(A(s(w/6))) &= F(u^8(w/6)) = j(w/6),\\
F\left(A\left(\frac{-3}{s(w/6)}\right)\right) &= F(u^8(w/2)) = j(w/2).
\end{align*}
Here,
$$A\left(\frac{-3}{x}\right) = B(x) = \frac{(x+3)(x-1)^3}{(x-3)(x+1)^3}.$$
Thus, if $\tau_3 = \left(\frac{\Omega_f/K}{\frak{p}_3}\right)$, we have $j(w/6)^{\tau_3} = j(w/2)$ and
$$F(A(s(w/6))^{\tau_3}) = F(A(-3/s(w/6))).$$
We would like to conclude that
\begin{equation}
A(s(w/6))^{\tau_3} = A(-3/s(w/6)).
\label{eqn:4.17}
\end{equation}
For $x= u(\tau), y= u(3\tau)$ we have the identity (Corollary \ref{cor:4.5}(b))
$$x^4+2x^3y^3 = y^4 +2xy.$$
Thus, $z=u^2(w/6), v=u^2(w/2)$ satisfy
\begin{equation}
z^4 -4v^3z^3 +  6v^2z^2  - 4vz + v^4 = 0.
\label{eqn:4.18}
\end{equation}
We also know from \eqref{eqn:3.18} that $v = y^2 = u^2(w/2) = \frac{-1}{\pi^\tau}$, where $\tau =  \left(\frac{\Omega_f/K}{\frak{p}_2}\right)$ and $\pi \cong \frak{p}_2$.
Putting $z=vt$ in \eqref{eqn:4.18} and dividing by $v^4$ gives the equation
$$ t^4 -4v^2 t^3+ 6t^2 - \frac{4}{v^2}t + 1= 0$$
for $t = z/v$.  This equation clearly implies that $t$ is a unit, so that $z \cong v \cong \frac{1}{\pi}$ and $x^8 \cong y^8 \cong \frac{1}{\pi^4}$.  Now, in the equation
for
$$F(x)-F(y) = \frac{16(x-y)(xy-1)}{x(x-1)^4y(y-1)^4} G(x,y)$$
from the proof of Proposition \ref{prop:3.7}, we replace $x, y$ by $x^8 = u^8(w/6)^{\tau_3}, y^8 = u^8(w/2)$ and note that all the terms in $G(x^8,y^8)$ except the constant term
$1$ have negative $\pi$-adic valuation.  This is because
$$752 = 2^4 \cdot 47, \ 2584 = 2^3 \cdot 17 \cdot 19, \ 13348 = 2^2 \cdot 47 \cdot 71.$$
Hence, $(x^8y^8-1)G(x^8, y^8) \neq 0$, implying that $x^8 = y^8$, i.e.
\begin{equation}
x^8 = u^8(w/6)^{\tau_3} = u^8(w/2) = y^8;
\label{eqn:4.19}
\end{equation}
on noting \eqref{eqn:4.16}, this yields \eqref{eqn:4.17}.  This implies further that
$$A(s(w/6)) = A(-3/s(w/6))^{\tau_3^{-1}} = \left(\frac{1+\xi^2}{1-\xi^2}\right)^{2\tau_3^{-1}}$$
from \eqref{eqn:4.16}.  Since $\xi^{\tau_3^{-1}}$ is a root of $b_d(x)$, this proves by \eqref{eqn:4.8} and \eqref{eqn:4.9} that $s(w/6)$ 
is a root of $\mathcal{Q}_d(y)$.  Then \eqref{eqn:4.10} and \eqref{eqn:4.11} show that 
$\frac{-s(w/6)+3}{s(w/6)+1} = s\left(\frac{w}{6}+\frac{1}{2}\right)$ is also a root of $\mathcal{Q}_d(x)$. \medskip

To prove the final assertion, assume $\left(\frac{d_K}{3}\right) = 1$, so that $s(w/6) \cong  s\left(\frac{w}{6}+\frac{1}{2}\right) \cong \mathfrak{p}_3'$, where 
$\mathfrak{p}_3 \mid w$.  We use the polynomial $\mathcal{C}(x,y)$, for which we know that
\begin{align*}
&\mathcal{C}(s(w/6), s(w/2)) = 0,\\
&\mathcal{C}\left(s\left(\frac{w}{6}+\frac{1}{2} \right), s\left(\frac{3w}{6}+\frac{3}{2}\right)\right) = 0.
\end{align*}
(Lemma \ref{lem:4.4}(b).) Since the exponents of $q$ in the $q$-expansion of $s(\tau)$ are integral, we know $s(\tau+1) = s(\tau)$.  Hence,
$$\mathcal{C}\left(s\left(\frac{w}{6}+\frac{1}{2} \right), s\left(\frac{w}{2}+\frac{1}{2}\right)\right) = 0,$$
or
$$\mathcal{C}\left(t\left(s\left(\frac{w}{6}\right)\right), t\left(s\left(\frac{w}{2}\right)\right)\right) = 0,$$
by Proposition \ref{prop:2.5}(a).  Now apply complex conjugation to this equation, using part (b) of the same proposition:
$$\mathcal{C}\left(\overline{t\left(s\left(\frac{w}{6}\right)\right)}, s\left(\frac{w}{2}\right)\right) = 0.$$
Letting $x_1 = s(w/6), x_2 = \overline{t\left(s\left(\frac{w}{6}\right)\right)}, y = s(w/2)$, we have that
$$\mathcal{C}(x_1,y) = \mathcal{C}(x_2,y) = 0.$$
However,
\begin{equation*}
\mathcal{C}(x_1,y)-\mathcal{C}(x_2,y) = -(x_1-x_2) [(x_1^2 + x_1 x_2 + x_2^2 - 6)y^2 -3(x_1 + x_2)y + 9].
\end{equation*}
Since $x_1 \cong \mathfrak{p}_3'$ and $x_2 = \overline{s\left(\frac{w}{6}+\frac{1}{2}\right)} \cong \mathfrak{p}_3$, we know that $x_1 \neq x_2$.  Hence, we
must have
\begin{equation}
(x_1^2 + x_1 x_2 + x_2^2 - 6)y^2 -3(x_1 + x_2)y + 9 = 0.
\label{eqn:4.20}
\end{equation}

Now assume $x_1 = s(w/6)$ and $\overline{x_2} = s\left(\frac{w}{6}+\frac{1}{2}\right)$ are conjugate over $\mathbb{Q}$.  Then $x_1$ and $x_2$ are also conjugate over $\mathbb{Q}$, and the automorphism $\sigma: x_1 \rightarrow x_2$ satisfies $\mathfrak{p}_3'^\sigma = \mathfrak{p}_3$, so that $\sigma$ does not fix 
the quadratic field $K$.  Since $\textrm{Gal}(\Omega_f/\mathbb{Q})$ is generalized dihedral, $\sigma$ must have order $2$.  Let $L$ denote the fixed field 
of $\sigma$ in $\Omega_f$.  Then $[L:\mathbb{Q}] = h(-d)$.  But now \eqref{eqn:4.20}, whose coefficients are fixed by $\sigma$ and lie in $L$, implies 
that $y = s(w/2)$ has degree $2$ over $L$, and therefore at most degree $2h(-d)$ over $\mathbb{Q}$.  But we know that $\mathbb{Q}(s(w/2)) = \Omega_{3f}$, so that 
$s(w/2)$ actually has degree $4h(-d)$ over $\mathbb{Q}$.  This contradiction shows that $x_1 = s(w/6)$ and $\overline{x_2} = s\left(\frac{w}{6}+\frac{1}{2}\right)$ 
cannot be conjugate, completing the proof.
\end{proof}

\begin{conj}
If $\left(\frac{-d}{3}\right) = 0$ or $1$, then $\pm 3 r\left(\frac{w}{6}\right)$ (for some choice of the sign) is conjugate to $s\left(\frac{w}{6}\right)$ and $\pm 3 r\left(\frac{w}{6}+\frac{1}{2}\right)$ (for some choice of the sign) is conjugate to $s\left(\frac{w}{6}+\frac{1}{2}\right)$ over $\mathbb{Q}$.
\label{conj:3}
\end{conj}

\begin{prop}
The resultant $\mathcal{Q}_d(y)$ in \eqref{eqn:4.9} factors into $3, 2$ or $1$ distinct, irreducible polynomials over $\mathbb{Q}$, according as $\left(\frac{-d}{3}\right) = 1, 0$ or $-1$.  The value $y = s(w/2)$ is always a root of $\mathcal{Q}_d(y)$, and is a root of an irreducible factor of largest degree dividing $\mathcal{Q}_d(y)$.  All the roots of $\mathcal{Q}_d(y)$ are periodic points for $k(x,y) = 0$.
\label{prop:4.9}
\end{prop}

\begin{proof}
From the proof of Proposition \ref{prop:3.8}, we know that
\begin{equation}
\frac{(s(w/2)-1)(s(w/2)+3)^3}{(s(w/2)+1)(s(w/2)-3)^3} = \left(\frac{1+\xi^2}{1-\xi^2}\right)^2.
\label{eqn:4.21}
\end{equation}
Hence,
$$\mathcal{Q}(\xi,s(w/2)) = 0 \ \wedge \ b_d(\xi) = 0 \ \Rightarrow \ \mathcal{Q}_d(s(w/2)) = 0.$$
This proves the first assertion concerning $y = s(w/2)$.  Since $\textrm{deg}(b_d(x)) = 2h(-d)$ and $\mathcal{Q}(x,y)$ has degree $4$ in $y$, it follows 
from \eqref{eqn:4.8} that
$$\textrm{deg}(\mathcal{Q}_d(y)) = 8h(-d).$$
This immediately implies that $\mathcal{Q}_d(y)$ is irreducible when $\left(\frac{-d}{3}\right) = -1$, since $s(w/2)$ generates $\Omega_{3f}$ over $\mathbb{Q}$ and 
$[\Omega_{3f}:\mathbb{Q}] = 8h(-d)$. \medskip

When $\left(\frac{-d}{3}\right) = 1$, Theorem \ref{thm:4.8} shows that $\mathcal{Q}_d(y)$ has the non-conjugate roots $s(w/6)$ and $s\left(\frac{w}{6}+\frac{1}{2}\right)$, whose minimal polynomials have degree $2h(-d)$.  Since $[\mathbb{Q}(s(w/2)):\mathbb{Q}] = 4h(-d)$, this shows that
$$\mathcal{Q}_d(y) = l^{(1)}_d(y) \cdot  l^{(2)}_d(y) \cdot q_d(y),$$
where the factors are irreducible over $\mathbb{Q}$ and $q_d(s(w/2)) = 0$, so that
$$\textrm{deg}(l^{(i)}_d(y)) = 2h(-d), \ \ \textrm{deg}(q_d(y)) = 4h(-d).$$

When $\left(\frac{-d}{3}\right) = 0$, $\mathcal{Q}_d(y) = q_d(y) l_d(y)$, where $q_d(s(w/2)) = 0$, so that
$$\textrm{deg}(q_d(y)) = 6h(-d).$$  
(See the proof of Proposition \ref{prop:3.4}.)  Hence, $\textrm{deg}(l_d(y)) = 2h(-d)$, and it is reasonable to guess that $l_d(y)$ belongs to the field $\Omega_f$ over $K$.  This follows from the equation
 following (3.17):
 $$\pi^{4\tau} = \left(\frac{1-\xi^2}{1+\xi^2}\right)^2, \ \ \tau = \left(\frac{\Omega_f/K}{\frak{p}_2}\right).$$
 We know that $\pi^4$ generates $\Omega_f$ over $\mathbb{Q}$, and that some root $\alpha$ of $l_d(y)$ satisfies
 $$\frac{(\alpha-1)(\alpha+3)^3}{(\alpha+1)(\alpha-3)^3} = \left(\frac{1+\xi^2}{1-\xi^2}\right)^{2\sigma} = \pi^{-4\tau \sigma},$$
 for some automorphism $\sigma \in \textrm{Gal}(\Omega_f/\mathbb{Q})$.  It follows from this that $\alpha$ must generate $\Omega_f$ over $\mathbb{Q}$.
Hence, $l_d(y)$ is irreducible and belongs to $\Omega_f$.  This proves the first assertion of the proposition and that $s(w/2)$ is always a root of the irreducible
factor of $\mathcal{Q}_d(y)$ of largest degree. \medskip

The last assertion is also clear when $\left(\frac{d_K}{3}\right) = \pm 1$, since each irreducible factor has a root which is a periodic point for $k(x,y)$
(Theorem \ref{thm:4.8}).  Now if $\left(\frac{d_K}{3}\right) = 0$ and $3 \nmid f$, then the periodic point $y = s(w/6)$ is a root of $l_d(y)$ by Theorem \ref{thm:4.8}.  
Assume that $3 \mid f$.  Replacing $d$ by $9d$, we will show that $\mathcal{Q}_{9d}(\pm 3 r(w/2)) = 0$, with $w = \frac{a+\sqrt{-d}}{2}$,
 and for this we must show that
\begin{equation}
y = \pm 3 r(w/2) \ \Rightarrow \ \frac{(y-1)(y+3)^3}{(y+1)(y-3)^3} = \frac{(1+x^2)^2}{(1-x^2)^2}, \ \ b_{9d}(x) = 0.
\label{eqn:4.22}
\end{equation}
We set $y = -3/z$, so that $z = \frac{\mp 1}{r(w/2)} = s(w/2)^\psi$, for some automorphism $\psi \in \textrm{Gal}(\Omega_{3f}/\mathbb{Q})$ 
(by Theorem \ref{thm:4.3}).  Then we must show that
$$\frac{(z+3)(z-1)^3}{(z-3)(z+1)^3} =  \frac{(1+x^2)^2}{(1-x^2)^2}, \ \ b_{9d}(x) = 0;$$
and for this it suffices to show that
$$\frac{(s(w/2)+3)(s(w/2)-1)^3}{(s(w/2)-3)(s(w/2)+1)^3} =  \frac{(1+x^2)^2}{(1-x^2)^2}, \ \ b_{9d}(x) = 0.$$
By the identity \eqref{eqn:4.13} this is equivalent to
$$u^8(3w/2) = \frac{(1+x^2)^2}{(1-x^2)^2}, \ \ b_{9d}(x) = 0.$$
But if $w = \frac{a+\sqrt{-d}}{2}$, then $3w/2 = \frac{3a+\sqrt{-9d}}{4}$ is an ideal basis quotient for the order $\textsf{R}_{-9d}$, 
so \eqref{eqn:3.18a} says in this situation that
$$u^8(3w/2) = \frac{(1+\xi^2)^2}{(1-\xi^2)^2},$$
for a solution of $\pi^4+\xi^4 = 1$ in $\Omega_{3f}$, where $\mathbb{Q}(\xi) = \Omega_{3f}$ and $\xi$ is a root of $b_{9d}(x)$.  This proves \eqref{eqn:4.22}.  Since $[\mathbb{Q}(r(w/2)):\mathbb{Q}] = [\Omega_{3f}:\mathbb{Q}] = 2h(-9d)$, this shows that $\pm 3r(w/2)$ is a root of $l_{9d}(y)$.
\medskip

Since $\pm 3r(w/2)$ is a periodic point for $k(x,y)$, by Theorem \ref{thm:4.3}, this completes the proof.
\end{proof}

\begin{thm}
If $\left(\frac{d_K}{3}\right) = 0$ or $1$ and $w = \frac{a+\sqrt{-d}}{2}$, with $a^2 + d \equiv 0$ mod $24$, where $-d = d_Kf^2, 3 \nmid f$, 
and $\frak{p}_2 \mid w$, then 
$$\upsilon = \frac{(s(w/6)-1)(s(w/6)-3)}{(s(w/6)+1)(s(w/6)+3)}$$
is the fourth power of a unit in the ring class field $\Omega_f$ over $K = \mathbb{Q}(\sqrt{-d})$.
\label{thm:4.10}
\end{thm}

\begin{proof}
From \cite[pp. 1978-79]{lm} we take the fact that
$$\pi = \frac{\beta}{\zeta_8^j \alpha} = \frac{\beta}{\alpha_1}, \ \ \zeta_8 = e^{2\pi i/8},$$
for some odd value of $j$, where $\alpha_1, \beta \in \Omega_f$ and $\alpha, \beta$ satisfy the equation
$$16\alpha^4 + 16\beta^4 = \alpha^4 \beta^4.$$
Moreover, $K(\sqrt{\beta})$ is an abelian extension of $K$, by \cite[Prop. 8.9, p.1990]{lm}.  From \eqref{eqn:4.19} we know that
$$u^8(w/6)^{\tau_3} = u^8(w/2), \ \ \tau_3 = \left(\frac{\Omega_f/K}{\frak{p}_3}\right),$$
and 
$$u^2(w/2) = \frac{-1}{\pi^\tau}, \ \ \tau = \left(\frac{\Omega_f/K}{\frak{p}_2}\right).$$
Corollary \ref{cor:4.5}, equation \eqref{eqn:4.14} gives that
$$\frac{u^4(w/2)}{u^4(w/6)} = \frac{(s(w/6)-1)(s(w/6)-3)}{(s(w/6)+1)(s(w/6)+3)}.$$
Thus, we must show that for
$$\frac{u(w/2)}{u(w/6)} = \varepsilon,$$
some $4$-th root of unity times $\varepsilon$ is a unit contained in $\Omega_f$.  Equation \eqref{eqn:3.18}, together with \eqref{eqn:4.19}, shows that
$$\varepsilon^8 = \left(\frac{u(w/2)}{u(w/6)} \right)^8 = \frac{1/\pi^{4\tau}}{1/\pi^{4\tau \tau_3^{-1}}} = \frac{\pi^{4\tau \tau_3^{-1}}}{\pi^{4\tau}}.$$
Now we know that $\varepsilon^2 \in \Omega_f$, since $u^2(w/6) \in \Omega_f$ by the first identity in Corollary
\ref{cor:4.5} and Lemma \ref{lem:4.6}(iii).  This implies that
$$\varepsilon^4 = \left(\frac{u(w/2)}{u(w/6)} \right)^4 = \frac{\pi^{2\tau \tau_3^{-1}}}{\pi^{2\tau}},$$
where the plus sign must hold, since $\sqrt{-1}$ is not contained in $\Omega_f$.
It will suffice to show that $\pi^{\tau_3}/\pi$ is the square of a unit in $\Omega_f$.  For this we use the above representation in terms of $\beta$ and $\alpha_1$:
\begin{equation*}
\frac{\pi^{\tau_3}}{\pi} = \frac{\beta^{\tau_3}}{\beta} \left(\frac{\alpha_1}{\alpha_1^{\tau_3}}\right).
\end{equation*}
From \cite[Prop. 8.10, pp. 1990-91]{lm} we know that $\alpha_1$ is a square in $\Omega_f$.  Thus, the term in parenthesis a square.  Furthermore, the fact that
$K(\sqrt{\beta})$ is abelian over $K$ implies that $K(\sqrt{\beta}) = K(\sqrt{\beta})^{\tau_3} = K(\sqrt{\beta^{\tau_3}})$, and hence that $\beta^{\tau_3}/\beta$ is a square in $\Omega_f$, by Kummer theory.  This shows that $\pi^{\tau_3}/\pi$ is a square in $\Omega_f$, from which we conclude that $\varepsilon^4$ is a $4$-th 
power in $\Omega_f$.  Since $\pi^{\tau \tau_3^{-1}} \cong \pi^\tau$, it is clear that $\varepsilon$ is a unit.  This proves the theorem.
\end{proof}

\noindent {\bf Example.} \medskip

\noindent The polynomial
\begin{align*}
p_{111}(x) &= x^{16} + 21x^{15} + 207x^{14} - 504x^{13} - 927x^{12} + 11826x^{11} - 22275x^{10}\\
& \ \  - 99387x^9 + 173988x^8 + 298161x^7 - 200475x^6 - 319302x^5 - 75087x^4\\
& \ \ + 122472x^3 + 150903x^2 - 45927x + 6561
\end{align*}
is a factor of $\mathcal{Q}_d(x)$, for $d = 111 = 3 \cdot 37$, whose root $s(w/6)$ generates the Hilbert class field $\Omega_1 = \Sigma$ 
of $K = \mathbb{Q}(\sqrt{-111})$, and which is a periodic point of period $8$ for $k(x,y)$.  This can be checked by using the fact that
\begin{align*}
b_{111}(x) &= x^{16} - 58x^{15} + 992x^{14} - 340x^{13} + 3786x^{12} - 5320x^{11} + 5496x^{10}\\
& \ \  - 7132x^9 + 109x^8 - 1502x^7 + 4464x^6 - 7952x^5 + 12288x^4 - 9440x^3\\
& \ \ + 5504x^2 - 896x + 256.
\end{align*}
Also, the resultant
$$\textrm{Res}_s((s+1)(s+3)x^4-(s-1)(s-3),p_{111}(s))$$
is divisible by the polynomial
\begin{align*}
q_{111}(x) &= x^{16} + x^{15} + 11x^{14} - 26x^{13} + 33x^{12} + 40x^{11} + 57x^{10} + 67x^9 + 52x^8\\
& \ \  - 67x^7 + 57x^6 - 40x^5 + 33x^4 + 26x^3 + 11x^2 - x + 1,
\end{align*}
whose roots $\varepsilon$ are units in $\Sigma$.  This polynomial is reverse reciprocal, meaning that $x^{16}q_{111}(-1/x) = q_{111}(x)$.  Hence 
$\varepsilon \rightarrow -1/\varepsilon$ is an automorphism of $\Sigma = \mathbb{Q}(\varepsilon)$.

\begin{conj}
If $\tau = w/2$, the quantity
$$\upsilon = \frac{(s(\tau)-1)(s(\tau)-3)}{(s(\tau)+1)(s(\tau)+3)}$$
is the fourth power of a unit in $\Omega_{3f}$.
\label{conj:4}
\end{conj}

\begin{prop}
If $d \neq d'$, then $(\mathcal{Q}_d(y), \mathcal{Q}_{d'}(y)) = 1$ in $\mathbb{Z}[y]$.
\label{prop:4.11}
\end{prop}

\begin{proof}
Suppose that $\alpha$ is a common root of $\mathcal{Q}_d(y)$ and $\mathcal{Q}_{d'}(y)$.  Then by \eqref{eqn:4.8},
$$A(\alpha) = \left(\frac{1+\xi^2}{1-\xi^2}\right)^2 =  \left(\frac{1+\xi'^2}{1-\xi'^2}\right)^2,$$
where
$$b_d(\xi) = 0, \ \ b_{d'}(\xi') = 0.$$
Then 
$$\frac{\xi^2+\frac{1}{\xi^2}+2}{\xi^2+\frac{1}{\xi^2}-2} = \frac{\xi'^2+\frac{1}{\xi'^2}+2}{\xi'^2+\frac{1}{\xi'^2}-2}.$$
As in the last part of the proof of Proposition \ref{prop:3.8}, this implies that $\xi^2 = \xi'^2$ and $\xi = \pm \xi'$.  By \cite[Thm. 14]{mdeu}, 
$-\xi'$ is a pre-periodic point of the same multi-valued function $\hat F(z)$ that $\xi$ is a periodic point for.  Hence, $\xi = -\xi'$ is impossible
and $\xi = \xi'$, implying that $b_d(x) = b_{d'}(x)$ and therefore $d = d'$.  This proves the proposition.
\end{proof}

\section{Results for $c(\tau)$.}

We set
$$g(x,y) = x^2+2xy^2-y,$$
for which we know that $g(c(\tau),c(2\tau)) = 0$ from \cite[Prop. 1]{akm}.

\begin{lem}
If $\alpha$ is a periodic point for $g(x,y) =0$ of period $n$, then $\omega^i \alpha$ is also a periodic point for $g(x,y) =0$, for $i = 1,2$, and its 
period divides $m = \textrm{lcm}(2,n)$.
\label{lem:5.1}
\end{lem}

\begin{proof}
Note first that
$$g(\omega x, \omega^2 y) = \omega^2 x^2 +2\omega x (\omega^2y)^2 - \omega^2 y = \omega^2 g(x,y).$$
To show $\omega \alpha$ is a periodic point, let $\alpha_0 = \alpha, \alpha_1, \dots, \alpha_{n-1}, \alpha_n = \alpha$ be a sequence for which
$$g(\alpha_i, \alpha_{i+1}) = 0, \ \ 0 \le i \le n-1.$$
Then the above identity shows that
$$g(\omega^{2^i}\alpha_i,\omega^{2^{i+1}} \alpha_{i+1}) = 0, \ \ 0 \le i \le n-1.$$
This chain of equations starts with $\omega \alpha$ and ends with $\omega^{2^n} \alpha$.  If $n$ is even, then $\omega^{2^n} = \omega$ and $\omega \alpha$
has period $n$.  If $n$ is odd, $\omega^{2^n} = \omega^2$, and there is a similar chain of $n$ equations beginning with $\omega^2 \alpha$ and ending with
$\omega^{2^{n+1}} \alpha = \omega \alpha$.  In this case, $\omega \alpha$ has period $2n$.
\end{proof}

\begin{thm}
If $-d \equiv 1$ (mod $8$) and $w = \frac{k+\sqrt{-d}}{2}$, with $k^2+d \equiv 0$ (mod $2^3$), then
$c(w/2)$ is a periodic point for the algebraic function defined by the relation
\begin{equation}
g(x,y) = x^2+2xy^2-y = 0.
\label{eqn:5.1}
\end{equation}
Moreover, $g(c(w/2), c(w/2)^{\bar \tau}) = 0$.
\label{thm:5.2}
\end{thm}

\begin{proof}
First we calculate the resultant
\begin{align*}
\textrm{Res}_t(\textrm{Res}_s(k(s,t),&1-s^2-8x^3), 1-t^2-8y^3) =\\
& (x^2 + 2xy^2 - y)^2(4x^2y^4 - 2x^3y^2 + x^4 + 2xy^3 + x^2y + y^2)^2.
\end{align*}
Since $k(s(w/2),s(w/2)^{\bar \tau})= 0$ and $c(w/2) \in \Omega_{3f}$ (Theorem \ref{thm:3.3}), this calculation shows that 
$(x,y) = (c(w/2), c(w/2)^{\bar \tau})$ satisfies either \eqref{eqn:5.1} or
$$T(x,y) = 4x^2y^4 - 2x^3y^2 + x^4 + 2xy^3 + x^2y + y^2 = 0.$$
Using Lemma \ref{lem:3.11}, which shows that $x, y = x^{\bar \tau}$ are integral for $\mathfrak{p}_2$, and the fact that
$\bar \tau$ is the Frobenius automorphism for $\mathfrak{p}_2$ on $\Omega_{3f}$, we have that
$$T(x,y) \equiv T(x,x^2) = 3x^4 \equiv x^4 \ (\textrm{mod} \ \mathfrak{p}_2).$$
Therefore, Lemma \ref{lem:3.11} implies that $T(x,y) \neq 0$.  This implies $g(x,y) = 0$ and proves the theorem.
\end{proof}

\begin{cor}
With notation and assumptions as in the theorem, the relation
$$c(w/4)^{\bar \tau} = c(w/2)$$
holds, where $\bar \tau = \left(\frac{\Omega_{3f}/K}{\mathfrak{p}_2}\right)$.
\label{cor:5.3}
\end{cor}

\begin{proof}
Theorem \ref{thm:3.12} implies that
$$1-8c^3(w/4)^{\bar \tau} = s^2(w/4)^{\bar \tau} = s^2(w/2) = 1-8c^3(w/2);$$
hence, that $c^3(w/4)^{\bar \tau} = c^3(w/2)$.  This yields that $c(w/4)^{\bar \tau} = \omega c(w/2)$, for some 
cube root of unity $\omega$. If $\omega \neq 1$, then using $\omega^{\bar \tau} = \omega^2$ gives
\begin{align*}
0 = & \ g(c(w/4),c(w/2)) = g(\omega^2 c(w/2)^{\bar \tau^{-1}}, c(w/2))\\
& \Rightarrow \ g(\omega c(w/2), c(w/2)^{\bar \tau}) = 0.
\end{align*}
However, if $x=c(w/2), y = c(w/2)^{\bar \tau}$, we have the resultant
$$\textrm{Res}_y(g(x,y),g(\omega x,y)) = -6(2x^3 - 1)x^3 \neq 0.$$
Since $2$ is unramified in $\Omega_{3f}$, it follows that $c(w/2)^3 \neq \frac{1}{2}$ and $\omega =1$, giving the assertion.
\end{proof}

We can now prove

\begin{thm}
With $w = \frac{a+\sqrt{-d}}{2} \in \textsf{R}_{-d}$ and $a^2+d \equiv 0$ (mod $8$), we have that $\mathbb{Q}(c(w/2)) = \Omega_{3f}$.
\label{thm:5.4}
\end{thm}

\begin{proof}
We first show that the normal closure of $L = \mathbb{Q}(c(w/2))$ is $\Omega_{3f}$.  Let $\bar L$ denote this normal closure.  Corollary \ref{cor:5.3} shows that
$\bar L$ contains $c(w/2)$ and $c(w/4)$, and therefore also $r(w/2)$, by the formula
\begin{equation*}
r(w/2) = \frac{c^2(w/2)}{c(w/4)+c^2(w/2)},
\end{equation*}
which is \eqref{eqn:2.2}.  But $r(w/2)$ generates $\Omega_{3f}$ over $\mathbb{Q}$, so that $\Omega_{3f} \subseteq \bar L$ and therefore $\bar L = \Omega_{3f}$, 
by Theorem \ref{thm:3.3}.  Now suppose that $L = \mathbb{Q}(c(w/2))$ is a proper subfield of $\Omega_{3f}$.  Then the relation
$$s^2(w/2) = 1-8c^3(w/2)$$
implies that $[\Omega_{3f}:L] = 2$, since $s(w/2)$ satisfies a quadratic equation over $L$.  Furthermore $s^2(w/2) \in L$, so that the Galois group of $\Omega_{3f}$
over $L$ is $\langle \sigma \rangle$, where $\sigma: s(w/2) \rightarrow -s(w/2)$.  Equation \eqref{eqn:3.4} implies that $r^2(w/2) \in L$ and $r(w/2)^\sigma = -r(w/2)$,
 since $r(w/2)$ generates $\Omega_{3f}$.  But the above relation for $r(w/2)$ implies that
$$-\frac{c^2(w/2)}{c(w/4)+c^2(w/2)} = \frac{c^2(w/2)}{c(w/4)^\sigma+c^2(w/2)},$$
and solving for $c(w/4)^\sigma$ gives that
$$c(w/4)^\sigma = -2c^2(w/2)-c(w/4),$$
or
$$2c(w/4)^\sigma+2c(w/4) = -4c^2(w/2).$$
Now Corollary \ref{cor:5.3} and Lemma \ref{lem:3.11} show that $2c(w/4) \cong 2c(w/2) \cong \mathfrak{p}_2$, so that the last relation implies that $\mathfrak{p}_2 \mid 2c(w/4)^\sigma$.  But this gives that $\mathfrak{p}_2^\sigma = \mathfrak{p}_2$.  Hence, $\sigma$ fixes $K$, from which we see that $K \subset L$ and $\sigma \in \textrm{Gal}(\Omega_{3f}/K)$.  Now the structure of 
 $\Gamma = \textrm{Gal}(\Omega_{3f}/\mathbb{Q})$ implies that $\sigma$ lies in the center of $\Gamma$, whence it follows that $L = \bar L$ is normal over $\mathbb{Q}$.  But then $L = \Omega_{3f}$, contradicting that $[\Omega_{3f}:L] = 2$.  Alternatively, $K \subset L \Rightarrow \Omega_{3f} = K(j(3w/2)) \subseteq L$ 
 from \cite[Prop. 5, (3.6)]{akm}.  Hence, we must have $L = \Omega_{3f}$.
\end{proof}

\noindent {\bf Remark.} The same proof shows that $L = \mathbb{Q}(s^2(w/2))$ is not a proper subfield of $\Omega_{3f} = \mathbb{Q}(s(w/2))$.  
Hence, $s(w/2)$ is not conjugate to $s(w/4)$.  If it were, then by Theorem \ref{thm:3.12} $s(w/2)$ would be conjugate to $-s(w/2)$ and 
$L \subsetneq \Omega_{3f}$ would be the fixed field of the automorphism $\sigma: s(w/2) \rightarrow -s(w/2)$.  \medskip

From \eqref{eqn:5.1} we know that $(x,y) = (2c(w/2),2c(w/2)^{\bar \tau})$ satisfies
$$g_1(x,y) = xy^2 - 2y + x^2 = 0,$$
so that
$$2c(w/2)^{\bar \tau} = \frac{1 \pm \sqrt{1-x^3}}{x} = \frac{1}{x} \pm \frac{1}{x} \sum_{n=0}^\infty{\binom{1/2}{n} (-1)^n x^{3n}}, \ \ x = 2c(w/2),$$
where the series is convergent in the completion $\Omega_{3f,\frak{P}}$, and $\frak{P}$ is a prime divisor of $\frak{p}_2$ in $\Omega_{3f}$.  This uses the 
fact that $x \cong \frak{p}_2 \cong 2$ in $\Omega_{3f,\frak{P}}$ and
$$\binom{1/2}{n} (-1)^n x^{2n-1} \cong \binom{1/2}{n} 2^{2n-1} = (-1)^{n-1} C_{n-1} \in \mathbb{Z},$$
where $C_{n-1}$ is the Catalan number.  The minus
sign must hold here, since
\begin{equation}
\frac{1}{x} + \frac{1}{x} + \frac{1}{x} \sum_{n=1}^\infty{\binom{1/2}{n} (-1)^n x^{3n}} = \frac{2}{x} + \sum_{n=1}^\infty{\binom{1/2}{n} (-1)^n x^{2n-1} x^{n}},
\label{eqn:5.2} 
\end{equation}
where the last series is divisible by $x \cong \Pi$, a prime element for $\frak{P}$ in $\Omega_{3f,\frak{P}}$.  However, $2/\Pi$ is a unit in this
completion, so that \eqref{eqn:5.2} represents a unit.  Now, $\bar \tau$ fixes $\frak{p}_2$ and $2c(w/2)^{\bar \tau} \cong \frak{p}_2$.  This shows that
the minus sign must hold above, i.e.,
\begin{equation}
2c(w/2)^{\bar \tau} = \frac{1 - \sqrt{1-x^3}}{x} = \frac{1}{x} - \frac{1}{x} \sum_{n=0}^\infty{\binom{1/2}{n} (-1)^n x^{3n}}, \ \ x = 2c(w/2).
\label{eqn:5.3}
\end{equation}
This expression gives the unique solution in $\Omega_{3f,\frak{P}}$ of the equation $g_1(x,y) = 0$, for which $y \cong x \cong 2$.  This also holds in any 
finite extension of $\mathbb{Q}_2$.

\begin{prop}
The minimal period of $2c(w/2)$ with respect to the algebraic function defined by $g_1(x,y) = 0$ is the order of the Frobenius automorphism
$$\bar \tau = \left(\frac{\Omega_{3f}/K}{\frak{p}_2}\right)$$
in $\textrm{Gal}(\Omega_{3f}/K)$.
\label{prop:5.5}
\end{prop}

\begin{proof}
Suppose that $\alpha_i, 0 \le i \le n$ satisfy $\alpha_0 = \alpha_n = 2c(w/2)$ and
\begin{equation}
g_1(\alpha_i,\alpha_{i+1}) = 0, \ \ 0 \le i \le n-1.
\label{eqn:5.4}
\end{equation}
Then each $\alpha_i$ is an algebraic integer, since $\alpha_i \mid \frak{p}_2$ by induction, using that $\alpha_{i+1} \mid \alpha_i$ from the equation
$$\alpha_{i+1}^2-\frac{2}{\alpha_i} \alpha_{i+1} + \alpha_i = 0.$$
If $w_\frak{P}(\alpha_i)$ were strictly smaller than the ramification index of a prime divisor $\frak{P}$ of $\frak{p}_2$ in some extension of $K$, then this equation shows that the same would be true of $w_\frak{P}(\alpha_{i+1})$; which cannot happen since $\alpha_n = 2c(w/2) \cong \frak{p}_2$.  Hence, $\alpha_i \cong \frak{p}_2$ for each $i$.  Considering \eqref{eqn:5.4} in a finite extension of the completion
$\Omega_{3f,\frak{P}}$ and using the above observation about \eqref{eqn:5.3} shows that $\alpha_1 = 2c(w/2)^{\bar \tau}$, from which we find inductively that 
$\alpha_{i} = \alpha_{i-1}^{\bar \tau} = 2c(w/2)^{\bar \tau^i}$.  Thus $2c(w/2) = 2c(w/2)^{\bar \tau^n}$, which implies that $n$ is the order of $\bar \tau$ in $\textrm{Gal}(\Omega_{3f}/\mathbb{Q})$, since $\mathbb{Q}(c(w/2)) = \Omega_{3f}$.
\end{proof}

\begin{cor}
The minimal period of $c(w/2)$ as a periodic point for $g(x,y) = 0$ is always even.
\label{cor:5.6}
\end{cor}

\begin{proof}
By the Artin map, the minimal period of $c(w/2)$ is equal to the order of $\frak{p}_2$ in the ideal class group 
of $\textsf{R}_{-9d}$, and this, in turn, is equal to the smallest integer $n \ge 1$, for which the equation
$$4 \cdot 2^n = x^2+9dy^2, \ \ 3 \nmid x,$$
has a primitive solution $(x,y) \in \mathbb{Z}^2$.  Considering this equation modulo $3$ shows that the left side is $2^n \equiv 1$ mod $3$, from
which it is clear that $n$ must be even. 
\end{proof}

\begin{prop}
Assume $\left(\frac{d_K}{3}\right) = 0$ or $1$.  Let $w = \frac{a+\sqrt{-d}}{2}$, where $-d = d_K f^2 \equiv 1$ (mod $8$), $24 \mid a^2+d$.  Setting $\gamma_2(\tau) = j(\tau)^{1/3}$, we have that $c(w/6) \gamma_2(w/3)^{-1} \in \Omega_f$.  Further, $c(w/6)^3 \in \Omega_f$, and $\omega^i c(w/6) \in \Omega_f$ for some $i$ if and only if $j(w/3)$ is a cube in $\Omega_f$.  In any case, $c(w/6) \in \Omega_{3f}$.
\label{prop:5.7}
\end{prop}

\noindent {\bf Remark.} By the results of \cite[Thm. 1, Eq. (2.10)]{mring}, $j(w/3)$ is a cube in $\Omega_f$ whenever $-d \equiv 1$ (mod $3$). See Tables 
\ref{tab:4}, \ref{tab:5}. \medskip

\begin{proof}
Replace $w$ in \eqref{eqn:3.5} and \eqref{eqn:3.6} by $w/3$, yielding 
\begin{align*}
c^8(w/6) &= \left(\frac{\eta(w/6)}{\eta(w/3)}\right)^8 \left(\frac{\eta(w)}{\eta(w/2)}\right)^{24}\\
&= \left(\frac{\eta(w/6)}{\eta(w/3)}\right)^8 \gamma_2(w/3)\\
& \ \ \times \left(\frac{\eta(w)}{\eta(w/2)}\right)^{24} \gamma_2(w)^{-3} \times \frac{\gamma_2(w)^3}{\gamma_2(w/3)}.
\end{align*}
In this case, the arguments of $\eta(\tau)$ in this formula are basis quotients of ideals in $\textsf{R}_{-d}$.  Furthermore,
\begin{equation*}
c(w/6) = \frac{c^9(w/6)}{c^8(w/6)} = \gamma_2(w/3) \frac{c^9(w/6)}{\gamma_2(w/3) c^8(w/6)}.
\end{equation*}
From the proof of Proposition \ref{prop:4.9} it follows that $\frac{1-s^2(w/6)}{9-s^2(w/6)} = r^2(w/6)$ and $s^2(w/6)$ lie in $\Omega_f$ (whether $3 \mid f$ or not). Since $c^3(w/6) = (1-s^2(w/6))/8$ and $\gamma_2(w/3) c^8(w/6)$ both lie in $\Omega_f$ (the latter by Schertz \cite[p. 159]{sch}), it follows that $c(w/6) \gamma_2(w/3)^{-1} \in \Omega_f$.  The assertions concerning $c^3(w/6)$ 
and $\omega^i c(w/6)$ are immediate from $\gamma_2^3(w/3) = j(w/3)$.  To show that $c(w/6) \in \Omega_{3f}$, recall from the proof of 
Theorem \ref{thm:3.3} that $\gamma_2(w) \in \Omega_{3f}$.  But $j(w)$ and $j(w/3)$ are conjugate over $\mathbb{Q}$, whence it follows that
$\gamma_2(w)$ and some $\omega^j \gamma_2(w/3)$ are conjugate over $\mathbb{Q}$.  Since $\Omega_{3f}$ is normal over $\mathbb{Q}$ and $\omega \in \Omega_{3f}$, we get that $\gamma_2(w/3) \in \Omega_{3f}$.  This implies the final assertion.
\end{proof}

\begin{prop}
For any discriminant $-d = d_K f^2 \equiv 1$ (mod $8$), and $w = \frac{a+\sqrt{-d}}{2} \in \textsf{R}_{-d}$, with $a^2+d \equiv 0$ (mod $8$), 
the quantity $c(w/6)$ is a periodic point for $g(x,y) = 0$.  The same holds for $c(w/6+1/2)$.
\label{prop:5.8}
\end{prop}

\begin{proof}
We use the fact that $c(w/2)$ is a periodic point for $g(x,y)$.  Let $N$ be the normal closure of $L = \mathbb{Q}(c(w/2),c(w/6)) = \Omega_{3f}(c(w/6))$ over $K$.  
($N= \Omega_{3f}$ if $\left(\frac{d_K}{3}\right) = 0, 1$, by Theorem \ref{thm:5.4} and Proposition \ref{prop:5.7}.) First, if $x=c(w/6)$, $y = c(w/2)$, we have the
modular identity
$$g_3(x,y) = x^3(4y^2 + 2y + 1) - y(y^2 - y + 1) = 0$$
from \cite[p. 1716]{akm}.  Furthermore, writing $x' = 1/x, y' = 1/y$, we have
\begin{equation*}
0 = -x'^3 y'^3 g_3\left(\frac{1}{x'},\frac{1}{y'}\right) = (y'^2 - y' + 1)x'^3 - y'^3 - 2y'^2 - 4y'.
\end{equation*}
By Lemma \ref{lem:3.11}, $y' = 1/c(w/2) \cong \mathfrak{p}_2'$.  The last equation shows that $x' = 1/c(w/6)$ is integral 
for $\mathfrak{p}_2, \mathfrak{p}_2'$ (or use \cite[Cor., p. 1719]{akm}) and that $\mathfrak{p}_2' \mid x' = 1/c(w/6)$.  It also shows that no prime divisor 
$\mathfrak{P}$ of $\mathfrak{p}_2$ in $L$ divides $x'$, since otherwise $\mathfrak{P} \mid y' \cong \mathfrak{p}_2'$.  Also, 
from the above equation and $y' \cong \mathfrak{p}_2'$ it is clear that
$$((y'^2-y'+1)(y'^3 + 2y'^2 + 4y'), \mathfrak{p}_2) = 1, \ \textrm{so}$$
$$\textrm{disc}_{X}\left((y'^2 - y' + 1)X^3 - (y'^3 + 2y'^2 + 4y')\right) = -27(y'^2 - y' + 1)^2 y'^2 (y'^2 + 2y' + 4)^2$$
is relatively prime to $\mathfrak{p}_2$.  From this it follows that the prime divisors of $\mathfrak{p}_2$ do not ramify in the extension 
$L = \mathbb{Q}(c(w/2),c(w/6))$ over $\Omega_{3f}$ or in $N/K$.\medskip

Now consider the double resultant
\begin{align*}
\textrm{Res}_y(\textrm{Res}_x(g(x,y),&g_3(z,x)),g_3(w,y)) = 729(z^2 + 2zw^2 - w)^3\\
& \ \ \ \times (4w^4z^2 - 2w^2z^3 + 2w^3z + z^4 + wz^2 + w^2)^3\\
& = 729g(z,w)^3 l(z,w)^3.
\end{align*}
We take $x = c(w/2), y = c(w/2)^{\bar \tau}, z = c(w/6), w = c(w/6)^{\bar \tau}$, where we denote by $\bar \tau = \left[\frac{N/K}{\overline{\mathfrak{P}}}\right]$ any extension of the automorphism $\bar \tau = \left(\frac{\Omega_{3f}/K}{\mathfrak{p}_2}\right)$ of $\Omega_{3f}/K$ to $N/K$, with a prime divisor 
$\overline{\mathfrak{P}} \mid \mathfrak{p}_2$ in $N$.  This makes sense, since the prime divisors of $\mathfrak{p}_2$ don't ramify in $N/K$.  \medskip

Then we have
$$g(x,y) = g_3(z,x) = g_3(w,y) = 0,$$
which implies that the above resultant is $0$.  We wish to show $g(z,w) = 0$.  To do this we show $l(z, w) \neq 0$.  Set $z' = 1/z, w' = 1/w$.  If $l(z,w) = 0$, 
then
\begin{align*}
0 = z'^4 w'^4 l\left(\frac{1}{z'},\frac{1}{w'}\right) &= w'^4 + w'^3 z'^2 + (z'^4 - 2z')w'^2 + 2w'z'^3 + 4z'^2\\
& \equiv w'^2(z'^4+w'z'^2+w'^2) \ (\textrm{mod} \ \mathfrak{p}_2).
\end{align*}
Furthermore, $w' = z'^{\bar \tau} \equiv z'^2$ (mod $\overline{\mathfrak{P}}$), since $z'$ is integral for $\mathfrak{p}_2$.  From this we obtain that
$$0 \equiv z'^4 (z'^4+z'^4+z'^4) \equiv z'^8 \ (\textrm{mod} \ \overline{\mathfrak{P}}).$$
But this contradicts the fact that $z' = 1/c(w/6)$ (formerly $x'$) is not divisible by any prime divisor of $\mathfrak{p}_2$.  This proves that $l(z,w) \neq 0$ 
and therefore
$$g(z,w) = g(c(w/6), c(w/6)^{\bar \tau}) = 0,$$
which implies that $c(w/6)$ is a periodic point for $g(x,y)$.  Replacing $w$ by $w+3$ shows that the same is true for $s(w/6+1/2)$.
\end{proof}

\noindent {\bf Remark.} Note that $c(w/6)$ is a periodic point for $g(x,y)$, irrespective of whether $3 \mid N(w)$ or not.  This is the case even when 
$\left(\frac{d_K}{3}\right) = -1$, when $w/6$ is not the basis quotient of a proper ideal in $\textsf{R}_{-d}$. \medskip

\section{Determining all periodic points for $k(x,y)$.}

\begin{thm}
We have the identity
\begin{equation}
2c(w/2)^{\bar \tau} \cdot 2c(w/2) = 1+s(w/2).
\label{eqn:5.5}
\end{equation}
\label{thm:6.1}
\end{thm}

\begin{proof}
Let $x = 2c(w/2), y = 2c(w/2)^{\bar \tau}, z = s(w/2), w = s(w/2)^{\bar \tau}$.  Then
\begin{equation*}
(1-w^2)(1-z^2)-(1+z)^3 = (z+1)(w^2z-w^2-z^2-3z) = -(z+1)k(z,w) = 0.
\end{equation*}
Thus
$$(1-s(w/2)^{2\bar \tau})(1-s(w/2)^2) = (1+s(w/2))^3,$$
which implies that
$$8c(w/2)^{3\bar \tau} \cdot 8c(w/2)^3 = (1+s(w/2))^3$$
and therefore
$$2c(w/2)^{\bar \tau} \cdot 2c(w/2) = \omega^i (1+s(w/2)),$$
for some $i \in \{0,1,2\}$.  In particular,
\begin{equation}
\frak{p}_2^2 \cong 1+s(w/2) \ \Rightarrow \ s(w/2) \equiv -1 \ (\textrm{mod} \ \frak{p}_2^2).
\label{eqn:5.6}
\end{equation}
If $\mathfrak{P}$ is a prime divisor of $\mathfrak{p}_2$ in $\Omega_{3f}$, then in the completion $\Omega_{3f,\frak{P}}$ we have
\begin{align}
\notag s^2(w/2) &= 1-8c^3(w/2) = 1-x^3,\\
\label{eqn:5.7} s(w/2) &= \pm \sqrt{1-x^3} = \pm \sum_{n=0}^\infty{\binom{1/2}{n}(-1)^n x^{3n}}\\
\notag & = \pm \left(1-\frac{x}{2}x^2 - \frac{x^3}{8} x^3 + O(x^5)\right)\\
\notag & \equiv \pm 1 \ (\textrm{mod} \ \frak{p}_2^2).
\end{align}
Then \eqref{eqn:5.6} implies that the leading sign must be negative, i.e.
$$s(w/2) = -\sqrt{1-x^3},\ \ x = 2c(w/2),$$
where the square-root is defined in $\Omega_{3f, \frak{P}}$ by Newton's series.  On the other hand, \eqref{eqn:5.3} yields that
$$2c(w/2)^{\bar \tau} = \frac{1-\sqrt{1-x^3}}{x} = \frac{1+s(w/2)}{2c(w/2)}.$$
This proves the theorem.
\end{proof}

This theorem relates to the identity
$$g_1\left(x,\frac{1+z}{x}\right) = \frac{x^3 +z^2 -1}{x}.$$
This illustrates the connection between the curve $g_1(x,y) = 0$ and the curve $\mathcal{E}$ considered in Theorem \ref{thm:2.1}.
\medskip

From \eqref{eqn:5.7} we see that
\begin{align}
\notag s(w/2)^{\bar \tau} &= -\sum_{n=0}^\infty{\binom{1/2}{n}(-1)^n (x^{\bar \tau})^{3n}},\\
\notag &= -\sum_{n=0}^\infty{\binom{1/2}{n}(-1)^n \left(\frac{1+s(w/2)}{2c(w/2)}\right)^{3n}}\\
\notag &= -\sum_{n=0}^\infty{\binom{1/2}{n}(-1)^n \frac{(1+s(w/2))^{3n}}{(1-s^2(w/2))^n}}\\
\label{eqn:6.3a} &= -\sum_{n=0}^\infty{\binom{1/2}{n}(-1)^n \frac{(1+s(w/2))^{2n}}{(1-s(w/2))^n}}.
\end{align}
As above, we have that $1+s(w/2) \cong \frak{p}_2^2$, while $1-s(w/2) \cong \frak{p}_2$, a fact which follows from
$$(1+s(w/2))(1-s(w/2)) = 8c^3(w/2) \cong \frak{p}_2^3.$$

Denoting by $\frak{o}_\frak{P}$ the ring of integers in $\Omega_{3f,\frak{P}}$, we let $x \in \textsf{D}_2 = \{x \in \frak{o}_\frak{P}: |x|_2 \le 2^{-2}\}$ and set
$$F(x) = -\sum_{n=0}^\infty{\binom{1/2}{n}(-1)^n \frac{x^{2n}}{(2-x)^n}}.$$
Then
\begin{align*}
F(4x) &=  -\sum_{n=0}^\infty{\binom{1/2}{n}(-1)^n \frac{4^{2n} x^{2n}}{(2-4x)^n}}\\
&=  -\sum_{n=0}^\infty{\binom{1/2}{n}(-1)^n 2^{3n} \frac{x^{2n}}{(1-2x)^n}}\\
&= -\sum_{n=0}^\infty{\binom{1/2}{n}(-1)^n 2^{3n} \left(\frac{x^2}{1-2x}\right)^n}.
\end{align*}
Expressing the function $\frac{x^2}{1-2x}$ as a convergent power series in $x$, for $|x|_2 \le 1$, we see that $F(4x) \in \mathbb{Z}_2[[x]]$ 
has 2-adic integral coefficients.  In fact, the series for $F(4x)$ shows that
$$F(4x) = -\left(1-\frac{8x^2}{1-2x}\right)^{1/2} = - 1 + 4x^2 + 8x^3 + 24x^4 + 64x^5 + 192x^6 + 576x^7 + \cdots .$$
Furthermore, solving for $y$ in 
$$k(4x-1,y) = (2-4x)y^2 + 16x^2 + 4x-2 = 0$$
gives that
$$y = \pm \left(1-\frac{8x^2}{1-2x}\right)^{1/2} = \pm F(4x).$$
We conclude that
$$k(x-1,F(x)) = 0, \ \ \textrm{identically, for} \ x \in \textsf{D}_2.$$
Moreover, for such $x$, $y = F(x)$ is the unique
solution of $k(x-1,y) = 0$ for which $F(x) \equiv -1$ mod $4$, and \eqref{eqn:6.3a} shows that
$$F(1+s(w/2)) = s(w/2)^{\bar \tau}.$$
Setting $\widetilde{F}(x) = F(x+1)$, this implies that
\begin{equation}
\widetilde{F}^n(s(w/2)) = s(w/2)^{\bar \tau^n}.
\label{eqn:5.8}
\end{equation}

\begin{thm}
If $\frak{P}$ is a prime divisor of $\mathfrak{p}_2$ in $\Omega_{3f}$ over $K = \mathbb{Q}(\sqrt{-d})$, then the value $s(w/2) \in \Omega_{3f,\frak{P}}$, the completion of $\Omega_{3f}$ with respect to $\mathfrak{P}$, is a periodic point in the ordinary sense of the function $\widetilde{F}(x) = F(x+1)$, which maps the disk
$$\widetilde{\textsf{D}}_2 = \{x \in \frak{o}_\frak{P}: |x+1|_2 \le 2^{-2}\}$$
to itself.
\label{thm:6.2}
\end{thm}

Since $\mathbb{Q}(s(w/2)) = \Omega_{3f}$, it is clear that the minimal period of $s(w/2)$ with respect to the algebraic function
$\widetilde{F}(x)$ is equal to the order of the automorphism $\bar \tau$ in $\textrm{Gal}(\Omega_{3f}/\mathbb{Q})$. \medskip

The next lemma is easily proved using induction.

\begin{lem}
Define, inductively,
\begin{align*}
R^{(1)}(x,x_1) &= k(x,x_1),\\
R^{(n)}(x,x_n) &= \textrm{Res}_{x_{n-1}}(R^{(n-1)}(x,x_{n-1}),k(x_{n-1},x_n)), \ \ n \ge 2.
\end{align*}
Then
$$R^{(n)}(x,x_n) \equiv (x_n^{2^n}+x)(x+1)^{2^n-1} \ \textrm{mod} \ 2.$$
Further, if $R_n(x) = R^{(n)}(x,x)$, then
$$R_n(x) \equiv (x^{2^n}+x)(x+1)^{2^n-1} \ \textrm{mod} \ 2.$$
\label{lem:6.3}
\end{lem}

\begin{prop}
The polynomial $R^{(n)}(x,x_n)$ defined in the Lemma satisfies
$$R^{(n)}(x_n,x)=-(x_n-1)(x_n+1)^{2^n-2}x^{2^n}+S_n(x_n,x),$$
where $S_n(x_n,x)$ is monic in $x_n$ with $\textrm{deg}_{x_n} \ S_n(x_n,x) = 2^n$, $\textrm{deg}_x \ S_n(x_n,x) \le 2^n-2$, and $S_n(x_n,0)=x_n(x_n+3)(x_n-3)^{2^n-2}$.
\label{prop:6.4}
\end{prop}

\begin{proof}
We temporarily define $h(x,y) = k(y,x)$ and note that the periodic points for $h(x,y)$ are the same as the periodic points for $k(x,y)$.  Moreover, the polynomial $R^{(n)}(x,x_n)$ for $k$ is the same as the corresponding polynomial for $h$, with $x$ and $x_n$ interchanged.  We use this to prove the proposition by induction.  \medskip

When $n=1$, we have
$$R^{(1)}(x_1,x)=h(x,x_1)=x^2+x_1^2+3x_1-x^2x_1=-(x_1-1)x^2+x_1(x_1+3),$$
with $S_1(x_1,x)=x_1(x_1+3)$ satisfying the assertions.  Assume the result is true for $n-1$, for $n\ge2$.  Then
\begin{align*}
R^{(n)}(x_n,x)&=\textrm{Res}_{x_{(n-1)}}\big(R^{(n-1)}(x_{n-1},x),h(x_{n-1},x_n)\big)\\
&=\textrm{Res}_{x_{n-1}}\big(h(x_{n-1},x_n),R^{(n-1)}(x_{n-1},x)\big)\\
&=(x_n-1)^{2^{n-1}}\prod_{i=1}^2{R^{(n-1)}(\rho_i,x)};
\end{align*}
where $\rho_1,\rho_2$ are the roots of $h(x_{n-1},x_n)=-(x_n-1)x_{n-1}^2+x_n(x_n+3)$ as a polynomial in $x_{n-1}$, given explicity as $\rho_1=\sqrt{\frac{x_n(x_n+3)}{(x_n-1)}}$ and $\rho_2=-\sqrt{\frac{x_n(x_n+3)}{(x_n-1)}}$. Thus, $\sigma_1 = \rho_1+\rho_2=0$ and $\sigma_2 = \rho_1\rho_2=-\frac{x_n(x_n+3)}{(x_n-1)}$. Hence,
$$R^{(n)}(x_n,x)=(x_n-1)^{2^{n-1}}\prod_{i=1}^2{\Big[-(\rho_i-1)(\rho_i+1)^{2^{n-1}-2}x^{2^{n-1}}+S_{n-1}(\rho_i,x)\Big]}.$$
The coefficient of the leading term in $x$ (namely, $x^{2^n}$) in the above product is
\begin{align*}
(x_n-1)^{2^{n-1}}&(\rho_1-1)(\rho_2-1)(\rho_1+1)^{2^{n-1}-2}(\rho_2+1)^{2^{n-1}-2}\\
&=(x_n-1)^{2^{n-1}}(\rho_1\rho_2+1)(\rho_1\rho_2+1)^{2^{n-1}-2}\\
&=(x_n-1)^{2^{n-1}}\bigg[-\frac{x_n(x_n+3)}{(x_n-1)}+1\bigg]^{2^{n-1}-1}\\
&=(x_n-1)\big[-x_n(x_n+3)+(x_n-1)\big]^{2^{n-1}-1}\\
&=-(x_n-1)\big[(x_n+1)^2\big]^{2^{n-1}-1}\\
&=-(x_n-1)(x_n+1)^{2^n-2}.
\end{align*}
Furthermore,
\begin{align*}
S_n(x_n,x) &= (x_n-1)^{2^{n-1}} \big[-(\rho_1-1)(\rho_1+1)^{2^{n-1}-2}S_{n-1}(\rho_2,x)\\
& \ \ \ -(\rho_2-1)(\rho_2+1)^{2^{n-1}-2}S_{n-1}(\rho_1,x)\big]x^{2^{n-1}}\\
& \ \ \ + (x_n-1)^{2^{n-1}}S_{n-1}(\rho_1,x)S_{n-1}(\rho_2,x).
\end{align*}
By the induction assumption, $\textrm{deg}_x S_n(x_n,x) \le 2^{n-1}-2+2^{n-1} = 2^n-2$.  Also, the leading term in $x_n$ of the bracketed term 
$[*]$ in the last expression is
\begin{align*}
(x_n-1)^{2^{n-1}}&\mathcal{L}\big[-(\rho_1-1)(\rho_1+1)^{2^{n-1}-2}\rho_2^{2^{n-1}} -(\rho_2-1)(\rho_2+1)^{2^{n-1}-2}\rho_1^{2^{n-1}}\big]\\
&= (x_n-1)^{2^{n-1}}\left(-\rho_1^{2^{n-1}-1}\rho_2^{2^{n-1}} -\rho_2^{2^{n-1}-1}\rho_1^{2^{n-1}}\right)\\
&=(x_n-1)^{2^{n-1}} (\rho_1 \rho_2)^{2^{n-1}-1}(-\rho_2-\rho_1)\\
&= (x_n-1)(x_n^2+3x_n)^{2^{n-1}-1}(\rho_2+\rho_1) = 0,
\end{align*}
while the other terms in $[*]$ contribute expressions
\begin{align*}
&(x_n-1)^{2^{n-1}}( \rho_1^j \rho_2^k+\rho_1^k \rho_2^j) x^r, \ \ \ \ 2^{n-1}-1 \ge j \ge k\\
&= (x_n-1)^{2^{n-1}} (\rho_1 \rho_2)^k \big[\rho_1^{j-k}+\rho_2^{j-k}\big] x^r\\
&= (x_n-1)^{2^{n-1}-k}(-x_n)^k(x_n+3)^k \sigma_1^{j-k} x^r + \ \textrm{lower degree terms};
\end{align*}
where the degree of the first term in $x_n$ is $\le 2^{n-1}+k \le 2^n-1$ and is only nonzero when $j=k$. \medskip

On the other hand, the sum of the terms without $x$ in $S_n(x_n,x)$ is given by
\begin{align*}
S_n(x_n,0)&=(x_n-1)^{2^{n-1}}\big[(S_{n-1}(\rho_1,0)\times S_{n-1}(\rho_2,0)\big]\\
&=(x_n-1)^{2^{n-1}}\big[\rho_1\rho_2(\rho_1+3)(\rho_2+3)(\rho_1-3)^{2^{n-1}-2}(\rho_2-3)^{2^{n-1}-2}\big]\\
&=(x_n-1)^{2^{n-1}}\big[\rho_1\rho_2(\rho_1\rho_2+9)(\rho_1\rho_2+9)^{2^{n-1}-2}\big]
\end{align*}
\begin{align*}
&=(x_n-1)^{2^{n-1}}\bigg[-\frac{x_n(x_n+3)}{(x_n-1)}\bigg]\bigg[-\frac{x_n(x_n+3)}{(x_n-1)}+9\bigg]^{2^{n-1}-1}\\
&=x_n(x_n+3)\big[x_n(x_n+3)-9(x_n-1)\big]^{2^{n-1}-1}\\
&=x_n(x_n+3)(x_n-3)^{2^n-2}.
\end{align*}
The same calculation as above shows that the other terms in the product $(x_n-1)^{2^{n-1}}S_{n-1}(\rho_1,x)S_{n-1}(\rho_2,x)$ have degree
at most $2^n-1$ in $x_n$.  This yields that $\textrm{deg} \ S_n(x_n,0) = 2^n$.
\end{proof}

\begin{cor}
For $n \ge 1$, $\textrm{deg} \ R_n(x) = 2^{n+1}-1$, and all the periodic points for $k(x,y) = 0$ are algebraic integers.
\label{cor:6.5}
\end{cor}

\begin{thm}
All the periodic points of the multivalued function $y = \textsf{F}(x)$ defined by $k(x,\textsf{F}(x)) = 0$ in the algebraic closure $\overline{\mathbb{Q}}_2$ 
lie in the maximal unramified extension $\textsf{K}_2$ of the $2$-adic field $\mathbb{Q}_2$.
\label{thm:6.6}
\end{thm}

\begin{proof}
Appealing to an argument similar to the one given in \cite[p. 335]{mdeu}, the last congruence in Lemma \ref{lem:6.3} and Hensel's Lemma imply 
that $R_n(x)$ has at least $2^n-1$ distinct roots in $\textsf{K}_2$, one of which is $x = 0$, since
$$R_1(x) = -x(x+1)(x-3)$$
is a divisor of $R_n(x)$.  The other $2^n-2$ roots reduce (mod $2$) to elements of $\mathbb{F}_{2^n}$ which have degree at least $2$ over the prime field,
and are therefore units in $\textsf{K}_2$.  Moreover, the congruence
$$k(x,y) \equiv (x+1)(y^2+x) \ (\textrm{mod} \ 2)$$
implies that if $x \rightarrow \bar x$ denotes the reduction map (mod $2$) from $\textsf{K}_2$ to $\overline{\mathbb{F}}_2$, and
\begin{equation}
k(\alpha_0,\alpha_1) = \cdots = k(\alpha_i,\alpha_{i+1}) = \cdots = k(\alpha_{n-1}, \alpha_0) = 0,
\label{eqn:5.9}
\end{equation}
where $\bar \alpha_0 \notin \mathbb{F}_2$, then $\bar \alpha_i \notin \mathbb{F}_2$, for all $i$.  It follows that the elements of the orbit of $\alpha_0$ 
are among the $2^n-2$ roots we have found so far.  Now we use the identity
\begin{equation}
(x+1)^2(y+1)^2 k\left(\frac{-x+3}{x+1},\frac{-y+3}{y+1}\right) = 16k(y,x).
\label{eqn:5.10}
\end{equation}
This identity implies that for every chain of equations \eqref{eqn:5.9} of length $n$, there is a corresponding chain, with elements
$$\beta_i = \frac{-\alpha_i+3}{\alpha_i+1}, \ \ 0 \le i \le n-1,$$
for which
$$k(\beta_0,\beta_{n-1}) = k(\beta_{n-1},\beta_{n-2}) = \cdots = k(\beta_1,\beta_0) = 0.$$
Hence, every periodic point $\alpha_0$ corresponds to another periodic point $\beta_0$.  Also, 
$$\beta_i + 1 = \frac{4}{\alpha_i+1}.$$
Thus, $\alpha_i + 1$, having a reduction whose degree is at least $2$, must be a unit, and we see that $\beta_i \equiv -1$ (mod $4$) in $\textsf{K}_2$.
Since $\alpha_i \rightarrow \beta_i$ is a 1-1 map, $\beta_i$ is a root of $R_n(x)$ whenever $\alpha_i$ is, and the $\beta$'s are distinct from all of the $\alpha$'s.  
Thus, $R_n(x)$ has at least $2^{n+1}-4$ roots whose reductions lie in $\mathbb{F}_{2^n}-\mathbb{F}_2$, and therefore at least $2^{n+1}-1$ distinct roots altogether.  
But since $\textrm{deg}(R_n(x)) = 2^{n+1}-1$, this gives all the roots of $R_n(x)$, and all of these roots lie in $\textsf{K}_2$.
\end{proof}

\begin{cor}
If a periodic point $\alpha \in \textsf{K}_2$ for $k(x,y) = 0$ has minimal period $n \ge 2$, and its reduction $\bar \alpha = \alpha \ (\textrm{mod} \ 2)$ does not lie in $\mathbb{F}_2$, then $\mathbb{F}_2(\bar \alpha) = \mathbb{F}_{2^n}$. Furthermore, if $\sigma \in \textrm{Gal}(\mathbb{Q}_2(\alpha)/\mathbb{Q}_2)$ is the element corresponding to $\bar \sigma: \bar x \rightarrow \bar x^2$ in the isomorphism
$$\textrm{Gal}(\mathbb{Q}_2(\alpha)/\mathbb{Q}_2) \cong \textrm{Gal}(\mathbb{F}_2(\bar \alpha)/\mathbb{F}_2),$$
then $k(\alpha^\sigma, \alpha) = 0$, $k(\beta, \beta^\sigma) = 0$, and the elements of the orbit of any periodic point $\alpha$ are conjugates of $\alpha$ over $\mathbb{Q}_2$.
\label{cor:6.7}
\end{cor}

\begin{proof}
The congruence for $k(x,y)$ in the proof shows that if $k(\alpha, \alpha_1) = 0$ in $\textsf{K}_2$, then $\bar \alpha_1^2 \equiv \bar \alpha$ (mod $2$).  This shows 
that $\bar \alpha$ and $\bar \alpha_1$ are conjugate over $\mathbb{F}_2$.  Then, the elements $\alpha_i$ in the proof satisfy
$$\bar \alpha_i^{2^i} \equiv \bar \alpha \ (\textrm{mod} \ 2), \ \ i = 0, \dots, n.$$
Moreover,  $i = n$ is the first index for which $\bar \alpha_i = \bar \alpha$, since the reduction of $R_n(x)$ (mod $2$) has no multiple roots other than $1$; and the periodic points of minimal period $n$ are those which are roots of
$$\textsf{P}_n(x) = \prod_{d \mid n}{R_d(x)^{\mu(n/d)}},$$
whose reductions are therefore roots of the product
$$ \prod_{d \mid n}{(x^{2^d} -x)^{\mu(n/d)}}.$$
Hence, the minimal polynomial of $\bar \alpha$ over $\mathbb{F}_2$ has degree $n$, proving the first assertion.  This also shows that the elements of the orbit of $\alpha$ must be conjugates of $\alpha$ over $\mathbb{Q}_2$, since they reduce (mod $2$) to the roots of an irreducible polynomial of degree $n$.  Similarly, since the $\alpha_i$ in the above proof are conjugates of $\alpha$, the elements $\beta_i$ of the orbit of $\beta = \frac{-\alpha+3}{\alpha+1}$ must be conjugates of $\beta$.  The relations $k(\alpha,\alpha_1) = 0$ and $\sigma: \alpha_1 \rightarrow \alpha$ imply that $k(\alpha, \alpha^{\sigma^{-1}}) = 0$ and therefore $k(\alpha^\sigma,\alpha) = 0$.  The identity \eqref{eqn:5.10} now shows that $k(\beta, \beta^\sigma) = 0$.  For $n=1$ the last statement is obvious, since the fixed points lie in $\mathbb{Q}$.
\end{proof}

\begin{thm}
The minimal period of $s(w/2)$ as a periodic point for $k(x,y) = 0$ is the order of the automorphism $\displaystyle \bar \tau = \left(\frac{\Omega_{3f}/K}{\frak{p}_2}\right)$ in $\textrm{Gal}(\Omega_{3f}/K)$, and is always even.  The same holds for the minimal period of $\pm 3 r(w/2)$.
\label{thm:6.8}
\end{thm}

\begin{proof}
By the proof of Theorem \ref{thm:6.6}, an orbit of minimal size $n \ge 2$ consists of elements of $\textsf{K}_2$, which are either all congruent 
to $-1$ mod $4$ or all units whose reductions have degree at least $2$ in $\overline{\mathbb{F}}_2$.  \medskip

On the other hand, $s(w/2)$ satisfies the congruence $s(w/2) \equiv -1$ (mod $\frak{p}_2^2$) in $\Omega_{3f}$ and therefore 
$$s(w/2) \equiv -1 \ \textrm{mod} \ 4 \ \ \textrm{in} \ \Omega_{3f,\frak{P}} \subset \textsf{K}_2.$$
Hence, the orbit of $\beta_0 = s(w/2)$ is $\{\beta_0, \beta_1, \dots, \beta_{n-1}\}$, where each $\beta_i \equiv -1$ (mod $4$) in $\textsf{K}_2$.  But, the observation above concerning the value of $F(x)$ shows that $\beta_1 = F(1+\beta_0) = \widetilde{F}(\beta_0)$, and by induction, $\beta_i = \widetilde{F}^i(\beta_0)$.  Now,
\eqref{eqn:5.8} shows that $\beta_n = \beta_0 = \widetilde{F}^n(\beta_0) = \beta_0^{\bar \tau^n}$, i.e., $s(w/2) = s(w/2)^{\bar \tau^n}$.  Hence $n$, being minimal, must be the order of $\bar \tau$.  \medskip

Now, Theorem \ref{thm:4.3} states that $\mp \frac{1}{r(w/2)}$ is a conjugate of $s(w/2)$ over $\mathbb{Q}$, so it has the same period as $s(w/2)$.  Then \eqref{eqn:4.7} shows that
$\pm 3r(w/2)$ has the same period as $\mp \frac{1}{r(w/2)}$.  \medskip

By the Artin map, this common period is equal to the order of $\frak{p}_2$ in the ideal class group 
of $\textsf{R}_{-9d}$, and this, in turn, is equal to the smallest integer $n \ge 1$, for which the equation
$$4 \cdot 2^n = x^2+9dy^2, \ \ 3 \nmid x,$$
has a primitive solution $(x,y) \in \mathbb{Z}^2$.  Considering this equation modulo $3$ shows that the left side is $2^{n+2} \equiv 1$ mod $3$, from
which it is clear that $n$ must be even. 
\end{proof}

We define
\begin{equation}
\mathfrak{g}(x,y) = x^2 y^2 + y^2 + 2x,
\label{eqn:6.1}
\end{equation}
as in \cite{mdeu}.  With this we note the resultant factorization
\begin{align}
\notag \textrm{Res}_w(\textrm{Res}_z&((k(z,w),\mathcal{Q}(x,z)), \mathcal{Q}(y,w)) = 2^{48} 3^{18} (x^2y^2 + y^2 + 2x)^4 \\
\notag & \ \ \ \ \ \times  (x^2y^2 + y^2 - 2x)^4(2xy^2 + x^2 + 1)^4(-2xy^2 + x^2 + 1)^4\\
\label{eqn:6.2} & =  2^{48} 3^{18} \mathfrak{g}(x,y)^4 \cdot \mathfrak{g}(-x,y)^4 \cdot y^8 \mathfrak{g}\left(x,\frac{1}{y}\right)^4 \cdot y^8 \mathfrak{g}\left(-x,\frac{1}{y}\right)^4.
\end{align}

We also note that the elliptic curve $\mathcal{E}_g: \mathfrak{g}(x,y) = 0$ is isomorphic to the curve
$$\textsf{E}_{-1}: \ X^2Y^2+Y^2-X^2+1 = 0$$
introduced just before Theorem \ref{thm:3.3}, by the birational map $\Phi: \mathcal{E}_g \rightarrow \textsf{E}_{-1}$ given by
$$\Phi: (x,y) \rightarrow \left(\frac{x-1}{x+1},y \right).$$
These are curves with $j$-invariant $j = 1728$. 
\medskip

\begin{thm}
The only periodic points for $k(x,y)$ are $0, -1, 3$, and the values and conjugates over $\mathbb{Q}$ of $s(w/2), s(w/6), s(w/6 + 1/2)$, and $\pm 3 r(w/2)$ (for some discriminant $-d \equiv 1$ mod $8$), the latter
sign depending on whether $-1/r(w/2)$ or $+1/r(w/2)$ is conjugate to $s(w/2)$ over $\mathbb{Q}$.  All periodic points for $k(x,y)$ are roots of $\mathcal{Q}_d(X)$, for suitable discriminants $-d \equiv 1$ (mod $8$).
\label{thm:6.9}
 \end{thm}
 
\begin{proof}
To prove this, let $\alpha$ be any periodic point for $k(x,y) = 0$, whose reduction mod $2$ does not lie in the prime field.  Then by Corollary \ref{cor:6.7}, $k(\alpha^\sigma, \alpha) = 0$, where $\sigma$ is a lift of the Frobenius automorphism to $\mathbb{Q}_2(\alpha)$.  Now put $z = \alpha^\sigma, w = \alpha$ and let 
 $y$ be a solution of $\mathcal{Q}(y,w) = \mathcal{Q}(y,\alpha) = 0$, so that $x = y^\sigma$ satisfies $\mathcal{Q}(x,z) = \mathcal{Q}(y^\sigma,\alpha^\sigma) = 0$.  (If necessary, we may extend the lift $\sigma$ on $\mathbb{Q}_2(\alpha)$ to the field $L = \mathbb{Q}_2(\alpha,y)$.)  Then the right side of \eqref{eqn:6.2} must be $0$.  Writing the equation $\frac{1}{4} \mathcal{Q}(y,w) = 0$ in the form
 $$4w^3\left(y+\frac{1}{y}\right)^2 = (w-1)(w+3)^3$$
 shows that $\mathcal{Q}(y,w) = 0$ has, for a given $w$, the solutions $\pm y, \frac{\pm 1}{y}$.  Thus we may choose the solution $y$ so that it is integral. \medskip
 
 If $y \in \mathfrak{o}_L$, the ring of integers in $L$, we claim that the last two factors in \eqref{eqn:6.2} cannot be $0$.  Let $\pi$ be a prime element in $L$.  If $2xy^2+x^2+1 = 0$ then
 $$0 = 2xy^2+x^2+1 \equiv (x+1)^2 \ (\textrm{mod} \ \pi),$$
which gives that $x \equiv y \equiv 1$ (mod $\pi$).  Now
$$\frac{1}{4}\mathcal{Q}(y,w) = 4y^4w^3 + (-w^4 - 18w^2 + 27)y^2 + 4w^3 = 0,$$
and $y \equiv 1$ (mod $\pi$) gives that $(w+1)^4 = (\alpha+1)^4 \equiv 0$ (mod $\pi$), which is impossible, since it would imply 
$\bar \alpha \in \mathbb{F}_2 \subset \frak{o}_L/(\pi)$.  \medskip

Hence, $\mathfrak{g}(x,y) = 0$ or $\mathfrak{g}(-x,y) = 0$.  But if $\mathfrak{g}(-x,y) = 0$, then also $\mathfrak{g}(-x,-y) = 0$ by \eqref{eqn:6.1}, 
and $-x = (-y)^\sigma$.  Thus, we may assume $x$ and $y$ 
satisfy $\mathfrak{g}(x,y) = 0$, i.e. $\mathfrak{g}(x,x^{\sigma^{-1}}) = 0$.  But then $x$ must be a periodic point for the algebraic function defined by 
$\mathfrak{g}(X,Y) = 0$.  By the results of \cite[Thm. 2, pp.326-327] {mdeu}, $x$ must be $0$ or $-1$ or a conjugate of one of the elements $\pi$ or $\xi$ satisfying 
$\pi^4+\xi^4 = 1$ in a ring class field $\Omega_f$ over the quadratic field $K$, 
for which $-d = d_K f^2 \equiv 1$ (mod $8$).  Now, $x,y \in \{0, -1\}$ implies that $w \in \{0,-1,3\}$.  Otherwise, $x, y$ are roots of the polynomial $b_d(x)$ and $\alpha$ is a root of the polynomial
$$\mathcal{Q}_d(X) = \textrm{Res}_y(\mathcal{Q}(y,X),b_d(y))$$
from \eqref{eqn:4.9}.  Note that this implies that $x, y \in \textsf{K}_2$ also. \medskip

Now assume $\beta$ is a periodic point for $k(x,y)$ satisfying $\beta \equiv -1$ (mod $4$) in $\textsf{K}_2$.  
By the proof of Theorem \ref{thm:6.6}, $\beta$ is 
given by the formula
$$\beta = \frac{-\alpha+3}{\alpha+1}, \ \textrm{with} \ \bar \alpha \notin \mathbb{F}_2.$$
By the above proof, $\alpha$ is a root of $\mathcal{Q}_d(X)$; and by \eqref{eqn:4.10} and \eqref{eqn:4.11} we find that
$$\mathcal{Q}_d(\beta) = \mathcal{Q}_d\left(\frac{-\alpha+3}{\alpha+1}\right) = 0.$$
This proves that all periodic points for $k(x,y)$ are roots of $\mathcal{Q}_d(X)$, for suitable discriminants $-d \equiv 1$ (mod $8$).  \medskip

Now the assertions follow from Theorems \ref{thm:4.3}, \ref{thm:4.8}, and the statement and last paragraph of the proof of Proposition \ref{prop:4.9}.
\end{proof}

See Table 2, which gives the minimal polynomials dividing $R_n(x)$ for periods $n = 2 - 6$ along with a root of each irreducible factor 
in terms of $s(\tau)$ and $r(\tau)$.  \medskip

Based on Theorem \ref{thm:6.9} we prove the following.

\begin{prop}
For $w = \frac{a+\sqrt{-d}}{2}, -d \equiv 1$ (mod $8$) and $8 \mid a^2+d$, the value $\alpha = -s(w/2)$ is not a periodic point for $k(x,y) = 0$.  The same is true for
$s(w/4)$, when $2^5 \mid a^2+d$.
\label{prop:6.10}
\end{prop}

For this we need a lemma.

\begin{lem}
Suppose that $-d = d_K f^2$ and $-d' = d_{K'} f'^2 $ are both $\equiv 1$ (mod $8$), for distinct discriminants $d_K \neq d_{K'}$ of $K = \mathbb{Q}(\sqrt{-d})$ and 
$K' = \mathbb{Q}(\sqrt{-d'})$.  Then the ring class fields $\Omega_f$ and $\Omega_{f'}$ are distinct number fields.
\label{lem:6.11}
\end{lem}

\begin{proof}
Assume, with the given conditions, that $L = \Omega_f = \Omega_{f'}$. Let $G = \textrm{Gal}(L/\mathbb{Q})$ and let $H, H'$ denote the subgroups of $G$ corresponding to $K$ and $K'$ in the Galois correspondence.  Then $H$ and $H'$ are both abelian over $\mathbb{Q}$ and $[G:H] = [G:H'] = 2$.  
Let $\psi \in H' \setminus H$.  Then the fact that
$G$ is generalized dihedral (over $\mathbb{Q}\subset K$) implies that $\psi^2 = 1$ in $G$ and $\psi \sigma \psi^{-1} = \sigma^{-1}$, $\forall \sigma \in H$.  
Furthermore, for $\sigma \in H \cap H'$ we have $\psi \sigma \psi^{-1} = \sigma$, since $H'$ is abelian.  Hence, $\sigma^2 = 1$, for $\sigma \in H \cap H'$.  
It follows that $\psi^2 = 1, \forall \psi
\in H'$.  By symmetry the same is true for $\sigma \in H$.  But then $\sigma^{-1} = \sigma$ implies that $\psi$ and $\sigma$ commute, for all $\psi \in H', \sigma \in H$.  Thus, $G = HH'$ is abelian.  Now $L = \Omega_f = \Omega_{2f}$, since $f$ is odd 
and $d_K \equiv 1$ (mod $8$).  Then $\mathbb{Q}(j(\sqrt{-d})) \subset \Omega_f$ being abelian over $\mathbb{Q}$ and Weber's Theorem in Kani \cite[Thm. 32]{ka} imply that $d \equiv 7$ mod $8$ is an idoneal (or ``convenient" \cite[pp. 53-56]{cox}) number, and
Grube's result \cite[Cor. 8]{ka} implies that $d = 7$ or $15$.  Since the same primes ramify in $\Omega_f$ and $\Omega_{f'}$, we must also have $d' = 7$ or $15$, respectively, and $K = K'$, a contradiction.
\end{proof}

\begin{proof}(of Proposition \ref{prop:6.10})
Suppose that $\alpha = -s(w/2)$ is a periodic point for $k(x,y) = 0$.  Then $\alpha$ must be a root of $\mathcal{Q}_{d'}(y)$, for some discriminant $-d' \equiv 1$ (mod $8$).  By Theorem
\ref{thm:3.10}, $\mathbb{Q}(-s(w/2)) = \Omega_{3f}$.  If $d' \neq d$, then by the proof of Proposition \ref{prop:4.9}, we must also have $\mathbb{Q}(-s(w/2)) = \Omega_{f'}$ or $\Omega_{3f'}$
over $K = \mathbb{Q}(\sqrt{-d'})$.  But the lemma shows that $\Omega_{3f} \neq \Omega_{f'}, \Omega_{3f'}$ unless $K = K'$.  Thus, $K = K'$.  Since the class field theoretic conductor 
over $K$ of the field $\Omega_{3f}$ is just $\mathfrak{f} = (3f)$ \cite[Problem 9.20, p. 177]{cox}, we must have $3f = f'$ or $3f'$.  The latter possibility implies that $f = f'$ and $d = d'$.  
But $-s(w/2)$ has the same degree over $\mathbb{Q}$ as $s(w/2)$, and the latter is a root of the unique factor of $\mathcal{Q}_d(y)$ of largest degree.  Hence $-s(w/2)$ and $s(w/2)$ 
would have to be conjugate.  But this is false, as in the remark after the proof of Theorem \ref{thm:5.4}.  We conclude that $3f = f'$ and $9d = d'$.  In that case, $-\alpha = s(w/2)$ being a root of $\mathcal{Q}_d(y)$ and $\alpha = -s(w/2)$ being a root of $\mathcal{Q}_{9d}(y)$ imply that
\begin{align*}
\frac{(\alpha-1)(\alpha+3)^3}{(\alpha+1)(\alpha-3)^3} &= \frac{(s(w/2)+1)(s(w/2)-3)^3}{(s(w/2)-1)(s(w/2)+3)^3} = \left(\frac{1-\xi^2}{1+\xi^2}\right)^2,\\
\frac{(\alpha-1)(\alpha+3)^3}{(\alpha+1)(\alpha-3)^3} &= \left(\frac{1+\xi'^2}{1-\xi'^2}\right)^2,
\end{align*}
for some roots $\xi, \xi'$ of $b_d(x)$ and $b_{d'}(x) = b_{9d}(x)$.  But the relation
$$\frac{(1-x^2)^2}{(1+x^2)^2}-\frac{(1+y^2)^2}{(1-y^2)^2} = -4\frac{(x^2 + y^2)(x^2y^2 + 1)}{(x^2+1)^2(y^2-1)^2}$$
implies then that
$$\xi^2+\xi'^2 = 0 \ \ \textrm{or} \ \ \xi^2 \xi'^2+1 = 0.$$
These equations are impossible, since they both imply that $-1$ is a square in $\Omega_{3f}$, an extension of $\mathbb{Q}$ which is unramified at the prime $p = 2$.
Hence, $-s(w/2)$ cannot be a root of any polynomial $\mathcal{Q}_{d'}(y)$, and is therefore not a periodic point for $k(x,y) = 0$.  The last statement is a 
consequence of Theorem \ref{thm:3.12}.
\end{proof}

From Theorem \ref{thm:6.9} we obtain the following.

\begin{thm}
The only periodic points for the algebraic function defined by $g_1(x,y) = 0$ are $0,\omega^i,-2\omega^i$ for $i \in \{0,1,2\}$, and the values and conjugates
over $\mathbb{Q}$ of
\begin{equation*}
2\omega^i c(w/2), \ 2\omega^i c(w/6), \ 2\omega^i c(w/6+1/2), \ 2\omega^ic(-1/3w),
\end{equation*}
where
$$w = \frac{k+\sqrt{-d}}{2}, \ \ k^2+d \equiv 0  \ \textrm{mod} \ 2^3,$$
and $-d = d_Kf^2$ runs over all negative quadratic discriminants satisfying $-d \equiv 1$ (mod $8$).
The last periodic point can also be written as
$$2\omega^ic(-1/3w) = \omega^i \frac{1-2c(3w/2)}{1+c(3w/2)} = \omega^i \mathcal{A}(2c(3w/2)) = \omega^i B\left(\frac{-1}{c(3w/2)}\right),$$
where
$$\mathcal{A}(x) = \frac{2-2x}{2+x}, \ \ B(x) = \frac{x+2}{x-1},$$
and $\frac{-1}{c(3w/2)}$ and $2c(3w/2)$ are periodic points for $g_1(x,y)=0$.
\label{thm:6.12}
\end{thm}

\noindent {\bf Remark.} It is clear that $\mathcal{A}(x)$ takes periodic points for $g_1(x,y)$ to other periodic points, by the identity
\begin{equation}
(2+x)^2 (2+y)^2 g_1\left(\mathcal{A}(x),\mathcal{A}(y)\right) = 6^2 g_1(y,x).
\label{eqn:6.10}
\end{equation}
Similarly,
\begin{equation}
(x-1)^2 (y-1)^2 g_1(B(x),B(y)) = 9g_1(x,y).
\label{eqn:6.11}
\end{equation}

\begin{proof}
This will follow from the resultant calculation
\begin{align}
\notag \textrm{Res}_w(\textrm{Res}_z(&g_1(z,w),1-x^2-z^3),1-y^2-w^3)\\
\notag &= (-xy^2 + x^2 + y^2 + 3x)^3 (xy^2 + x^2 + y^2 - 3x)^3\\
\label{eqn:6.12} &= k(x,y)^3k(-x,y)^3.
\end{align}
Let $\gamma= \gamma_0 = \gamma_n$ be a periodic point for $g_1(x,y) = 0$ and $\gamma_i, 1 \le i \le n-1$, be numbers for which
$$g_1(\gamma_i,\gamma_{i+1}) = 0, \ \ 0 \le i \le n-1.$$
Find $\alpha_i, 0 \le i \le n$, so that $\alpha_i^2 = 1- \gamma_i^3$.  Then, for each $i$, \eqref{eqn:6.12} gives that
$$k(\alpha_i, \alpha_{i+1}) = 0 \ \ \textrm{or} \ \ k(-\alpha_i, \alpha_{i+1}) = 0.$$
By choosing the sign of $\alpha_0$, we can assume the first equation holds for $i = 0$.  If the second equation holds (for some $i \ge 1$) , 
we use that $k(x,y) = k(x,-y)$ and replace $\alpha_i$ by $-\alpha_i$, so that
$$k(\alpha_{i-1}, -\alpha_i) = k(-\alpha_i,\alpha_{i+1}) = 0.$$
Continuing in this way, the signs of the elements $\alpha_i$ can be chosen so that
$$k(\alpha_i, \alpha_{i+1}) = 0, \ \ 0 \le i \le n-1, \ \alpha_n = \alpha_0.$$
Then $\alpha = \alpha_0 = \alpha_n$ is a periodic point for $k(x,y) = 0$; and Theorem \ref{thm:6.9} implies that some conjugate $\alpha'$ of $\alpha$ satisfies
$$\alpha' \in \{0, -1, 3, s(w/2), s(w/6), s(w/6+1/2), \pm 3r(w/2)\},$$
for some discriminant $-d \equiv 1$ (mod $8$).
Then the relation $x^2 = 1-z^3$ in \eqref{eqn:6.12} implies that $\gamma'^3 = 1 - \alpha'^2$, so some conjugate $\gamma'$ of $\gamma$ satisfies
$$\gamma' \in \{0,\omega^i,-2\omega^i, 2\omega^i c(w/2), \ 2\omega^i c(w/6), \ 2\omega^i c(w/6+1/2), \ 2\omega^ic(-1/3w)\};$$
where $i \in \{0,1,2\}$ and the last value follows from Lemma \ref{lem:2.2} and 
$$\gamma'^3 = 1-9r^2(w/2) = 1-s^2(-1/3w).$$
The value $2c(-1/3w)$ is certainly a periodic point for $g_1(x,y)$, by \eqref{eqn:6.12}, since $\pm 3r(w/2)$, for some sign, is a periodic point for
$k(x,y)$.  The other values for $\gamma'$ are also periodic points for $g_1(x,y) = 4g(x/2,y/2)$ by Theorem \ref{thm:5.2} and Proposition \ref{prop:5.8}.  
Furthermore, \cite[pp. 1714-1715]{akm} gives that
$$2c(-1/3w) = \frac{1-2c(3w/2)}{1+c(3w/2)} = \mathcal{A}(2c(3w/2)).$$
Applying $B(x)$ to this formula and noting $B^2(x) = x$ also yields that
$$B(2c(-1/3w)) = \frac{-1}{c(3w/2)} \ \Rightarrow \ 2c(-1/3w) = B(-1/c(3w/2)).$$
Since $\mathcal{A}$ and $B$ take periodic points for $k(x,y)$ to other periodic points, it is clear that $-1/c(3w/2)$ and $2c(3w/2)$ are periodic points.
 This completes the proof.
\end{proof}

\noindent {\bf Example.} From Table \ref{tab:3} we know that
\begin{align*}
p_{15}(x) &= x^{12} - 3x^{11} - 15x^{10} + 39x^9 + 51x^8 + 9x^7 - 180x^6 - 837x^5 + 2079x^4\\
 & \ - 1278x^3 + 135x^2 + 54x + 9
 \end{align*}
is the minimal polynomial of $\alpha = -3r(w/2)$ for $\mathbb{Q}(\sqrt{-15})$, which generates $\Omega_3$ over $\mathbb{Q}$.  Here,
$w = \frac{23+\sqrt{-15}}{2}$.  Note that the period of $\alpha$ is $n = 6$.  We calculate that
\begin{align*}
\textrm{Res}_x&(p_{15}(x),1-x^2-z^3) = z^{36} + 27z^{33} + 198z^{30} - 328z^{27} - 4473z^{24}\\
&  + 37053z^{21} + 193947z^{18} - 296424z^{15} - 286272z^{12} + 167936z^9 + 811008z^6\\
& - 884736z^3 + 262144 = q_{15}(z).
\end{align*}
The root $\gamma = \mathcal{A}(2c(3w/2))$ of $q_{15}(z)$ generates the same field over $\mathbb{Q}$ as does $2c(3w/2)$, which is $\Omega_{3 \cdot 3f} = \Omega_{9f}$, by Theorem \ref{thm:3.10}.  Note that $\gamma$ is not conjugate to $2c(3w/2)$, since the latter is a root of
\begin{align*}
&2^{-18}(x+2)^{36} q_{15}\left(\mathcal{A}(x)\right) = p_{135}(x) = x^{36} + 18x^{35} + 192573x^{34} + 3460038x^{33} \\
&  + 33153696x^{32} - 611690337x^{31} + 4655053353x^{30} - 24787366965x^{29}\\
& + 104419059804x^{28} - 326112084643x^{27}+ 928846349469x^{26}
\end{align*}
\begin{align*}
& - 1954061339535x^{25} + 4354422142788x^{24} - 6369632847495x^{23} \\ 
 & + 13079097127845x^{22} - 12124557064791x^{21} + 30983586793713x^{20}\\
 & - 10410555305610x^{19} + 64162661792703x^{18} + 20821110611220x^{17}\\
 &   + 123934347174852x^{16} + 96996456518328x^{15}+ 209265554045520x^{14} \\
 & + 203828251119840x^{13} + 278683017138432x^{12} + 250119851460480x^{11} \\
 & + 237784665464064x^{10} + 166969387337216x^9 + 106925117239296x^8 \\
 & + 50764527544320x^7 + 19067098533888x^6 + 5010967240704x^5  \\
 & + 543190155264x^4 - 113378525184x^3 + 12620464128x^2 - 2359296x + 262144.
\end{align*}
Since
\begin{align*}
q_{15}(z) & \equiv z^{18} (z^{18} + z^{15} + z^6 + z^3 + 1) \ (\textrm{mod} \ 2)\\
p_{135}(x) & \equiv x^{18}(x^{18} + x^{16} + x^{13} + x^{12} + x^{11} + x^9 + x^8 + x^7 + x^5 + x^4\\
& \ + x^3 + x^2 + 1) \ (\textrm{mod} \ 2),
\end{align*} \medskip
 Corollary \ref{cor:6.7} implies that $\gamma = \mathcal{A}(2c(3w/2))$ and $2c(3w/2)$ both have period $3n = 18$ with respect to $g_1(x,y) = 0$.
 
\begin{table}
  \begin{center}
  \caption{Fields belonging to factors of $\textsf{P}_n(x)$} \label{tab:2}
\hspace*{-1.3cm}
\begin{tabular}{| c | c | c |}
\hline
  &   &    \\
$n$   &   factors of $\textsf{P}_n$   &   roots, field generated  \\
  &  & $\Sigma$ = Hilbert class field \\
\hline
&  &  \\
 2 &  $x^4 + 3x^3 - 9x + 9$ & $s(w/6), \Sigma$ for $\mathbb{Q}(\sqrt{-15})$   \\
  &  &  $w = \frac{3+\sqrt{-15}}{2}$\\
   &  &   \\
3 & $x^6 + 6x^5 + 2x^4 - 25x^3 + 6x^2 - 9x + 27$ & $s(w/6+1/2), \Sigma$ for $\mathbb{Q}(\sqrt{-23})$ \\
 &  &  $w = \frac{13+\sqrt{-23}}{2}$\\
 &  &   \\
 3 & $x^6 + x^5 + 2x^4 + 25x^3 + 6x^2 - 54x + 27$ & $s(w/6), \Sigma$ for $\mathbb{Q}(\sqrt{-23})$ \\
 3 & sextics map to each other under $x \rightarrow t(x)$ & \\
 &  &   \\
4  &  $x^8 + 15x^7 + 108x^6 + 270x^5 + 270x^4 - 405x^3 $ & $s(w/2)$, $\frac{-3s(w/2)}{(1-2\omega^i c(w/2))^2}$\\
 & $- 972x^2 + 729$ &  $\Omega_3$ for $\mathbb{Q}(\sqrt{-7})$  \\
 &   &   \\
 4 & $x^8 - 12x^6 + 15x^5 + 30x^4 - 90x^3 + 108x^2 $ & $3r(w/2)$, $\Omega_3$ for $\mathbb{Q}(\sqrt{-7})$  \\
 & $- 45x + 9$  &  $w = \frac{1+\sqrt{-7}}{2}$ \\
 &  &    \\
 4 & $x^8 - 3x^7 - 9x^6 + 99x^5 + 144x^4 - 297x^3 - 81x^2 $ & $s(w/6), \Sigma$ for $\mathbb{Q}(\sqrt{-39})$  \\
  & $+ 81x + 81$  & $w = \frac{3+\sqrt{-39}}{2}$ \\
 4 & all octics map to themselves under $x \rightarrow t(x)$ & \\
 &  &   \\
 5 & $x^{10} - 6x^9 - 30x^8 + 147x^7 + 446x^6 - 340x^5 - 807x^4$ & $s(w/6 + 1/2), \Sigma$ for $\mathbb{Q}(\sqrt{-47})$ \\
 & $ + 27x^3 + 675x^2 - 324x + 243$ & $w = \frac{1+\sqrt{-47}}{2}$ \\
 5 & $x^{10} + 4x^9 + 25x^8 - 3x^7 - 269x^6 + 340x^5+ 1338x^4$ & $s(w/6), \Sigma$ for $\mathbb{Q}(\sqrt{-47})$\\
 &  $ - 1323x^3 - 810x^2 + 486x + 243$ & \\
 5 & $10$-th degrees map to each other under $x \rightarrow t(x)$ & \\
 &  &  \\
 5 & $x^{20} + 36x^{19} + 376x^{18} + 532x^{17} - 1604x^{16} - 7195x^{15}$ & $s(w/6+1/2)$, \\
 &  $- 34456x^{14} + 48433x^{13} + 320874x^{12} - 635035x^{11} $ & $\Sigma$ for $\mathbb{Q}(\sqrt{-7 \cdot 17})$\\
 & $ - 1279627x^{10} + 4012587x^9 + 3393657x^8 - 10257489x^7 $ & $w = \frac{1+\sqrt{-7\cdot 17}}{2}$\\
 &  $- 1980126x^6 + 10119735x^5 - 1482057x^4 $ & \\
  &  $- 3864429x^3 + 1528713x^2 + 59049x + 59049$ & \\
 &  &  \\
 5 & $x^{20} - 3x^{19} + 233x^{18} + 1767x^{17} - 2033x^{16} - 41645x^{15}$ & $s(w/6)$, \\
 & $- 24446x^{14} + 379907x^{13} + 377073x^{12} - 1337529x^{11}$ & $\Sigma$ for $\mathbb{Q}(\sqrt{-7 \cdot 17})$ \\
 & $ - 1279627x^{10} + 1905105x^9 + 2887866x^8 - 1307691x^7 $ & \\
 & $ - 2790936x^6 + 1748385x^5 - 1169316x^4 - 1163484x^3 $ & \\
 & $+ 2466936x^2 - 708588x + 59049$ & \\
 5 & $20$-th degrees map to each other under $x \rightarrow t(x)$ & \\
 &  &. \\
 \hline
\end{tabular}
\end{center}
\end{table}

\begin{table}
  \begin{center} 
  \caption{Fields belonging to factors of $\textsf{P}_n(x)$} \label{tab:3}
\hspace*{-1.3cm}
\begin{tabular}{| c | c | c |}
\hline
  &   &    \\
$n$   &   factors of $\textsf{P}_n$   &   roots, field generated  \\
  &  &  \\
\hline
&  &  \\
 6 &  $x^{12} - 7x^{11} + 3x^{10} + 105x^9 - 286x^8 - 182x^7 + 1535x^6 $ & $-3r(w/2)$, $\Omega_3$ for $\mathbb{Q}(\sqrt{-23})$   \\
  & $ - 1736x^5+ 604x^4- 175x^3 + 222x^2 - 21x + 1$ & $w = \frac{13 + \sqrt{-23}}{2}$ \\
  &  &  \\
6 & $x^{12} - 18x^{11} + 135x^{10} + 189x^9 - 3159x^8 + 729x^7 $ & $s(w/6), \Sigma$ for $\mathbb{Q}(\sqrt{-3 \cdot 29})$ \\
& $ + 21870x^6 - 2187x^5 - 28431x^4 - 5103x^3 $ & $w = \frac{3+\sqrt{-87}}{2}$\\
 &  $+ 10935x^2 + 4374x + 729$ &  \\
 &  &   \\
6 & $x^{12} + 63x^{11} + 1998x^{10} + 4725x^9 + 48924x^8 + 421848x^7$ & $s(w/2)$, $\Omega_3$ for $\mathbb{Q}(\sqrt{-23})$ \\
& $+ 1119015x^6 + 398034x^5 - 1876446x^4 - 2066715x^3 $ & $w = \frac{13+\sqrt{-23}}{2}$\\
& $ + 177147x^2 + 1240029x + 531441$ & \\
 &  &   \\
6  &  $x^{12} - 18x^{11} + 135x^{10} + 3834x^9 + 18711x^8 + 22599x^7 $ & $s(w/2)$, $\Omega_3$ for $\mathbb{Q}(\sqrt{-15})$\\
 & $ - 14580x^6- 2187x^5 + 37179x^4 - 85293x^3 - 98415x^2 $ & $w = \frac{23+\sqrt{-15}}{2}$   \\
  & $+ 59049x + 59049$ & \\
 &   &   \\
6 & $x^{12} - 3x^{11} - 15x^{10} + 39x^9 + 51x^8 + 9x^7 - 180x^6 $ & $-3r(w/2)$, $\Omega_3$ for $\mathbb{Q}(\sqrt{-15})$  \\
& $ - 837x^5 + 2079x^4- 1278x^3 + 135x^2 + 54x + 9$ & \\
6 & $12$-th degrees map to themselves under $x \rightarrow t(x)$ & \\
 &  &    \\
6 & $x^{24} + 99x^{23} + 4914x^{22} + 31023x^{21} + 147798x^{20} $ & $s(w/6), \Sigma$ for $\mathbb{Q}(\sqrt{-255})$ \\
& $ - 815913x^{19} - 5654853x^{18} + 6884919x^{17} $ & $w = \frac{3+\sqrt{-255}}{2}$ \\
& $  + 66944313x^{16} - 10292022x^{15} - 445515957x^{14} $ & \\
& $  - 117358794x^{13} + 1841605632x^{12} + 352076382x^{11} $ & \\
& $ - 4009643613x^{10} + 277884594x^9 + 5422489353x^8 $ & \\
& $ - 1673035317x^7- 4122387837x^6 $ & \\
 & $+ 1784401731x^5 + 969702678x^4 - 610625709x^3$ &  \\
  &  $+ 290166786x^2 - 17537553x + 531441$ & \\
 &  &   \\
6 & $x^{24} - 63x^{23} + 1917x^{22} + 34587x^{21} + 73035x^{20}$ &  $s(w/6), \Sigma$ for $\mathbb{Q}(\sqrt{-231})$\\
& $  - 928179x^{19} - 2353698x^{18} + 7754859x^{17} $ & $w = \frac{3+\sqrt{-231}}{2}$\\
& $ + 18825453x^{16} - 9797760x^{15} - 49760811x^{14}$ &  \\
& $ - 127064700x^{13} + 87525198x^{12} + 381194100x^{11}$ & \\
& $ - 447847299x^{10} + 264539520x^9 + 1524861693x^8$  & \\
& $- 1884430737x^7 - 1715845842x^6 + 2029927473x^5 $ & \\
& $+ 479182635x^4 - 680775921x^3$ & \\
& $ + 113196933x^2 + 11160261x + 531441$ &  \\
6 & $24$-th degrees map to themselves under $x \rightarrow t(x)$ & \\
 &  &. \\
 \hline
\end{tabular}
\end{center}
\end{table}

\begin{table}[htp]
\caption{Minimal polynomials for periodic points of $g_1(x,y) = 0$.}  \label{tab:4}
\begin{center}
\hspace*{-1.3cm}
\begin{tabular}{|c|c|c|c|}
\hline
&  &  \\
$n$ & Factors of $\mathcal{P}_n(x)$ &roots, field generated\\
&  & $\Sigma$ = Hilbert class field \\
\hline
&  &  \\
$2$ & $x^2+x+1$ & $\omega$\\
&  &  \\
$2$ & $x^2-2x+4$ & $-2\omega$\\
&  &  \\
$3$ & $x^6+5x^4-3x^3+12x^2+4x+8$ & $2\omega c(w/6)$, \\
&  &   $\Sigma$ for $\mathbb{Q}(\sqrt{-23})$\\
&  &  \\
$3$ & $x^6-x^5+6x^4+3x^3+10x^2+8$ & $2\omega^2 c(w/6+1/2)$, \\
&  &   $\Sigma$ for $\mathbb{Q}(\sqrt{-23})$\\
&  &  \\
$4$ & $x^8+x^7+16x^6+16x^5+7x^4-32x^3+64x^2-8x+16$ & $2\omega c(w/2)$,\\
&  &   $\Omega_3$ for $\mathbb{Q}(\sqrt{-7})$\\
&  &  \\
$4$ & $x^8-5x^7+7x^6+7x^5+37x^4+58x^3-8x^2-32x+16$ & $2\omega^2 c(w/2)$,\\
&  &   $\Omega_3$ for $\mathbb{Q}(\sqrt{-7})$\\
&  &  \\
$4$ & $x^8+4x^7-2x^6-29x^5+37x^4-14x^3+28x^2+40x+16$ &  $2c(w/2)$\\
&  &   $\Omega_3$ for $\mathbb{Q}(\sqrt{-7})$\\
&  &  \\
$5$ & $x^{10}-4x^9+13x^8-2x^7+14x^6+35x^5+98x^4$ & $2\omega^2 c(w/6+1/2)$,\\
& $+48x^3+24x^2-16x+32$ & $\Sigma$ for $\mathbb{Q}(\sqrt{-47})$\\
&  &  \\
$5$ & $x^{10}+x^9+3x^8-12x^7+49x^6-35x^5+28x^4$ & $2\omega c(w/6)$,\\
& $+8x^3+104x^2+64x+32$ & $\Sigma$ for $\mathbb{Q}(\sqrt{-47})$\\
&  &  \\
$5$ & $x^{20}+9x^{19}+50x^{18}+48x^{17}+218x^{16}+698x^{15}+386x^{14}+463x^{13}$ & $2\omega c(w/6)$\\
& $+4370x^{12}-6295x^{11}+10855x^{10}-19578x^9+22392x^8-22808x^7$ & $\Sigma$ for $\mathbb{Q}(\sqrt{-7\cdot 17})$\\
& $+32560x^6-7520x^5+27200x^4+640x^3+10752x^2+3584x+1024$ & \\
&  &  \\
$5$ & $x^{20}-7x^{19}+42x^{18}-5x^{17}+425x^{16}+235x^{15}+2035x^{14}+2851x^{13}$ & $2\omega^2 c(w/6+1/2)$,\\
& $+5598x^{12}+9789x^{11}+10855x^{10}+12590x^9+17480x^8-3704x^7$ & $\Sigma$ for $\mathbb{Q}(\sqrt{-7\cdot 17})$\\
& $+6176x^6-22336x^5+13952x^4-6144x^3+12800x^2-4608x+1024$ & \\
&  &  \\
\hline
\end{tabular}
\end{center}
\end{table}

\begin{table}[htp]
\caption{Minimal polynomials for periodic points of $g_1(x,y) = 0$, period $6$.}  \label{tab:5}
\begin{center}
\hspace*{-1.3cm}
\begin{tabular}{|c|c|c|c|}
\hline
&  &  \\
$n$ & Factors of $\mathcal{P}_6(x)$ &roots, field generated\\
&  &  \\
\hline
&  &  \\
$6$ & $x^{12} + x^{11} - 5x^{10} + 12x^9 + 29x^8 - 38x^7 - 35x^6 - 22x^5$ & $2c(w/6 + 1/2)$\\
&  $+ 52x^4 + 48x^3 - 80x^2 + 64$ & $\Omega_{3}$ for $\mathbb{Q}(\sqrt{-23})$ \\
&  &  \\
$6$ & $x^{12} + 15x^{11} + 30x^{10} - 688x^9 + 1737x^8 - 2061x^7 + 3381x^6$ & $2\omega c(w/2)$\\
& $+ 4500x^5 + 7920x^4 + 1616x^3 + 2496x^2 + 672x + 64$ & $\Omega_{3}$ for $\mathbb{Q}(\sqrt{-23})$\\
&  &  \\
$6$ & $x^{12} - 21x^{11} + 156x^{10} - 202x^9 + 1980x^8 - 2250x^7 + 3381x^6$ & $2c(w/2)$\\
& $+ 4122x^5 + 6948x^4 + 5504x^3 + 480x^2 - 480x + 64$ & $\Omega_{3}$ for $\mathbb{Q}(\sqrt{-23})$ \\
&  &  \\
$6$ & $x^{12} - 5x^{10} - 6x^9 + 13x^8 + 11x^7 - 35x^6 + 76x^5 + 116x^4$ & $2c(w/6)$ \\
& $- 96x^3 - 80x^2 - 32x + 64$ & $\Omega_{3}$ for $\mathbb{Q}(\sqrt{-23})$ \\
&  &  \\
$6$ & $x^{12} + 6x^{11} + 165x^{10} + 905x^9 + 1980x^8 - 1089x^7 - 6609x^6$ & $2\omega^2 c(w/2)$\\
& $+ 2178x^5 + 7920x^4 - 7240x^3 + 2640x^2 - 192x + 64$ & $\Omega_{3}$ for $\mathbb{Q}(\sqrt{-23})$  \\
&  &  \\
$6$ & $x^{12} + 5x^9 + 51x^6 - 40x^3 + 64$ & $2c(w/6), w = \frac{3+\sqrt{-15}}{2}$\\
&  & $\Omega_3$ for $\mathbb{Q}(\sqrt{-15})$ \\
&  &  \\
$6$ & $x^{12} - 3x^{11} + 66x^{10} - 175x^9 + 495x^8 - 198x^7 + 1383x^6$ & $2c(w/2), , w = \frac{3+\sqrt{-15}}{2}$ \\
& $ + 396x^5 + 1980x^4 + 1400x^3 + 1056x^2 + 96x + 64$ &   $\Omega_3$ for $\mathbb{Q}(\sqrt{-15})$\\
&  &  \\
$6$ & $x^{12} - 12x^{11} + 48x^{10} + 14x^9 + 90x^8 + 801x^7 + 1005x^6$ & $2c(w/2), w = \frac{23+\sqrt{-15}}{2}$\\
& $ - 306x^5 - 1584x^4 - 112x^3 + 1200x^2 - 480x + 64$ & $\Omega_3$ for $\mathbb{Q}(\sqrt{-15})$\\
&  &  \\
$6$ & $x^{12} + 15x^{11} + 75x^{10} + 14x^9 - 396x^8 + 153x^7 + 1005x^6$ & $2\omega^2 c(w/2), w = \frac{23+\sqrt{-15}}{2}$\\
& $- 1602x^5 + 360x^4 - 112x^3 + 768x^2 + 384x + 64$ & $\Omega_3$ for $\mathbb{Q}(\sqrt{-15})$\\
&  &  \\

\hline
\end{tabular}
\end{center}
\end{table}

\section{Appendix: Corrections and additions to Part I.}

\begin{thm}
Let $c_1(\tau) = c(\tau +3/2)$.  If $-d = d_K f^2 \equiv 1$ (mod $3$) and
$$w = k + \frac{\sqrt{-d}}{2} \ \ \textrm{or} \ \ \frac{k+\sqrt{-d}}{2}, \ \ 9 \mid N(w), \ k \equiv 1 \ (\textrm{mod} \ 6),$$
according as $-d$ is even or odd, then $2c(w/3), \frac{c_1(w/6)}{c(w/3)}=\frac{-1}{c(w/6)}$ and $\frac{-1}{c_1(w/6)}$ are periodic points of the algebraic function defined by
$$f(x,y) = x^3(y^2 + y + 1) - y(y^2 - 2y + 4) = 0.$$
Along with the values in the set $\{0,1,-2,\pm \sqrt{-2}\}$, these numbers and their conjugates over $\mathbb{Q}$ are the only periodic points 
of the algebraic function defined by $f(x,y) = 0$, as $-d$ ranges over negative quadratic discriminants congruent to $1$ (mod $3$).  The periodic points 
$2c(w/3), \frac{c_1(w/6)}{c(w/3)}=\frac{-1}{c(w/6)}$ and $\frac{-1}{c_1(w/6)}$ corresponding to the discriminant $-d$ all lie in the ring class field $\Omega_{2f}$.
\label{thm:7.1}
\end{thm}

\begin{proof}
All assertions but part of the last are proved in \cite{akm}, and the last assertion is partly implicit there.  The fact that $2c(w/3)$ lies in $\Omega_{2f}$ 
is proved in \cite[Thm. 1, pp. 1734-1735]{akm}, and this implies the same for the values $-1/c(w/6), -1/c_1(w/6)$, when $\left(\frac{-d}{2}\right) = 0, -1$; all three are roots 
of a cubic polynomial over $\Omega_f$, so at least one of the latter values is conjugate to $2c(w/3)$ over $\Omega_f$.  If $\left(\frac{-d}{2}\right) = +1$, 
the same arguments used in the proof of Theorem 1 in \cite{akm} apply to the value $c(w/6)$, since $w/6$ is then the ratio of a basis for a proper ideal in the order 
$\textsf{R}_{-d}$.  Then the relation $c(w/6) c_1(w/6) = -c(w/3)$ from \cite[p. 1730]{akm} shows that $c_1(w/6) \in \Omega_{2f}$, as well.
\end{proof}

\begin{lem}
If $q(x)$ is the irreducible factor of $C_d(x) = p_d\left(\frac{x^3+2}{x}\right)$ having $2c(w/3)$ as a root, and $\delta = \textrm{deg}(q(x))$, then
$$(x-1)^\delta q(B(x)) = q(1)q(x), \ \ B(x) = \frac{x+2}{x-1}.$$
Also, the map $x \rightarrow B(x)$ permutes the roots of the polynomial $C_d(x)$.
\label{lem:7.2}
\end{lem}

\begin{proof}
This is equation (6.1) in \cite[p. 1753]{akm}, which follows from the proof of Prop.18(a) on p. 1751.  That proof is (mistakenly) only given there 
for the case when $d$ is even, but it actually applies in general.  The only details that need to be changed are that
$$C_{4d}(x) = \hat q(x) \tilde q(x), \ \ \textrm{deg}(\hat q(x)) = 2h(-4d) = 2h(-d)\left(2-\left(\frac{-d}{2}\right) \right);$$
and since $2c(w/3)$ generates $\Omega_{2f}$, that it has degree $2h(-d),4h(-d)$ or $6h(-d)$ over $\mathbb{Q}$.  This coincides
with $\textrm{deg}(\hat q(x))$, while the cofactor $\tilde q(x)$ has degree $4h(-4d)$.  The rest of the details are correct.  The second assertion
follows from
\begin{equation*}
\frac{B(x)^3+2}{B(x)} = \frac{(3x^3 + 18x + 6)}{(x + 2)(x - 1)^2} = \frac{3(x^3+2)/x+18}{(x^3+2)/x-3},
\end{equation*}
where $\sigma_1(x) = (3x+18)/(x-3)$ permutes the roots of $p_d(x)$, by \cite[Prop. 3.1, p. 863]{mring}.
\end{proof}

\noindent Dept. of Mathematics, University of Missouri – Columbia \smallskip

\noindent 304 Math Sci Bldg, 810 Rollins St., Columbia, MO, 65211 \smallskip

\noindent email: akkarapakams@missouri.edu \bigskip

\noindent Dept. of Mathematical Sciences \smallskip

\noindent Indiana University at Indianapolis (IUI) \smallskip

\noindent LD 270, 402 N. Blackford St., Indianapolis, IN, 46202 \smallskip

\noindent e-mail: pmorton@iu.edu


\medskip

\begin{thebibliography}{WWW}

\bibitem[1]{akm0} S. J. Akkarapakam and P. Morton, Periodic points of algebraic functions related to a continued fraction of Ramanujan, New York J. Math. {\bf 30} (2024), 783-827.

\bibitem[2]{akm} S. J. Akkarapakam and P. Morton, On Ramanujan's cubic continued fraction, Ramanujan J. {\bf 65} (2024), 1709-1756.

\bibitem[3]{ckp} B. Cho, J.K. Koo, and Y.K. Park, On Ramanujan's continued fraction as a modular function, Tohoku Math. J. {\bf 62} (2010), 579-603.

\bibitem[4]{coo} S. Cooper, {\it Ramanujan's Theta Functions}, Springer, 2017.

\bibitem[5]{cox} D.A. Cox, {\it Primes of the form $x^2+ny^2$, Fermat, Class Field Theory, and Complex Multiplication}, 2nd ed., Wiley, 2013.

\bibitem[6]{deu2} M. Deuring, Die Klassenk\"orper der komplexen Multiplikation, Enzyklop. d. math. Wissensch. I2.2 {\bf 23}, Heft 10 \ II, 1958.

\bibitem[7]{duk} W. Duke, Continued fractions and modular functions, Bull. Amer. Math. Soc. {\bf 42}, No. 2 (2005), 137-162.

\bibitem[8]{hbeg} H. Hasse, Neue Begr\"undung der komplexen Multiplikation. I. Einordnung in die allgemeine Klassenk\"orpertheorie, J. f\"ur reine u. angew. Math. {\bf 157} (1927), 115-139.  Paper {\bf 33} in {\it Helmut Hasse Mathematische Abhandlungen, Band II}, Walter de Gruyter, Berlin, 1975.

\bibitem[9]{hzer} H. Hasse, Das Zerlegungsgesetz f\"ur die Teiler des Moduls in den Ringklassenk\"orpern der komplexen Multiplikation, Monatshefte f\"ur Mathematik u. Physik {\bf 38} (1931), 331-344.  Paper {\bf 37} in {\it Helmut Hasse Mathematische Abhandlungen, Band II}, Walter de Gruyter, Berlin, 1975.

\bibitem[10]{ka} E. Kani, Idoneal numbers and some generalizations, Ann. Sci. Math. Québec {\bf 35}(2) (2011), 197-227.

\bibitem[11]{ko} G. K\"ohler, {\it Eta Products and Theta Series Identities}, Springer-Verlag Berlin Heidelberg, 2011.

\bibitem[12]{lm} R. Lynch and P. Morton, The quartic Fermat equation in Hilbert class fields of imaginary quadratic fields, Int. J. of Number Theory {\bf 11} (2015), 1961-2017.

\bibitem[13]{mcong} P. Morton, Explicit congruences for class equations, Functiones Approx. {\bf 51} (2014), 77-110.

\bibitem[14]{mring} P. Morton, Solutions of the cubic Fermat equation in ring class fields of imaginary quadratic fields (as periodic points of a $3$-adic algebraic function), Int. J. Number Theory {\bf 12} (2016), 853-902.

\bibitem[15]{mdeu} P. Morton, Periodic points of algebraic functions and Deuring's class number formula, The Ramanujan J. {\bf 50} (2019), 323-354.

\bibitem[16]{mnew} M. Newman, Construction and application of a class of modular functions (II), Proc. London Math. Soc. {\bf 9} (1959), 373-387.

\bibitem[17]{sch} R. Schertz, {\it Complex Multiplication}, Cambridge University Press, 2010.

\end{thebibliography}
\end{document}